\documentclass[oneside,english]{amsart}
\usepackage{tgtermes}
\usepackage[T1]{fontenc}
\usepackage[latin9]{inputenc}
\usepackage{color}
\usepackage{babel}
\usepackage{mathtools}
\usepackage{amsbsy}
\usepackage{amstext}
\usepackage{amsthm}
\usepackage{amssymb}
\usepackage{graphicx}
\usepackage[unicode=true,pdfusetitle,
 bookmarks=true,bookmarksnumbered=true,bookmarksopen=true,bookmarksopenlevel=3,
 breaklinks=false,pdfborder={0 0 0},pdfborderstyle={},backref=false,colorlinks=true]
 {hyperref}

\makeatletter

\newcommand*{\LyxTextAccent}[3][0ex]{%
  \hmode@bgroup\ooalign{\null#3\crcr\hidewidth
  \raise#1\hbox{#2}\hidewidth}\egroup}
\newcommand{\LyxAccentSize}[1][\sf@size]{%
  \check@mathfonts\fontsize#1\z@\math@fontsfalse\selectfont
}
\ProvideTextCommandDefault{\textcommabelow}[1]{
  \LyxTextAccent[-.31ex]{\LyxAccentSize,}{#1}}

\numberwithin{equation}{section}
\numberwithin{figure}{section}

\newtheorem{theorem}{Theorem}[section]
\newtheorem{proposition}[theorem]{Proposition}%
\newtheorem{lemma}[theorem]{Lemma}

\newtheorem{definition}[theorem]{Definition}
\newtheorem{remark}[theorem]{Remark}%
\newtheorem{claim}[theorem]{Claim}

\makeatother

\begin{document}
\title{Relative Ding Stability and an Obstruction to the Existence of Mabuchi Solitons}
\author{Yi Yao}
\email{yeeyoe@163.com}
\address{School of Mathematics, Hunan University, Changsha, 410082, CHINA.}

\begin{abstract}
Mabuchi solitons generalize K\"{a}hler-Einstein metrics on Fano manifolds,
which constitute a Yau-Tian-Donaldson type correspondence with relative
Ding stability. Comparing with K\"{a}hler-Ricci solitons, there is a distinct
necessary condition for the existence. We show this condition can
be implied by the uniformly relative Ding stability. For this we study
the inner product of $\mathbb{C}^{*}$-actions on equivariant test-configurations
and obtain an integration formula over the total space. To analyze
the uniform stability, by adapting Okounkov body construction to the
setting of torus action, we give a convex-geometry description for
the reduced non-Archimedean J-functionals.
\end{abstract}



\maketitle

\section{Introduction}

Since the Yau-Tian-Donaldson correspondence for Fano manifolds is
established, namely the equivalence between the existence of K\"{a}hler-Einstein
(KE) metrics and the K-polystability of Fano manifolds, it is desired
to extend this correspondence to other kinds of canonical metrics.
It is well-known that a Fano manifold $M$ with nonvanishing Futaki
invariant does not admit KE metrics, we can seek other canonical metrics
such as K\"{a}hler-Ricci(KR) solitons and extremal metrics. In \cite{Mabuchi},
Mabuchi introduced a new kind of generalized KE metrics now we call
\textit{Mabuchi solitons}, since they give self-similar solutions
of the gradient flow of Ding functional, see \cite{inv MA}. Let
\[
\omega=\sqrt{-1}g_{i\bar{j}}dz^{i}\wedge d\bar{z}^{j}\in2\pi c_{1}(M)
\]
be a K\"{a}hler metric with Ricci potential $h_{\omega}$ which satisfies
$Ric(\omega)-\omega=i\partial\bar{\partial}h_{\omega}$ and $\int e^{h_{\omega}}\omega^{n}=V\coloneqq\int\omega^{n}$.
By the definition in \cite{Mabuchi}, $\omega$ is called a Mabuchi
soliton if the $(1,0)$-gradient vector field
\[
\textrm{grad}_{\omega}e^{h_{\omega}}\coloneqq\sum_{i,j}\left(g^{i\bar{j}}\frac{\partial e^{h_{\omega}}}{\partial\bar{z}^{j}}\right)\frac{\partial}{\partial z^{i}}
\]
is holomorphic. It turns out this vector field must coincide with
the extremal vector field $Z$ defined in \cite{Futaki-Mabuchi},
which can be determined in advance by Futaki invariant and a chosen
maximal compact subgroup $K$ of $\textrm{Aut}^{0}(M)$. Mabuchi solitons
can be regarded as an algebraic analogue of KR-solitons which require
$\textrm{grad}_{\omega}h_{\omega}$ is holomorphic. From the view
of PDE, fixing a background metric $\omega_{0}$, let $\theta_{0}$
be the Hamiltonian function of $Z$ w.r.t. $\omega_{0}$ which satisfies
$\iota_{Z}\omega_{0}=i\bar{\partial}\theta_{0}$ and $\int\theta_{0}\omega_{0}^{n}=0$.
Then $\omega_{u}=\omega_{0}+i\partial\bar{\partial}u$ is a Mabuchi
soliton if and only if $u$ satisfies a Monge-Amp\`{e}re type equation
\begin{equation}
\left(1-\theta_{Z}(u)\right)\omega_{u}^{n}=e^{h_{\omega_{0}}-u}\omega_{0}^{n},\label{intro-MA equation}
\end{equation}
where $\theta_{Z}(u)=\theta_{0}+Z(u)$ is the Hamiltonian function
of $Z$ w.r.t. $\omega_{u}$, satisfying $\iota_{Z}\omega_{u}=i\bar{\partial}\theta_{Z}(u)$
and $\int\theta_{Z}(u)\omega_{u}^{n}=0$. As a comparison, a KR-soliton
$\omega_{u}$ should satisfy
\[
\omega_{u}^{n}=e^{h_{\omega_{0}}-\theta_{X}(u)-u}\omega_{0}^{n},
\]
where $X$ is a holomorphic vector field. A significant difference
is that (\ref{intro-MA equation}) admits a smooth solution only if
\[
\vartheta(M)\coloneqq\max_{M}\theta_{Z}(u)<1.
\]
Actually, this maximum is an invariant of Fano manifold $M$ (see
Def. \ref{def: VarTheta M} for details, we alter the notation $\alpha_{M}$
used in \cite{Mabuchi} to avoid confusion with Tian's $\alpha$-invariants).
It concerns with the Duistermaat-Heckman measure of the extremal action,
i.e. the $\mathbb{C}^{*}$-action generated by $Z$.

Mabuchi solitons also can be regarded as an analogue of extremal metrics
from the view of a new GIT (geometric invariant theory) model recently
constructed by Donaldson \cite{Dona Ding}, in which the moment map
is Ricci potential $1-e^{h_{\omega}}$, contrasting to scalar curvature
in the old model. The stability notion in the new model is \textit{Ding
stability} (abbr. as D-stability in below) which arised from the study
of limit slopes of Ding functional by Berman \cite{Berman}. As with
extremal metrics are critical points of Calabi's energy, Mabuchi solitons
are critical points of a similar functional
\[
Din(\omega)=\frac{1}{V}\int\left(1-e^{h_{\omega}}\right)^{2}\omega^{n},\ \omega\in2\pi c_{1}(M)
\]
called \textit{Ding energy}, see \cite{Yao}. Hence from the view
of a general framework in \cite{Gabor thesis}, as with extremal metrics
correspond to relative K-stability, Mabuchi solitons should correspond
to the relative version of D-stability.

When $M$ is a toric manifold, in \cite{Yao}, relative D-stability
is defined in terms of limit slopes of the modified Ding functional
along toric geodesic rays. In this situation, (uniformly) relative
D-stable is equivalent to $\vartheta(M)\leq1$ ($\vartheta(M)<1$),
and $M$ admits a Mabuchi soliton if and only if it is uniformly relative
D-stable.

Relative D-stability is extended to general Fano manifolds in \cite{inv MA,Hisa MS}.
Let $T\subset\textrm{Aut}^{0}(M)$ be a torus complexifying real torus
$S$. Taking a maximal compact subgroup $K\subset\textrm{Aut}^{0}(M)$
containing $S$, which determines an extremal vector field $Z$. Let
\[
\mathcal{H}_{\omega}^{Z}=\{u\in C^{\infty}(M)\mid\omega_{u}>0,\ \textrm{Im}Z.u=0\}
\]
be the space of $\textrm{Im}Z$-invariant potentials. The \textit{modified
Ding functional} $D_{Z}$ is defined on $\mathcal{H}_{\omega}^{Z}$
and takes Mabuchi solitons as the critical points,
\[
D_{Z}(u)\coloneqq-E_{Z}(u)-\log\left(\frac{1}{V}\int_{M}e^{h_{\omega}-u}\omega^{n}\right),
\]
see Definition \ref{Def of modi Ding} for $E_{Z}$. Follow the work
\cite{Berman} of Berman, for a $T$-equivariant test-configurations
$(\mathcal{X},\mathcal{L})$ of $(M,-K_{M})$, consider the limit
slope of $D_{Z}$ along the associated Phong-Sturm's geodesic ray.
The result is called the \textit{relative Berman-Ding invariant}:
\[
D_{Z}^{NA}(\mathcal{X},\mathcal{L})=D^{NA}(\mathcal{X},\mathcal{L})+\left\langle \alpha,\beta_{Z}\right\rangle ,
\]
where $D^{NA}$ is Berman-Ding invariant (\ref{eq: def of Ding inv})
and $\left\langle \alpha,\beta_{Z}\right\rangle $ is the inner product
of two $\mathbb{C}^{*}$-actions on $(\mathcal{X},\mathcal{L})$ which
will be discussed in below.

Then $M$ is said to be \textit{D-semistable relative to $T$} if
\begin{equation}
D_{Z}^{NA}(\mathcal{X},\mathcal{L})\geq0\label{semistable}
\end{equation}
for all $T$-equivariant test-configurations $(\mathcal{X},\mathcal{L})$.
$M$ is said to be \textit{uniformly} D-stable relative to $T$ if
there exists $\delta>0$ such that
\begin{equation}
D_{Z}^{NA}(\mathcal{X},\mathcal{L})\geq\delta\cdot J_{T}^{NA}(\mathcal{X},\mathcal{L})\label{Unif stable}
\end{equation}
for all $T$-equivariant test-configurations, where $J_{T}^{NA}$
is the \textit{reduced non-Archimedean(NA) J-functional} introduced
by Hisamoto \cite{Hisa toric}.

In \cite{LiY ZhouB}, it shows that the existence of Mabuchi solitons
is equivalent to the properness of $D_{Z}$. Then existence can imply
the above uniform stability by taking the limit slopes. For the converse
direction, Hisamoto \cite{Hisa MS} and Han-Li \cite{Li-HanJY} shows
that uniform stability and $\vartheta(M)<1$ can imply the existence,
based on the variational approach of \cite{BBJ}. Han-Li's work deals
with more general $g$-solitons, including Mabuchi solitons and KR-solitons.

With the expectation that all obstructions to the existence of canonical
metrics could be explained as some stability conditions, our main
result is that the additional assumption $\vartheta(M)<1$ can implied
by the uniform stability condition. This is easily obtained in the
toric case \cite{Yao}. More specifically, we have

\begin{theorem}[Main Theorem, see Theorem \ref{thm: semistable to theta 1}, \ref{thm: unif stable theta 1}]\label{thm: intro theta 1}Let
$M$ be a Fano manifold, $T\subset\textrm{Aut}^{0}(M)$ be a torus.
If $M$ is D-semistable relative to $T$ in the sense of (\ref{semistable}),
then $\vartheta(M)\leq1$. Furthermore, if $M$ is uniformly D-stable
relative to $T$ in the sense of (\ref{Unif stable}), then
\[
\vartheta(M)\leq1-\delta<1.
\]
\end{theorem}

Comparing with K-stability, a new feature of relative D-stability
is that uniformly stability may be not equivalent to stability. At
least there is an orbifold toric surfaces (see Example 5.14 \cite{Yao})
which is relative D-stable but not uniformly, but the existence of
smooth such examples is not known yet.

In the following, we outline the proof of Theorem \ref{thm: intro theta 1}
and introduce its byproducts which may benefit the study of relative
K-stability.

We construct a specific family of $T$-equivariant test-configurations
$(\mathcal{X},\mathcal{L}_{c})$ with a parameter $c$, that is deformation
to the normal cone of a $T$-fixed point. Then we need to analyze
the relative Berman-Ding invariant,
\[
D_{Z}^{NA}(\mathcal{X},\mathcal{L}_{c})=D^{NA}(\mathcal{X},\mathcal{L}_{c})+\left\langle \alpha,\beta_{Z}\right\rangle _{0}.
\]
The first term (Berman-Ding invariant) is easy to compute, so we focus
on $\left\langle \alpha,\beta_{Z}\right\rangle _{0}$, which is the
inner product of the structure $\mathbb{C}^{*}$-action $\alpha$
on $(\mathcal{X},\mathcal{L}_{c})$ and the fiberwise extremal action
$\beta_{Z}$. We introduce some new methods to deal with it.

\subsubsection*{Inner products of $\mathbb{C}^{*}$-actions}

Let $(M,L)$ be a polarized manifold with a lifted action $\beta:\mathbb{C}^{*}\rightarrow\textrm{Aut}(M,L)$
generated by vector field $X$. Let $(\mathcal{X},\mathcal{L})$ be
a $\mathbb{C}^{*}$-equivariant test-configuration for $(M,L)$ with
structure action $\alpha$ and fiberwise action $\beta$. In the following,
$\mathcal{X}$ always means the canonically compactified family over
$\mathbb{P}^{1}$. The inner product $\left\langle \alpha,\beta\right\rangle _{0}$
is defined by Sz\'{e}kelyhidi \cite{Gabor thesis} via the induced $\mathbb{C}^{*}$-actions
on $\textrm{H}^{0}(\mathcal{X}_{0},k\mathcal{L}_{0})$, it requires
$\mathcal{L}$ is relatively ample, see (\ref{eq: Gabor's def}).
In Section \ref{subsec:intersec def of inner prod}, we give a new
definition, denoted by $\left\langle \alpha,\beta\right\rangle $
see (\ref{eq:def inner product}), in terms of intersection numbers
over the associated bundles for action $\beta$. Naturally, $\left\langle \alpha,\beta\right\rangle $
is invariant under pulling-back of test-configurations and valid for
general $\mathcal{L}$. It coincides with $\left\langle \alpha,\beta\right\rangle _{0}$
when $\mathcal{L}$ is relatively ample, see Theorem \ref{thm: coincide of inner product}.

This new definition gives rise to a new proof for the limit slope
formula for the modified terms of energy functionals, e.g. modified
K-energy (\ref{eq: modif K and Ding}) for extremal metrics and modified
Ding functional (\ref{eq: modif Ding}) for Mabuchi solitons.

\begin{theorem}[Limit slopes of the modified terms, see Theorem \ref{thm: limit slope}]Let
$(\mathcal{X},\mathcal{L})$ be a $\mathbb{C}^{*}$-equivariant test-configuration
for $(M,L)$, where $\mathcal{L}$ is not necessarily ample. Let $\Phi$
be a $\alpha(\mathbb{S}^{1})\times\beta(\mathbb{S}^{1})$-invariant
metric on $\mathcal{L}$ satisfying condition \textbf{A} or \textbf{B}
in Definition \ref{def: type of metric}. It induces a ray of metrics
$\{u_{t}\}_{t\geq0}$ on $L$ with curvature form $\omega_{u_{t}}$.
Suppose the fiberwise action $\beta$ is generated by vector field
$X$ on $M$, which has normalized Hamiltonian function $\theta_{X}(u_{t})$
with respect to $\omega_{u_{t}}$. Then we have
\begin{equation}
\lim_{t\rightarrow+\infty}\frac{1}{V}\int_{M}\dot{u}_{t}\theta_{X}(u_{t})\omega_{u_{t}}^{n}=\left\langle \alpha,\beta\right\rangle .\label{eq: intro lim slope formu}
\end{equation}
\end{theorem}

Case \textbf{A} includes any smooth metric $\Phi$, not necessarily
having positive curvature (even along fibers). When $\mathcal{L}$
is ample and $\{u_{t}\}$ is Phong-Sturm's geodesic ray (contained
in case \textbf{B}), this have been obtained by Hisamoto \cite{Hisa ortho proj}
through approximation by Bergman geodesics.

Our method uses equivariant Hirzebruch-Riemann-Roch (HRR) formula,
so it is valid for general rays of metrics. It is similar to the method
used by Donaldson \cite{Dona lower} (Proposition 3) when dealing
with the norm $\left\Vert (\mathcal{X},\mathcal{L})\right\Vert _{2}$.
Firstly, by Stokes formula, we convert the limit in (\ref{eq: intro lim slope formu})
to an integral of Hamiltonian function of $X$ over the total space
$\mathcal{X}$. Then we apply equivariant HRR formula on $\mathcal{X}$,
it relates the integral to the leading coefficient of equivariant
Euler characteristic $\chi_{1}^{\beta}(\mathcal{X},k\mathcal{L})$
(Definition \ref{def: Euler num d1}). Finally, via Leray's spectral
sequences, $\chi_{1}^{\beta}$ can be further related to $\left\langle \alpha,\beta\right\rangle $.
There is a little trouble is that we need smooth total space $\mathcal{X}$
to apply equivariant HRR, but thanks to the pulling-back invariance
of $\left\langle \alpha,\beta\right\rangle $, we can apply HRR on
a resolution of $\mathcal{X}$ and then come back.

As a byproduct, we obtain an integral formula for inner products.

\begin{theorem}[Integral formula for inner products, see Theorem \ref{thm: limit slope} (2)]Assume
the total space $\mathcal{X}$ is smooth. We also denote by $X$ the
generating vector field on $\mathcal{X}$ for the fiberwise action
$\beta$. Given a smooth $2$-form $\Omega\in2\pi c_{1}(\mathcal{L})$
and function $\Theta$ on $\mathcal{X}$ satisfying $\iota_{X}\Omega=i\overline{\partial}\Theta$
and $\int_{\mathcal{X}_{1}}\Theta\Omega^{n}=0$ (equivalent to integrate
along any smooth fiber), then we have
\begin{equation}
\left\langle \alpha,\beta\right\rangle =\frac{1}{(n+1)L^{n}}\int_{\mathcal{X}}\Theta\left(\frac{\Omega}{2\pi}\right)^{n+1}=\mathrm{const}\cdot\int_{\mathcal{X}}\left(\Theta+\Omega\right)^{n+2}.\label{eq: intro integral formu}
\end{equation}
\end{theorem}

Note the middle integral only involves the Hamiltonian function for
action $\beta$. The rightmost formula hints us to compute $\left\langle \alpha,\beta\right\rangle $
by localization method on $\mathcal{X}$.

Go back to the proof of Theorem \ref{thm: intro theta 1}. Firstly,
by (\ref{eq: intro integral formu}) we express $\left\langle \alpha,\beta_{Z}\right\rangle $
as an integral over $\mathcal{X}$, then expand it in terms of $c$.
After evaluating the integrals in the coefficients (by a localization
argument), we found
\begin{equation}
D_{Z}^{NA}(\mathcal{X},\mathcal{L}_{c})=\frac{1-\vartheta(M)}{(n+1)c_{1}(M)^{n}}c^{n+1}+Ac^{n+2},\ 0<c\ll1,\label{eq: Intro expand of DNA}
\end{equation}
where $A$ is a constant independent of $c$. By this expansion, relative
D-semistability immediately implies $\vartheta(M)\leq1$. But for
the case of uniform stability, we need to expand $J_{T}^{NA}(\mathcal{X},\mathcal{L}_{c})$
in terms of $c$. This leads us to study the reduced non-Archimedean
J-functionals.

\subsubsection*{A convex-geometry description for the reduced non-Archimedean J-functionals}

Let $(\mathcal{X},\mathcal{L})$ be a $T$-equivariant test-configuration,
the reduced non-Archimedean(NA) J-functional $J_{T}^{NA}(\mathcal{X},\mathcal{L})$
was introduced by Hisamoto \cite{Hisa toric}. It is defined as the
infimum of values of NA J-functional $J^{NA}$ on the all twistings
of $(\mathcal{X},\mathcal{L})$ by $1$-parameter subgroups $\rho:\mathbb{C}^{*}\rightarrow T$,
see Definition \ref{def: twist test-config} for the twisting operation.

By \cite{Nystrom TC-body}, we know $(\mathcal{X},\mathcal{L})$ induces
a filtration $\mathcal{F}(\mathcal{X},\mathcal{L})$ on the section
ring $R(M,L)$, see Section \ref{subsec: Filtra of a TC}. It is more
natural to twist $\mathcal{F}(\mathcal{X},\mathcal{L})$, since this
also makes sense for irrational $\rho$. We denotes by $\mathcal{F}(\mathcal{X},\mathcal{L})^{\rho}$
(\ref{eq: twist filtra formula}) the twisted filtration, then $J_{T}^{NA}$
is defined by
\[
J_{T}^{NA}(\mathcal{X},\mathcal{L})\coloneqq\inf_{\rho\in\mathbb{R}^{m}}J^{NA}\left(\mathcal{F}(\mathcal{X},\mathcal{L})^{\rho}\right),
\]
see (\ref{eq: JNA on filtration}) for evaluating $J^{NA}$ on filtrations.
When $(\mathcal{X},\mathcal{L})$ dominates the trivial product, $J^{NA}$
can be expressed as intersection numbers over $\mathcal{X}$, see
(\ref{eq: J intersection}). But when we twist $(\mathcal{X},\mathcal{L})$,
it no longer dominates the product. Thus we need another method to
analyze $J_{T}^{NA}$.

\begin{figure}
\centering \includegraphics[clip,scale=0.5,viewport=250bp 280bp 690bp 630bp]{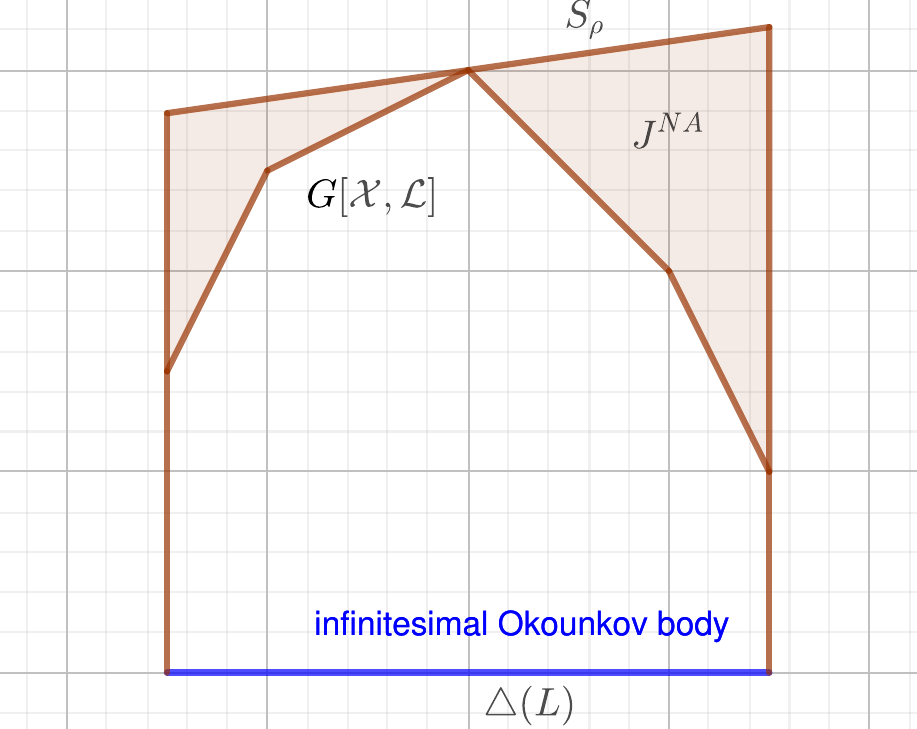}
\caption{\label{fig: NA J} The shaded area is $J^{NA}$}
\end{figure}

When $(M,L)$ is the toric manifold associated to Delzant polytope
$P$, and $(\mathcal{X},\mathcal{L})$ is a toric test-configuration
associated to a piecewise linear concave function $g$, Hisamoto \cite{Hisa toric}
gives an explicit formula
\[
J_{T}^{NA}(\mathcal{X},\mathcal{L})=\inf_{\ell}\frac{1}{\left\vert P\right\vert }\int_{P}\left(\sup_{P}(g+\ell)-(g+\ell)\right)dy,
\]
where $\ell$ runs over all affine functions. It is same to the minimum
volume of the region bounded by the graph of $g$ and its support
function, see the shadow area in Fig. \ref{fig: NA J}.

We extend this picture to general setting by replacing $P$ with the
Okounkov body $\triangle(L)$ \cite{Okounkov}. By \cite{Chen huayi},
filtration $\mathcal{F}(\mathcal{X},\mathcal{L})$ induces a concave
function $G[\mathcal{X},\mathcal{L}]$ on $\triangle(L)$ called the
\textit{concave transform} of $\mathcal{F}(\mathcal{X},\mathcal{L})$,
which takes over the role of $g$ in the toric setting. But for adapting
with the $T$-action, we employ the \textit{infinitesimal} Okounkov
body introduced by Lazarsfeld-Musta\c{t}\u{a} \cite{L-Mustata},
which is associated to a flag of subspaces of the tangent space at
a $T$-fixed point. This construction gives rise to a $\mathbb{R}^{m}$-valued
affine function $H_{\beta}$ on $\triangle(L)$, decided by $T$-action
around the fixed point, see (\ref{eq: weight and valu}).

\begin{theorem}[A convex-geometry description for $J_{T}^{NA}$] \label{Intro convex-geometry of JNA}Let
$T\subset\textrm{Aut}(M,L)$ be a torus, and $(\mathcal{X},\mathcal{L})$
a $T$-equivariant ample test-configuration for $(M,L)$. Let $\triangle(L)$
be the infinitesimal Okounkov body (Definition \ref{Def: infi Okounkov body})
associated to a $T$-fixed point and an infinitesimal flag. Let $G=G[\mathcal{X},\mathcal{L}]$
be the concave transform (see Theorem \ref{thm: push-out measure okounkov})
of the associated filtration $\mathcal{F}(\mathcal{X},\mathcal{L})$
on the section ring.

(1) For any $\rho\in\mathbb{R}^{m}$, the concave transform of the
twisted filtration $\mathcal{F}(\mathcal{X},\mathcal{L})^{\rho}$
(\ref{eq: twist filtra formula}) is $G+\left\langle \rho,H_{\beta}\right\rangle $.
We have
\begin{equation}
J^{NA}\left(\mathcal{F}(\mathcal{X},\mathcal{L})^{\rho}\right)=\sup_{\triangle(L)}\left(G+\left\langle \rho,H_{\beta}\right\rangle \right)-\int_{\triangle(L)}\left(G+\left\langle \rho,H_{\beta}\right\rangle \right)\frac{dx}{\left\vert \triangle(L)\right\vert },\label{eq: Intro JNA formula}
\end{equation}
and $J_{T}^{NA}(\mathcal{X},\mathcal{L})$ is the infimum of above
quantities over $\rho\in\mathbb{R}^{m}$.

(2) Let $c$ be the constant making the affine function $S_{\rho}\coloneqq c-\left\langle \rho,H_{\beta}\right\rangle $
to be a support function of $G$, i.e. $\inf_{\triangle}(S_{\rho}-G)=0$,
then
\[
J^{NA}\left(\mathcal{F}(\mathcal{X},\mathcal{L})^{\rho}\right)=\int_{\triangle(L)}\left(S_{\rho}-G\right)\frac{dx}{\left\vert \triangle(L)\right\vert },
\]
which is the shadow area in Fig. \ref{fig: NA J}. In particular,
$J_{T}^{NA}(\mathcal{X},\mathcal{L})$ is bounded below by the infimum
of above quantities when $S_{\rho}$ running over all support functions
of $G$. \end{theorem}

Go back to the proof of Theorem \ref{thm: intro theta 1}, with above
convex-geometry description of $J_{T}^{NA}$, we can show
\[
J_{T}^{NA}(\mathcal{X},\mathcal{L}_{c})=\frac{c^{n+1}}{(n+1)L^{n}},\ \textrm{when}\ c\ll1.
\]
Combing with (\ref{eq: Intro expand of DNA}), we see uniformly relative
D-stability implies $\vartheta(M)<1$. The proof of main theorem is
completed.

\subsubsection*{From Monge-Amp\`{e}re measures to Duistermaat-Heckman measures}

Let $\{u_{t}\}$ be the Phong-Sturm's geodesic ray associated to a
test-configuration $(\mathcal{X},\mathcal{L})$, Hisamoto \cite{Hisa measure}
showed the pushforward of Monge-Amp\`{e}re measure $V^{-1}\omega_{u_{t}}^{n}$
by $\dot{u}_{t}$ is equal to the Duistermaat-Heckman measure $\textrm{DH}(\mathcal{X},\mathcal{L})$,
see Section \ref{subsec:DH measure} for definition. In particular,
these pushforward measures are independent of $t$, but this no longer
holds for general rays. By similar techniques with the proof of (\ref{eq: intro lim slope formu}),
we extend Hisamoto's result to a more wider class of rays.

\begin{theorem} \label{thm: lim to DH them}Let $(\mathcal{X},\mathcal{L})$
be an ample test-configuration for $(M,L)$. Let $\Phi$ be a $\alpha(\mathbb{S}^{1})$-invariant
metric on $\mathcal{L}$ satisfying condition \textbf{A} or \textbf{B}
in Definition \ref{def: type of metric}. It induces a ray of metrics
$\{\phi_{t}=\psi_{0}+u_{t}\}_{t\geq0}$ on $L$. We assume $\omega_{u_{t}}\geq0$
in case \textbf{A}. Then the pushforward measure
\[
\left(\dot{u}_{t}\right)_{\#}\left(\frac{1}{V}\omega_{u_{t}}^{n}\right)\rightarrow\textrm{DH}(\mathcal{X},\mathcal{L})
\]
weakly when $t\rightarrow+\infty$. \end{theorem}

\subsubsection*{K\"{a}hler-Ricci solitons, Mabuchi solitons and extremal metrics}

For Fano manifolds, there are three kinds of canonical metrics in
$c_{1}(M)$ listed in the title (abbr. KRS, MS and extK). KRS and
MS are studied as instances of more general $g$-solitons in \cite{Li-HanJY},
where contains a complete discussion of existence, stability and properness
of energy functionals. For extK (in general class), the equivalence
of existence with properness of modified K-energy has been established
in \cite{Chen-Cheng}, also see \cite{HeWY}. From properness to uniform
stability can be obtained by the limit slope formulas in \cite{Chi Li cscK,Li-HanJY}.
The reverse direction has not been fully achieved, see \cite{Chi Li cscK}
for recent progress via BBJ's \cite{BBJ} approach.

It seems that KRS has no direct connection to MS except they satisfy
similar equations. Toric manifolds always admit KRS \cite{Wang Zhu}
but not for MS \cite{Yao}. Although MS and extK satisfy very different
equation (second v.s. fourth order), but the existence of MS can imply
the existence of extK due to inequality: $M_{Z}\geq D_{Z}-C$, see
Remark \ref{rem extK and MS}.

\subsubsection*{Some of related works}

For a $\mathbb{C}^{*}$-equivariant test-configuration $(\mathcal{X},\mathcal{L})$
with smooth $\mathcal{X}$, let $\Theta_{W}$ and $\tilde{\Theta}_{X}$
be the Hamiltonian function on $\mathcal{X}$ for action $\alpha$
and $\beta$ respectively. In \cite{Dervan relative}, $\left\langle \alpha,\beta\right\rangle _{0}$
is defined to be $\int_{\mathcal{X}_{0}}\Theta_{W}\tilde{\Theta}_{X}\Omega^{n}$,
then (\ref{eq: intro lim slope formu}) is obtained easily. Our method
avoids to integrate along $\mathcal{X}_{0}$, since it may be singular
and nonreduced. Instead, formula (\ref{eq: intro integral formu})
integrating over $\mathcal{X}$ which can be smooth after taking a
resolution. We also can use (\ref{eq: intro integral formu}) to define
inner products of actions on test-configurations for general K\"{a}hler
manifolds in the sense of \cite{Dervan relative}. This method that
converting limit slopes to an integral over $\mathcal{X}$ had been
employed by Sj\"{o}str\"{o}m Dyrefelt \cite{Dyrefelt} to extend K-stability
to transcendental K\"{a}hler classes. Finally, the twisting of filtrations
is also introduced in the work of Li \cite{Chi Li G-uniform}.

\subsubsection*{Organization}

In Section \ref{sec:Preliminaries}, we review main notions and tools,
including test-configurations, Hamiltonian functions and equivariant
HRR formula. In Section \ref{sec:Inner prod}, we study inner products
of $\mathbb{C}^{*}$-actions. In Section \ref{sec:From-Monge-Ampere-measures},
we prove Theorem \ref{thm: lim to DH them} (this section can be skipped
if the reader only cares the main theorem). In Section \ref{sec: Mabuchi-solitons-and Stability},
we review Mabuchi solitons, energy functionals and relative D-semistability.
In Section \ref{sec: theta<=00003D1}, we show semistability implies
$\vartheta(M)\leq1$. In Section \ref{sec: Unif stability}, we review
the twisting of test-configurations (filtrations), $J_{T}^{NA}$ and
uniformly relative D-stability, and give a convex-geometry description
for $J_{T}^{NA}$ via infinitesimal Okounkov body. Finally, in Section
\ref{sec: theta<1}, we show uniform stability implies $\vartheta(M)<1$.

\subsubsection*{Conventions and notations}

Let $V$ be a vector space with a $\mathbb{C}^{*}$-action. For $\tau\in\mathbb{C}^{*}$
and $x\in V$, $\tau.x$ denotes the action. We denote by $V_{\mu}\coloneqq\{x\in V\mid\tau.x=\tau^{\mu}x\}$
the weight subspace for each weight $\mu\in\mathbb{Z}$. We use $\mu$
to denote the weights of action $\alpha$ and $\nu$ for action $\beta$.
If $V$ carries two commuting $\mathbb{C}^{*}$-action $\alpha$ and
$\beta$, then it carries a $\mathbb{C}^{*}\times\mathbb{C}^{*}$-action
denoted by $\alpha\times\beta$. We denote by $V_{\mu,\nu}$ the weight
subspace with weight $\mu$ and $\nu$ w.r.t. $\alpha$ and $\beta$
respectively.

Let $(M,L)$ be a polarized manifold, set $N_{k}=h^{0}(M,kL)=\frac{L^{n}}{n!}k^{n}+O(k^{n-1})$.
Suppose $h$ is a Hermitian metric on $L$ and $s$ is a local section
of $L$, we usually use the local function $\phi=-\log\left\vert s\right\vert _{h}^{2}$
to denote $h$. Its curvature form is $Ric(h)=i\partial\bar{\partial}\phi\in2\pi c_{1}(L)$.
We set $V=\int\left(i\partial\bar{\partial}\phi\right)^{n}=(2\pi)^{n}L^{n}$.

In this paper, a $G$-equivariant test-configuration $(\mathcal{X},\mathcal{L})$
always means the canonically \textit{compactified} family over $\mathbb{P}^{1}$.
We use $(\mathcal{X},\mathcal{L})\vert_{\mathbb{C}}$ to denote the
restricted family over $\mathbb{C}$ which is the original definition.
The structure action is denoted by $\alpha:\mathbb{C}^{*}\rightarrow\textrm{Aut}(\mathcal{X},\mathcal{L})$,
and the fiberwise action is denoted by $\beta:G\rightarrow\textrm{Aut}(\mathcal{X},\mathcal{L})$.
The fibers are denoted by $\mathcal{X}_{\tau}\coloneqq\pi^{-1}(\tau)$
and $\mathcal{L}_{\tau}\coloneqq\mathcal{L}\vert_{\mathcal{X}_{\tau}}$
for $\tau=e^{-\frac{1}{2}(t+is)}\in\mathbb{P}^{1}$.

\section{Preliminaries \label{sec:Preliminaries}}

\subsection{Test-configurations\label{def TC}}

Test-configurations abstracts the one-parameter subgroups in Hilbert-Mumford's
criterion in geometric invariant theory.

\begin{definition} Let $(M,L)$ be a polarized manifold. A \textit{test-configuration}
$(\mathcal{X},\mathcal{L})$ for $(M,L)$ is constituted of a normal
variety $\mathcal{X}$ with a $\mathbb{Q}$-line bundle $\mathcal{L}$
and

(1) a flat morphism $\pi:\mathcal{X}\rightarrow\mathbb{C}$;

(2) a $\mathbb{C}^{*}$-action $\alpha$ (called the \textit{structure
action}) on $\mathcal{X}$ such that $\pi$ is equivariant w.r.t.
the multiplication action on $\mathbb{C}$;

(3) a lifting $\alpha:\mathbb{C}^{*}\rightarrow\textrm{Aut}(\mathcal{X},\mathcal{L})$,
also called a linearization of $\alpha$;

(4) an isomorphism $(\mathcal{X}_{1},\mathcal{L}_{1})\cong(M,L)$.

We say $(\mathcal{X},\mathcal{L})$ is \textit{(semi-)ample} if $\mathcal{L}$
is relatively (semi-)ample w.r.t. $\pi$. Another test-configuration
$(\mathcal{X}',\mathcal{L}')$ is called a \textit{pullback} of $(\mathcal{X},\mathcal{L})$
if there is an equivariant birational map $f:\mathcal{X}'\rightarrow\mathcal{X}$
which is isomorphism from $\mathcal{X}'\backslash\mathcal{X}'_{0}$
to $\mathcal{X}\backslash\mathcal{X}_{0}$ and $\mathcal{L}'=f^{*}\mathcal{L}$.
Following \cite{BHJ1} Def. 6.1, if $(\mathcal{X}',\mathcal{L}')$
and $(\mathcal{X},\mathcal{L})$ have a common pullback, we say they
are \textit{equivalent} with each other.

Let $G\subset\textrm{Aut}(M,L)$ be a reductive subgroup, a \textit{$G$-equivariant}
test-configuration is a test-configuration $(\mathcal{X},\mathcal{L})$
with a lifted $G$-action $\beta:G\rightarrow\textrm{Aut}(\mathcal{X},\mathcal{L})$
(called the fiberwise action) which preserves each fiber, commutes
with $\alpha$ and coincides with $G\subset\textrm{Aut}(M,L)$ on
$(\mathcal{X}_{1},\mathcal{L}_{1})\cong(M,L)$. \end{definition}

Structure action $\alpha$ and isomorphism $(\mathcal{X}_{1},\mathcal{L}_{1})\simeq(M,L)$
induce an equivariant trivialization
\[
(\mathcal{X},\mathcal{L})\vert_{\mathbb{C}^{*}}\simeq(M\times\mathbb{C}^{*},p_{1}^{*}L),
\]
where $\mathbb{C}^{*}$ trivially acts on the first factor of $M\times\mathbb{C}^{*}$.
In particular, $\mathcal{X}\backslash\mathcal{X}_{0}$ is smooth.
We can canonically compactify $(\mathcal{X},\mathcal{L})$ by gluing
it with the product $(M\times(\mathbb{P}^{1}\backslash\{0\}),p_{1}^{*}L)$
via the above trivialization. The result is a flat family over $\mathbb{P}^{1}$
with fibers are isomorphic to $(M,L)$ except over $0$. Note the
action $\alpha$ on the fiber over $\infty$ is trivial. Since we
almost always work on the compactified family, we use $(\mathcal{X},\mathcal{L})$
to denote the \textit{compactified} test-configurations, and use notation
$(\mathcal{X},\mathcal{L})\vert_{\mathbb{C}}$ for the original family
over $\mathbb{C}$.

\begin{remark} \label{rem: action on O(-1)}If we change the lifted
action $\alpha$ on $\mathcal{L}\vert_{\mathbb{C}}$ by a character,
\[
\alpha(\tau)\rightsquigarrow\tau^{c}\cdot\alpha(\tau),\ \tau\in\mathbb{C}^{*},\ c\in\mathbb{Z},
\]
then the compactified family $\mathcal{X}$ is unchanged, but $\mathcal{L}$
will change to $\mathcal{L}_{c}\coloneqq\mathcal{L}\otimes\pi^{*}\mathcal{O}_{\mathbb{P}^{1}}(c)$.
The action $\alpha$ on $\mathcal{L}_{c}$ is the product of the original
$\mathbb{C}^{*}$-action on $\mathcal{L}$ and the $\mathbb{C}^{*}$-action
on $\mathcal{O}_{\mathbb{P}^{1}}(c)$. The action on $\mathcal{O}_{\mathbb{P}^{1}}(-1)$
is defined by
\[
\tau.(z_{0},z_{1})=(\tau^{-1}z_{0},z_{1}),\ \textrm{for}\ (z_{0},z_{1})\in\mathcal{O}_{\mathbb{P}^{1}}(-1)\vert_{[z_{0},z_{1}]},
\]
note the action on fiber over $\infty=[0,1]$ is trivial. \end{remark}

\subsection{Duistermaat-Heckman measures of test-configurations \label{subsec:DH measure}}

Let $(\mathcal{X},\mathcal{L})$ be an ample test-configuration for
$(M,L)$, then the central fiber $(\mathcal{X}_{0},\mathcal{L}_{0})$
is a polarized scheme equipped with a $\mathbb{C}^{*}$-action. For
sufficiently large and divisible $k$, consider the induced $\mathbb{C}^{*}$-action
on $\textrm{H}^{0}(\mathcal{X}_{0},k\mathcal{L}_{0})$, the weight
distribution is defined by
\begin{equation}
\mathfrak{m}_{k}\coloneqq\frac{1}{N_{k}}\sum_{\mu\in\mathbb{Z}}\dim\textrm{H}^{0}(\mathcal{X}_{0},k\mathcal{L}_{0})_{\mu}\cdot\delta_{\mu/k}.\label{eq: weight measure}
\end{equation}
By Corollary 3.4 in \cite{BHJ1}, $\{\mathfrak{m}_{k}\}$ have uniformly
bounded support and converges weakly to a probability measure on $\mathbb{R}$,
denoted by $\textrm{DH}(\mathcal{X},\mathcal{L})$, called the \textit{Duistermaat-Heckman
(DH) measure} of $(\mathcal{X},\mathcal{L})$.

In Section \ref{subsec: Filtra of a TC}, $\textrm{DH}(\mathcal{X},\mathcal{L})$
is also equal to the limit measure of the associated filtration $\mathcal{F}(\mathcal{X},\mathcal{L})$
of section ring. For a semiample test-configuration, its DH measure
is defined to be that of its ample model, see Proposition 2.17 in
\cite{BHJ1} for details.

\subsection{Non-Archimedean functionals }

By the point of view of \cite{BHJ1}, a (equivalence class of) test-configuration
$(\mathcal{X},\mathcal{L})$ gives rise to a NA metric on the Berkovich
analytification $(M^{an},L^{an})$. Many functionals in K\"{a}hler geometry
have a NA counterpart which can be interpreted as the limit slopes
along geodesic rays, see \cite{BHJ2} for details.

(1) For any test-configuration $(\mathcal{X},\mathcal{L})$ of $(M,L)$,
the \textit{non-Archimedean Monge-Amp\`{e}re energy} is defined by
\[
E^{NA}(\mathcal{X},\mathcal{L})\coloneqq\frac{\mathcal{L}^{n+1}}{(n+1)L^{n}},
\]
where $\mathcal{L}^{n+1}$ is the intersection number on $\mathcal{X}$.
When $(\mathcal{X},\mathcal{L})$ is semiample, Lemma 7.3 \cite{BHJ1}
says it equals to the barycenter of DH measure,

\begin{equation}
E^{NA}(\mathcal{X},\mathcal{L})=\int_{\mathbb{R}}\lambda\ d\textrm{DH}(\mathcal{X},\mathcal{L})\label{eq: NA MA-energy}
\end{equation}

(2) The \textit{non-Archimedean J-functional} is defined by
\begin{equation}
J^{NA}(\mathcal{X},\mathcal{L})\coloneqq\sup\textrm{supp}\textrm{DH}(\mathcal{X},\mathcal{L})-\int_{\mathbb{R}}\lambda\ d\textrm{DH}(\mathcal{X},\mathcal{L}).\label{eq: NA J}
\end{equation}
It can be expressed as intersection numbers. If $\mathcal{X}$ dominates
the trivial product $M\times\mathbb{P}^{1}$, i.e. there exists $\Pi:\mathcal{X}\rightarrow M\times\mathbb{P}^{1}$
which is equivariant w.r.t. the trivial action on the target. Then
we have
\begin{equation}
J^{NA}(\mathcal{X},\mathcal{L})=(L^{n})^{-1}\mathcal{L}\cdot(\Pi^{*}p_{1}^{*}L)^{n}-E^{NA}(\mathcal{X},\mathcal{L}).\label{eq: J intersection}
\end{equation}
(3) Suppose $M$ is a Fano manifold, $(\mathcal{X},\mathcal{L})$
is an ample test-configuration for $(M,-K_{M})$. There is a unique
$\mathbb{Q}$-divisor $B$ supported on $\mathcal{X}_{0}$ such that
\[
\mathcal{L}+K_{\mathcal{X}/\mathbb{P}^{1}}=\mathcal{O}_{\mathcal{X}}(B).
\]
Let
\[
L^{NA}(\mathcal{X},\mathcal{L})\coloneqq\textrm{lct}(\mathcal{X},-B;\mathcal{X}_{0})-1,
\]
where
\[
\textrm{lct}(\mathcal{X},-B;\mathcal{X}_{0})\coloneqq\sup\{t\in\mathbb{R}\mid\textrm{pair}\ (\mathcal{X},-B+t\mathcal{X}_{0})\ \textrm{is\ log\ canonical}\}
\]
is the log canonical threshold measured on $\mathcal{X}$. Then the
\textit{Berman-Ding invariant} (or non-Archimedean Ding functional)
is defined as
\begin{equation}
D^{NA}(\mathcal{X},\mathcal{L})\coloneqq L^{NA}(\mathcal{X},\mathcal{L})-E^{NA}(\mathcal{X},\mathcal{L}).\label{eq: def of Ding inv}
\end{equation}
By \cite{Berman}, it is the limit slope of classical Ding functional
along the geodesic ray associated to $(\mathcal{X},\mathcal{L})$.

\subsection{Hamiltonian functions}

Let $M$ be a K\"{a}hler manifold with a line bundle $L$ with Hermitian
metric $h$. Its curvature form is $\omega=-i\partial\bar{\partial}\log h$.
Let $X$ be a $(1,0)$ holomorphic vector field on $M$. A real-valued
function $f$ is called a \textit{Hamiltonian function} of $X$ w.r.t.
$\omega$ if $\iota_{X}\omega=i\bar{\partial}f$. It implies $\mathcal{L}_{X}\omega=i\partial\bar{\partial}f$.
Note that $f$ is unique up to a constant, and $\textrm{Im}X$ preserves
$\omega$. If further $\int f\omega^{n}=0$, we say $f$ is \textit{normalized}
and denote it by $\theta_{X}(\omega)$.

Suppose we have a lifted action $\gamma:\mathbb{C}^{*}\rightarrow\textrm{Aut}(M,L)$
such that $\gamma(\mathbb{S}^{1})$ preserves $h$. Let holomorphic
vector field $X$ be the generator of $\gamma$ such that
\[
\gamma(e^{-\frac{1}{2}(t+is)})=\exp\left(t\cdot\textrm{Re}X-s\cdot\textrm{Im}X\right).
\]
Note that it automatically implies $\exp(4\pi\textrm{Im}X)=id_{M}$.
We can obtain a Hamiltonian function of $X$ by taking derivative
along the lifted action. Define
\[
f(x)\coloneqq-\frac{d}{dt}\left(\log\left\vert \gamma(e^{-\frac{t}{2}}).s\right\vert _{h}^{2}\right)\vert_{t=0},\ x\in M,\ s\in L_{x}.
\]
It can be verified that $\iota_{X}\omega=i\bar{\partial}f$, and change
the lifting $\gamma$ by a character will change $f$ by a constant.
As usual, we use $\phi=-\log h$ to denote the Hermitian metric $h$,
then the above formula can be shortly written as
\begin{equation}
f\coloneqq\frac{d}{dt}\left(\gamma(e^{-\frac{t}{2}})^{*}\phi\right)\vert_{t=0}.\label{eq: simple Hamil func formula}
\end{equation}
We will call $f$ the Hamiltonian function \textit{induced} by the
lifted action $\gamma$, which is not necessarily normalized.

\subsection{Equivariant Hirzebruch-Riemann-Roch formula for line bundles }

The equivariant HRR formula is our main tool to connect the integrals
of Hamiltonian functions to the equivariant Euler characteristic.
We only need the formula for equivariant line bundles, for its general
versions, see \cite{equi R-R} for the differential forms version
and the appendix of \cite{BHJ1} for the algebraic version.

Let $M$ be a $n$-dimensional complex manifold with a Hermitian line
bundle $(L,h)$, where the curvature form $\omega=-i\partial\bar{\partial}\log h$
is not necessarily positive. There is a lifted action $\gamma:\mathbb{S}^{1}\rightarrow\textrm{Aut}(M,L)$
preserving $h$, which is generated by holomorphic vector field $X$
such that $\gamma(e^{it/2})=\exp(t\cdot\textrm{Im}X)$. This lifting
induces a Hamiltonian function $f$ for $X$ such that $\iota_{X}\omega=i\bar{\partial}f$.

Lifted action $\gamma$ induces $\mathbb{S}^{1}$-actions on cohomology
groups $\textrm{H}^{q}(M,kL)$. For $k\geq1$, we define character function for $kL$ by
\[
\chi(kL,e^{it})\coloneqq\sum_{q\geq0}(-1)^{q}\textrm{Tr}\left(e^{it}\vert\textrm{H}^{q}(kL)\right)=\sum_{q\geq0}(-1)^{q}\sum_{\lambda\in\mathbb{Z}}e^{i\lambda t}\dim\textrm{H}^{q}(kL)_{\lambda}.
\]
Then the \textit{equivariant HRR formula} says
\[
\chi(kL,e^{it})=\int_{M}\textrm{ch}(kL,h^{k},t)\cdot\textrm{td}(TM,g,t),
\]
where $g$ is any $\mathbb{S}^{1}$-invariant Hermitian metric on
$TM$ (no relation with $\omega$), $\textrm{ch}(kL,h^{k},t)\coloneqq\exp k(\frac{\omega}{2\pi}+ift)$
is the equivariant Chern character form, and $\textrm{td}(TM,g,t)$
is the equivariant Todd form whose zero degree part is $1$. Expanding
both sides in terms of variable $t$, then comparing the coefficients
of $t^{d}$ ($d\geq0$), we obtain
\begin{equation}
\sum_{q\geq0}(-1)^{q}\sum_{\lambda\in\mathbb{Z}}\frac{\lambda^{d}}{d!}\dim\textrm{H}^{q}(M,kL)_{\lambda}=\int_{M}\frac{f^{d}}{d!}\frac{(\omega/2\pi)^{n}}{n!}\cdot k^{n+d}+O(k^{n+d-1}).\label{eq: HRR formula}
\end{equation}
What we really need in this paper is this expansion formula. The left
hand side is called \textit{degree-$d$ equivariant Euler characteristic},
since when $d=0$ it is the classical Euler characteristic.

\section{Inner product of $\mathbb{C}^{*}$-actions on a test-configuration
\label{sec:Inner prod}}

\subsection{The original definition by Sz\'{e}kelyhidi}

We recall the inner product of $\mathbb{C}^{*}$-actions defined by
Sz\'{e}kelyhidi \cite{Gabor thesis} in order to define relative K-stability.
In this section, $(M,L)$ is a polarized manifold, and test-configurations
are ample.

Let $(\mathcal{X},\mathcal{L})$ be a $\mathbb{C}^{*}$-equivariant
test-configuration for $(M,L)$ with structure action $\alpha$ and
fiberwise action $\beta$. Then for sufficiently divisible $k$, action
$\alpha$ and $\beta$ induce two commutative $\mathbb{C}^{*}$-actions
on $\textrm{H}^{0}(\mathcal{X}_{0},k\mathcal{L}_{0})$ (still denoted
by $\alpha$ and $\beta$). Hence $\textrm{H}^{0}(\mathcal{X}_{0},k\mathcal{L}_{0})$
becomes a $\mathbb{C}^{*}\times\mathbb{C}^{*}$-module.

In the following, we always denote the weights for $\alpha$ by $\mu\in\mathbb{Z}$
and $\nu\in\mathbb{Z}$ for $\beta$. For example, $\textrm{H}^{0}(\mathcal{X}_{0},k\mathcal{L}_{0})_{\mu}$
means the weight-$\mu$ subspace w.r.t. action $\alpha$, and $\textrm{H}^{0}(\mathcal{X}_{0},k\mathcal{L}_{0})_{\mu,\nu}$
means the weight-$(\mu,\nu)$ subspace w.r.t. action $\alpha\times\beta$.

In \cite{Gabor thesis}, the inner product of $\alpha$ and $\beta$
is defined by
\begin{eqnarray}
\left\langle \alpha,\beta\right\rangle _{0} & \coloneqq & \lim_{k\rightarrow\infty}\Bigl[\frac{1}{k^{2}N_{k}}\sum_{\mu,\nu}\mu\nu\cdot\dim\textrm{H}^{0}(\mathcal{X}_{0},k\mathcal{L}_{0})_{\mu,\nu}\nonumber \\
 &  & -\frac{1}{k^{2}N_{k}^{2}}\sum_{\mu}\mu\dim\textrm{H}^{0}(\mathcal{X}_{0},k\mathcal{L}_{0})_{\mu}\cdot\sum_{\nu}\nu\dim\textrm{H}^{0}(\mathcal{X}_{0},k\mathcal{L}_{0})_{\nu}\Bigl],\label{eq: Gabor's def}
\end{eqnarray}
where $N_{k}=h^{0}(M,kL)=\frac{L^{n}}{n!}k^{n}+O(k^{n-1})$. We use
subscript ``$0$'' to indicate this definition using the actions
on $(\mathcal{X}_{0},\mathcal{L}_{0})$.

We explicitly express $\left\langle \alpha,\beta\right\rangle _{0}$
by the coefficients of asymptotic expansions. First by Theorem 3.1
in \cite{BHJ1}, we know that
\begin{equation}
\sum_{\mu,\nu}\mu\nu\cdot\dim\textrm{H}^{0}(\mathcal{X}_{0},k\mathcal{L}_{0})_{\mu,\nu}=\frac{c_{0}}{(n+2)!}k^{n+2}+O(k^{n+1})\label{eq: def of c_0}
\end{equation}
is a polynomial of $k\gg1$ with degree at most $n+2$. Similarly,
we have
\[
\sum_{\mu}\mu\dim\textrm{H}^{0}(\mathcal{X}_{0},k\mathcal{L}_{0})_{\mu}=\frac{a_{0}}{(n+1)!}k^{n+1}+O(k^{n}).
\]
Since the weight distribution $\mathfrak{m}_{k}$ (\ref{eq: weight measure})
weakly converges to $\textrm{DH}(\mathcal{X},\mathcal{L})$, (\ref{eq: NA MA-energy})
implies $a_{0}=\mathcal{L}^{n+1}$.

Next consider action $\beta$. Since $\pi:\mathcal{X}\rightarrow\mathbb{P}^{1}$
is flat and $\mathcal{L}$ is relatively ample, by Corollary 12.9
\cite{Hartshorne}, when $k$ is large and divisible enough, $\pi_{*}(k\mathcal{L})$
will be a vector bundle on $\mathbb{P}^{1}$ and equipped with a fiberwise
action induced by $\beta$. In particular, $\textrm{H}^{0}(\mathcal{X}_{0},k\mathcal{L}_{0})$
is isomorphic to $\textrm{H}^{0}(M,kL)$ w.r.t. action $\beta$, this
implies
\begin{equation}
\sum_{\nu}\nu\dim\textrm{H}^{0}(\mathcal{X}_{0},k\mathcal{L}_{0})_{\nu}=\sum_{\nu}\nu\dim\textrm{H}^{0}(M,kL)_{\nu}=\frac{b_{0}}{(n+1)!}k^{n+1}+O(k^{n}).\label{eq: def of b_0}
\end{equation}
By these expansions, we have
\begin{equation}
\left\langle \alpha,\beta\right\rangle _{0}=\frac{1}{L^{n}/n!}\frac{c_{0}}{(n+2)!}-\frac{1}{\left(L^{n}/n!\right)^{2}}\frac{\mathcal{L}^{n+1}}{(n+1)!}\frac{b_{0}}{(n+1)!}.\label{eq:old def inner product}
\end{equation}

\subsection{An intersection-theoretic definition for inner products \label{subsec:intersec def of inner prod}}

Following \cite{BHJ1}, the pullbacks of an equivariant test-configuration
are regarded to be equivalent to the original one, it is natural to
expect they have same inner products. Instead of checking the invariance
of $\left\langle \alpha,\beta\right\rangle _{0}$, we prefer to give
another definition in terms of intersection numbers, whose invariance
is more apparent. Moreover, the new definition is valid for general
$\mathcal{L}$. We will see that it coincides with $\left\langle \alpha,\beta\right\rangle _{0}$
when $\mathcal{L}$ is relatively ample. In the following, $(\mathcal{X},\mathcal{L})$
is a $\mathbb{C}^{*}$-equivariant test-configuration for $(M,L)$,
where $\mathcal{L}$ is not necessarily ample.

Let $\mathbb{C}^{2}\backslash\{0\}\rightarrow\mathbb{P}^{1}$ be the
tautological $\mathbb{C}^{*}$-principal bundle, we consider the fiber
bundle over $\mathbb{P}^{1}$ associated to action $\beta$,
\begin{equation}
F:\mathcal{X}_{\beta}=\left(\mathbb{C}^{2}\backslash\{0\}\right)\times_{\beta}\mathcal{X}\coloneqq\left(\mathbb{C}^{2}\backslash\{0\}\right)\times\mathcal{X}/\sim_{\beta}\rightarrow\mathbb{P}^{1},\label{eq: projection F}
\end{equation}
where the equivalence relation is $(z,x)\sim_{\beta}(\tau.z,\beta(\tau)x)$,
for $z\in\mathbb{C}^{2}\backslash\{0\}$, $x\in\mathcal{X}$ and $\tau\in\mathbb{C}^{*}$.

In the same way, $\mathcal{L}$ equipped with action $\beta$ induces
a line bundle $\mathcal{L}_{\beta}$ on $\mathcal{X}_{\beta}$. By
the asymptotic Riemann-Roch formula (see Theorem 1.1.24 \cite{Lazarsfeld}),
we have expansion of Euler characteristic,
\[
\chi\left(\mathcal{X}_{\beta},k\mathcal{L}_{\beta}\right)=\frac{\mathcal{L}_{\beta}^{n+2}}{(n+2)!}k^{n+2}+O(k^{n+1}),
\]
where $\mathcal{L}_{\beta}^{n+2}$ is intersection number on $\mathcal{X}_{\beta}$
which is of dimension $n+2$. Then our new definition is just replacing
$c_{0}$ in (\ref{eq:old def inner product}) by the above leading
coefficient.

\begin{definition}[intersection-theoretic definition for inner products]Let
$(\mathcal{X},\mathcal{L})$ be a $\mathbb{C}^{*}$-equivariant test-configuration
for $(M,L)$ with structure action $\alpha$ and fiberwise action
$\beta$. The inner product of $\alpha$ and $\beta$ is defined by
\begin{equation}
\left\langle \alpha,\beta\right\rangle =\frac{1}{L^{n}/n!}\frac{\mathcal{L}_{\beta}^{n+2}}{(n+2)!}-\frac{1}{\left(L^{n}/n!\right)^{2}}\frac{\mathcal{L}^{n+1}}{(n+1)!}\frac{b_{0}}{(n+1)!},\label{eq:def inner product}
\end{equation}
where $b_{0}$ is given by (\ref{eq: def of b_0}). \end{definition}

In Theorem \ref{thm: coincide of inner product}, we will show $\left\langle \alpha,\beta\right\rangle $
coincides with $\left\langle \alpha,\beta\right\rangle _{0}$ when
$\mathcal{L}$ is ample. Thus it extends the original definition.
By above definition, immediately we have

\begin{proposition}\label{prop: invarience of inner product}

Let $(\mathcal{X},\mathcal{L})$ be a $\mathbb{C}^{*}$-equivariant
test-configuration for $(M,L)$ and $g:(\mathcal{Y},g^{*}\mathcal{L})\rightarrow(\mathcal{X},\mathcal{L})$
be a pullback, then we have $\left\langle \alpha,\beta\right\rangle _{\mathcal{X}}=\left\langle \alpha,\beta\right\rangle _{\mathcal{Y}}$.
\end{proposition}
\begin{proof}
Denote $\mathcal{K}=g^{*}\mathcal{L}$, the fiber bundle construction
induces a birational morphism $G:\mathcal{Y}_{\beta}\rightarrow\mathcal{X}_{\beta}$
such that $\mathcal{K}_{\beta}=G^{*}\mathcal{L}_{\beta}$, thus $\mathcal{K}_{\beta}^{n+2}=\mathcal{L}_{\beta}^{n+2}$.
We also have $\mathcal{K}^{n+1}=\mathcal{L}^{n+1}$, and $b_{0}$
is unchanged.
\end{proof}

\subsection{Modified terms of energy functionals converge to inner products}

When we consider the limit slope of modified energy functionals (e.g.
modified K-energy (\ref{eq: modif K and Ding}) or Ding functional
(\ref{eq: modif Ding})) along a ray of metrics $\{\phi_{t}\}$, it
will involve the limit of modified term $\int\dot{\phi}_{t}\theta_{X}(\phi_{t})\textrm{MA}(\phi_{t})$.
When the ray is induced from a test-configuration in the way given
below, we show this integral converges to the inner product. In the
sequel, for a test-configuration $(\mathcal{X},\mathcal{L})$ of $(M,L)$,
we always identify $(\mathcal{X}_{1},\mathcal{L}_{1})\cong(M,L)$.
Restricted by our method, we only consider the following two types
of rays of metrics on $(\mathcal{X},\mathcal{L})$.

\begin{definition} \label{def: type of metric}(Admissible rays of
metrics) Let $(\mathcal{X},\mathcal{L})$ be a test-configuration
for polarized manifold $(M,L)$. We consider a $\alpha(\mathbb{S}^{1})$-invariant
metric $\Phi$ on $\mathcal{L}$ satisfying one of the following conditions:

\textbf{(A)} $\Phi$ is smooth, not necessary with positive curvature
(even along fibers);

\textbf{(B)} $\Phi$ is a locally bounded psh metric (i.e. with positive
curvature current) and locally $C^{1,1}$ on $\mathcal{X}\backslash\mathcal{X}_{0}$.
The most important example is the metric on $\mathcal{L}$ gives Phong-Sturm's
geodesic ray \cite{Phong-Sturm}, its $C^{1,1}$-regularity is established
in \cite{Chu C11}.

Take a smooth reference metric $\psi_{0}$ on $L$. Since we identify
$(\mathcal{X}_{1},\mathcal{L}_{1})\cong(M,L)$, the pulling-back
\[
\{\phi_{t}=\psi_{0}+u_{t}\coloneqq\alpha(e^{-\frac{t}{2}})^{*}\Phi\mid t\geq0\}
\]
is a ray of metrics on $L$. We call this ray of metrics is \textit{induced}
by $\Phi$. Let $\Omega=i\partial\bar{\partial}\Phi$ be the curvature
current, we call $\{\phi_{t}\}$ is a \textit{subgeodesic} ray if
$\Omega$ is a positive current, and we call $\{\phi_{t}\}$ is a
\textit{geodesic} ray if further $\Omega^{n+1}=0$ holds. \end{definition}

With the assumptions in above definition, and we further assume that
$\mathcal{X}$ is smooth. Let $\omega_{0}=i\partial\bar{\partial}\psi_{0}$,
then the curvature form of $\phi_{t}$ is $\omega_{u_{t}}=\omega_{0}+i\partial\bar{\partial}u_{t}$,
and we have
\[
\omega_{u_{t}}=\alpha(e^{-t/2})^{*}\Omega,\ \mathrm{for}\ t\geq0.
\]
Suppose the structure $\mathbb{C}^{*}$-action $\alpha$ is generated
by holomorphic vector field $W$ on $\mathcal{X}$ in the sense that
\begin{equation}
\alpha(e^{-\frac{1}{2}(t+is)})=\exp\left(t\cdot\textrm{Re}W-s\cdot\textrm{Im}W\right).\label{eq: generator W}
\end{equation}
Since action $\alpha$ has been lifted to $\mathcal{L}$, by (\ref{eq: simple Hamil func formula}),
we obtain a Hamiltonian function $\Theta_{W}$ of $W$ such that
\begin{equation}
\iota_{W}\Omega=i\bar{\partial}\Theta_{W},\ \textrm{on}\ \mathcal{X}.\label{eq: identity for Theta W}
\end{equation}
Taking the time derivative of
\[
\phi_{t}=\alpha(e^{-t/2})^{*}\Phi=\alpha(e^{-t_{0}/2})^{*}\alpha(e^{-(t-t_{0})/2})^{*}\Phi
\]
at any $t\in\mathbb{R}$, we obtain
\begin{equation}
\dot{\phi}_{t}=\dot{u}_{t}=\alpha(e^{-t/2})^{*}\Theta_{W},\ \textrm{on}\ \mathcal{X}_{1}\cong M,\ \textrm{for}\ \forall t\in\mathbb{R}.\label{eq: pullback of Theta W}
\end{equation}

Next we assume $(\mathcal{X},\mathcal{L})$ is a $\mathbb{C}^{*}$-equivariant
test-configuration with fiberwise action $\beta$, and require $\Phi$
is $\beta(\mathbb{S}^{1})$-invariant. Suppose the action $\beta$
is generated by holomorphic vector field $X$ on $\mathcal{X}$ in
the sense that
\begin{equation}
\beta(e^{-\frac{1}{2}(t+is)})=\exp\left(t\cdot\textrm{Re}X-s\cdot\textrm{Im}X\right).\label{eq: generator X}
\end{equation}
Note that $X$ is tangent to each fiber except $\mathcal{X}_{0}$.
In the same way, we obtain a Hamiltonian function $\Theta_{X}$ for
$X$ from the lifted action $\beta$, it satisfies
\begin{equation}
\iota_{X}\Omega=i\bar{\partial}\Theta_{X},\ \textrm{on}\ \mathcal{X}.\label{eq: identity for Theta X}
\end{equation}

\begin{lemma} Let $f$ be any $C^{1}$ function on $\mathbb{R}$,
then the integral $\int_{\mathcal{X}_{\tau}}f(\Theta_{X})\Omega^{n}$
is independent of $\tau\in\mathbb{P}^{1}\backslash\{0\}$. \end{lemma}
\begin{proof}
Let $I(\tau)=\int_{\mathcal{X}_{\tau}}f(\Theta_{X})\Omega^{n}$, for
$\tau\in\mathbb{P}^{1}\backslash\{0\}$. By the property of fiber
integration (over $\mathbb{P}^{1}\backslash\{0\}$) and (\ref{eq: identity for Theta X}),
we have
\[
i\bar{\partial}I(\tau)=\int_{\mathcal{X}/\mathbb{P}^{1}}f'(\Theta_{X})\iota_{X}\Omega\wedge\Omega^{n}=\frac{1}{n+1}\int_{\mathcal{X}/\mathbb{P}^{1}}f'(\Theta_{X})\iota_{X}\left(\Omega^{n+1}\right).
\]
Since $X$ is tangent to each fiber, the last integral vanishes. Consider
that $I(\tau)$ is real-valued, thus it must be a constant.
\end{proof}
Let
\begin{equation}
\tilde{\Theta}_{X}=\Theta_{X}-\frac{1}{V}\int_{\mathcal{X}_{1}}\Theta_{X}\Omega^{n},\label{eq: modified Theta X}
\end{equation}
then above lemma implies $\int_{\mathcal{X}_{\tau}}\tilde{\Theta}_{X}\Omega^{n}=0$
for each $\tau\in\mathbb{C}^{*}$. We define
\begin{equation}
\theta_{X}(u_{t})\coloneqq\alpha(e^{-t/2})^{*}\tilde{\Theta}_{X},\ \textrm{on}\ \mathcal{X}_{1}\simeq M,\label{eq: pullback of Theta X}
\end{equation}
then $\int_{M}\theta_{X}(u_{t})\omega_{u_{t}}^{n}=0$, and (\ref{eq: identity for Theta X})
is translated into $\iota_{X}\omega_{u_{t}}=i\bar{\partial}\theta_{X}(u_{t})$.
Thus $\theta_{X}(u_{t})$ is the normalized Hamiltonian function for
$X$ w.r.t. $\omega_{u_{t}}$.

\begin{proposition} \label{Prop: limit to global integral}With the
assumptions and notations as above. We assume $\mathcal{X}$ is smooth
and $\Phi$ is a metric on $\mathcal{L}$ satisfying condition (A)
or (B) in Definition \ref{def: type of metric}. Then we have
\begin{equation}
\lim_{t\rightarrow+\infty}\frac{1}{V}\int_{M}\dot{u}_{t}\theta_{X}(u_{t})\omega_{u_{t}}^{n}=\frac{1}{2\pi(n+1)V}\int_{\mathcal{X}}\tilde{\Theta}_{X}\Omega^{n+1}.\label{eq: limit to global integral}
\end{equation}
\end{proposition}
\begin{proof}
First we suppose $\Phi$ satisfies condition (A). We define a function
on $\mathbb{P}^{1}\backslash\{0\}$,
\[
F(\tau)\coloneqq\int_{\mathcal{X}_{\tau}}\Theta_{W}\cdot\tilde{\Theta}_{X}\Omega^{n}.
\]
By the assumptions, $\Theta_{W}$ and $\Omega$ are $\alpha(\mathbb{S}^{1})$-invariant,
$\alpha$ and $\beta$ commutes with each other, thus $\tilde{\Theta}_{X}$
is also $\alpha(\mathbb{S}^{1})$-invariant. This follows that $F(\tau)$
only depends on $\left\vert \tau\right\vert $. By relations (\ref{eq: pullback of Theta W})
and (\ref{eq: pullback of Theta X}), we have
\[
F(e^{-t/2})=\int_{M}\dot{u}_{t}\theta_{X}(u_{t})\omega_{u_{t}}^{n},\ \textrm{for}\ t\in\mathbb{R}.
\]
Since action $\alpha$ on $(\mathcal{X},\mathcal{L})_{\infty}$ is
trivial, thus $\Theta_{W}\vert_{\mathcal{X}_{\infty}}\equiv0$, so
$F(\infty)=0$.

In terms of fiber integration along $\pi:\mathcal{X}\vert_{\mathbb{P}^{1}\backslash\{0\}}\rightarrow\mathbb{P}^{1}\backslash\{0\}$,
we have
\[
F=\int_{\mathcal{X}/(\mathbb{P}^{1}\backslash0)}\Theta_{W}\cdot\tilde{\Theta}_{X}\Omega^{n}
\]
Next we show
\begin{equation}
i\bar{\partial}\left(F(\tau)\frac{d\tau}{\tau}\right)=\frac{1}{n+1}\int_{\mathcal{X}/\mathbb{C}^{*}}\left(\tilde{\Theta}_{X}\Omega^{n+1}\right),\ \textrm{on}\ \mathbb{C}^{*},\label{eq: projective formula}
\end{equation}
where the RHS is fiber integration over $\mathbb{C}^{*}$. Exchange
$\bar{\partial}$-operator with fiber integration, and use relations
(\ref{eq: identity for Theta W}) and (\ref{eq: identity for Theta X}),
we have
\begin{eqnarray*}
i\bar{\partial}F(\tau)\wedge d\tau & = & \left(\int_{\mathcal{X}/\mathbb{C}^{*}}\iota_{W}\Omega\wedge\tilde{\Theta}_{X}\Omega^{n}+\int_{\mathcal{X}/\mathbb{C}^{*}}\Theta_{W}\cdot\iota_{X}\Omega\wedge\Omega^{n}\right)\wedge d\tau\\
 & = & \frac{1}{n+1}\int_{\mathcal{X}/\mathbb{C}^{*}}\tilde{\Theta}_{X}\iota_{W}(\Omega^{n+1})\wedge\pi^{*}(d\tau)\\
 & = & \frac{\tau}{n+1}\int_{\mathcal{X}/\mathbb{C}^{*}}\tilde{\Theta}_{X}\Omega^{n+1}.
\end{eqnarray*}
The second integral in the first row vanishes since $X$ is tangent
to each fiber. The second row is by projection formula. The third
row uses fact $\pi_{*}(W)=-\tau\frac{\partial}{\partial\tau}$.

Let $A_{r,R}\coloneqq\{\tau\mid r\leq\left\vert \tau\right\vert \leq R\}$,
applying Stokes formula on $A_{r,R}$, we have
\[
\int_{A_{r,R}}i\bar{\partial}\left(F(\tau)\frac{d\tau}{\tau}\right)=i\int_{A_{r,R}}d\left(F(\tau)\frac{d\tau}{\tau}\right)=2\pi\left(F(r)-F(R)\right).
\]
Let $r\rightarrow0$ and $R\rightarrow+\infty$, then (\ref{eq: limit to global integral})
follows from (\ref{eq: projective formula}).

When $\Phi$ satisfies condition (B), then $\Phi$ has $C^{1,1}$-regularity
on $\mathcal{X}\backslash\mathcal{X}_{0}$. To ensure the above arguments
still work, we note: (1) the fiber integration is taken over $\mathbb{P}^{1}\backslash\{0\}$;
$\Omega$ is a $2$-form with $L^{\infty}$-coefficients on $\mathcal{X}\backslash\mathcal{X}_{0}$;
the associated functions $\Theta_{W}$, $\tilde{\Theta}_{X}$ (locally
is $X\Phi$) are $C^{0,1}$ on $\mathcal{X}\backslash\mathcal{X}_{0}$;
(2) Over the total space $\mathcal{X}$, $\Omega^{n+1}$ is interpreted
as the nonpluripolar product. So the above arguments still work.
\end{proof}
Next we use the equivariant HRR formula to connect the RHS integral
of (\ref{eq: limit to global integral}) to the equivariant Euler
characteristic.

\begin{definition} \label{def: Euler num d1}Let $(\mathcal{X},\mathcal{L})$
be a $\mathbb{C}^{*}$-equivariant test-configuration for $(M,L)$
with fiberwise action $\beta$, the degree-$1$ equivariant Euler
characteristic is defined by
\[
\chi_{1}^{\beta}(\mathcal{X},k\mathcal{L})\coloneqq\sum_{q=0}^{n+1}(-1)^{q}\sum_{\nu\in\mathbb{Z}}\nu\cdot\dim\textrm{H}^{q}(\mathcal{X},k\mathcal{L})_{\nu},
\]
where $\textrm{H}^{q}(\mathcal{X},k\mathcal{L})_{\nu}$ is the weight-$\nu$
subspace of $\textrm{H}^{q}(\mathcal{X},k\mathcal{L})$ with respect
to $\mathbb{C}^{*}$-action $\beta$. \end{definition}

\begin{theorem} \label{thm: limit slope}(1) (modified terms converge
to inner products) Let $(\mathcal{X},\mathcal{L})$ be a $\mathbb{C}^{*}$-equivariant
test-configuration for $(M,L)$ with structure action $\alpha$ and
fiberwise action $\beta$ generated by $X$. $\mathcal{X}$ is not
necessarily smooth. Let $\Phi$ be a $\alpha(\mathbb{S}^{1})\times\beta(\mathbb{S}^{1})$-invariant
metric on $\mathcal{L}$ satisfying condition (A) or (B) in Definition
\ref{def: type of metric}. It induces a ray of metrics $\{\phi_{t}=\psi_{0}+u_{t}\}_{t\geq0}$
on $L$ with curvature form $\omega_{u_{t}}$. Let $\theta_{X}(u_{t})$
be the normalized Hamiltonian function of $X$ w.r.t. $\omega_{u_{t}}$.
Then we have
\begin{equation}
\lim_{t\rightarrow+\infty}\frac{1}{V}\int_{M}\dot{u}_{t}\theta_{X}(u_{t})\omega_{u_{t}}^{n}=\left\langle \alpha,\beta\right\rangle .\label{eq: limit to inner product}
\end{equation}
(2) (Integral formula for inner products) Let $(\mathcal{X},\mathcal{L})$
be a $\mathbb{C}^{*}$-equivariant test-configuration for $(M,L)$.
Assume $\mathcal{X}$ is smooth, given a smooth 2-form $\Omega\in2\pi c_{1}(\mathcal{L})$
and function $\Theta$ on $\mathcal{X}$ satisfying $\iota_{X}\Omega=i\overline{\partial}\Theta$
and $\int_{\mathcal{X}_{1}}\Theta\Omega^{n}=0$, then we have
\begin{equation}
\left\langle \alpha,\beta\right\rangle =\frac{1}{(n+1)L^{n}}\int_{\mathcal{X}}\Theta\left(\frac{\Omega}{2\pi}\right)^{n+1}.\label{eq: integral formula}
\end{equation}
\end{theorem}
\begin{proof}
If the total space $\mathcal{X}$ is singular, we take a $\mathbb{C}^{*}\times\mathbb{C}^{*}$-equivariant
resolution $p:\mathcal{X}'\rightarrow\mathcal{X}$ which is equivariant
and an isomorphism on $p^{-1}(\mathcal{X}\backslash\mathcal{X}_{0})$.
Let $\mathcal{L}'=p^{*}\mathcal{L}$, then $(\mathcal{X}',\mathcal{L}')$
has same inner product $\left\langle \alpha,\beta\right\rangle $
with $(\mathcal{X},\mathcal{L})$ by Proposition \ref{prop: invarience of inner product}.
Endowing $\mathcal{L}'$ with the pulling-back metric $p^{*}\Phi$,
it induces the same ray of metrics on $L$, so gives the same limit
on the LHS of (\ref{eq: limit to inner product}). So if we establish
(\ref{eq: limit to inner product}) for $(\mathcal{X}',\mathcal{L}')$,
then it also holds for $(\mathcal{X},\mathcal{L})$. In the following
we can assume $\mathcal{X}$ is smooth.

Since $\Omega=i\partial\bar{\partial}\Phi\in2\pi c_{1}(\mathcal{L})$,
$V=(2\pi)^{n}L^{n}$, by Proposition \ref{Prop: limit to global integral}
and (\ref{eq: modified Theta X}), the limit on the LHS of (\ref{eq: limit to inner product})
is equal to
\[
\frac{1}{L^{n}/n!}\int_{\mathcal{X}}\Theta_{X}\frac{(\Omega/2\pi)^{n+1}}{(n+1)!}-\frac{1}{\left(L^{n}/n!\right)^{2}}\frac{\mathcal{L}^{n+1}}{(n+1)!}\cdot\int_{\mathcal{X}_{1}}\Theta_{X}\frac{(\Omega/2\pi)^{n}}{n!}.
\]
Now we use equivariant HRR formula to relate these integrals to the
equivariant Euler characteristic.

For the integral over $\mathcal{X}_{1}$, apply (\ref{eq: HRR formula})
to datum $(M,L,\Phi\vert_{\mathcal{X}_{1}},\beta)$, it gives
\[
\sum_{\nu\in\mathbb{Z}}\nu\dim\textrm{H}^{0}(M,kL)_{\nu}=\int_{\mathcal{X}_{1}}\Theta_{X}\frac{(\Omega/2\pi)^{n}}{n!}\cdot k^{n+1}+O(k^{n}),\ \textrm{for}\ k\gg1.
\]
For the integral over $\mathcal{X}$, apply (\ref{eq: HRR formula})
to datum $(\mathcal{X},\mathcal{L},\Phi,\beta)$, it gives
\[
\chi_{1}^{\beta}(\mathcal{X},k\mathcal{L})=\int_{\mathcal{X}}\Theta_{X}\frac{(\Omega/2\pi)^{n+1}}{(n+1)!}\cdot k^{n+2}+O(k^{n+1}),\ \textrm{for}\ k\gg1.
\]
Combing these two expansions with the below lemma, which connects
$\chi_{1}^{\beta}(\mathcal{X},k\mathcal{L})$ to $\chi\left(\mathcal{X}_{\beta},k\mathcal{L}_{\beta}\right)$,
we have
\[
\mathcal{L}_{\beta}^{n+2}=(n+2)!\cdot\int_{\mathcal{X}}\Theta_{X}\frac{(\Omega/2\pi)^{n+1}}{(n+1)!}.
\]
Then (\ref{eq: limit to inner product}) follows from our Definition
(\ref{eq:def inner product}) of inner product.

For the part (2), choose a metric on $\mathcal{L}$ with curvature
form $\Omega$, then we obtain a Hamiltonian function $\Theta_{X}$
from the lifted action $\beta$. Its normalization $\tilde{\Theta}_{X}$
by (\ref{eq: modified Theta X}) must coincide with $\Theta$, then
(\ref{eq: integral formula}) follows in the same way as above.
\end{proof}
\begin{lemma} \label{lem: intersec num vs equi Euler num}Let $(\mathcal{X},\mathcal{L})$
be a $\mathbb{C}^{*}$-equivariant test-configuration for $(M,L)$
with fiberwise action $\beta$. $(\mathcal{X}_{\beta},\mathcal{L}_{\beta})$
is the fiber bundle (\ref{eq: projection F}) associated to $\beta$.
Then we have relation
\[
\chi\left(\mathcal{X}_{\beta},k\mathcal{L}_{\beta}\right)=\chi_{1}^{\beta}(\mathcal{X},k\mathcal{L})+\chi(\mathcal{X},k\mathcal{L}).
\]
Thus $\chi_{1}^{\beta}(\mathcal{X},k\mathcal{L})$ is a polynomial
of $k$ with degree at most $n+2$ when $k\gg1$. $\chi\left(\mathcal{X}_{\beta},k\mathcal{L}_{\beta}\right)$
and $\chi_{1}^{\beta}(\mathcal{X},k\mathcal{L})$ have same coefficient
of the leading term $k^{n+2}$. \end{lemma}
\begin{proof}
Consider the Leray's spectral sequence associated to $F:\mathcal{X}_{\beta}\rightarrow\mathbb{P}^{1}$
(see \cite{Harder} Section 4.6.3 for details),
\[
E_{2}^{p,q}=\textrm{H}^{p}(\mathbb{P}^{1},\textrm{R}^{q}F_{*}(k\mathcal{L}_{\beta}))\Rightarrow E^{p+q}=\textrm{H}^{p+q}(\mathcal{X}_{\beta},k\mathcal{L}_{\beta}).
\]
Since the Euler characteristic $\sum_{p,q}(-1)^{p+q}\dim E_{r}^{p,q}$
of each page are same with each other, we have
\[
\chi\left(\mathcal{X}_{\beta},k\mathcal{L}_{\beta}\right)=\sum_{q}(-1)^{q}\chi\left(\mathbb{P}^{1},\textrm{R}^{q}F_{*}(k\mathcal{L}_{\beta})\right).
\]
We claim that
\begin{equation}
\textrm{R}^{q}F_{*}(k\mathcal{L}_{\beta})=\bigoplus_{\nu\in\mathbb{Z}}\textrm{H}^{q}(\mathcal{X},k\mathcal{L})_{\nu}\otimes\mathcal{O}_{\mathbb{P}^{1}}(\nu).\label{eq: higher direct image}
\end{equation}
Decomposing $\textrm{H}^{q}(\mathcal{X},k\mathcal{L})$ with respect
to action $\beta$, we have
\begin{eqnarray*}
\textrm{R}^{q}F_{*}(k\mathcal{L}_{\beta}) & = & \left(\mathbb{C}^{2}\backslash\{0\}\right)\times_{\beta}\textrm{H}^{q}(\mathcal{X},k\mathcal{L})\\
 & = & \bigoplus_{\nu\in\mathbb{Z}}\left(\mathbb{C}^{2}\backslash\{0\}\right)\times_{\tau^{\nu}}\textrm{H}^{q}(\mathcal{X},k\mathcal{L})_{\nu},
\end{eqnarray*}
where the $\mathbb{C}^{*}$-action on $\textrm{H}^{q}(\mathcal{X},k\mathcal{L})_{\nu}$
is multiplication by character $\tau^{\nu}$. Since $\left(\mathbb{C}^{2}\backslash\{0\}\right)\times_{\tau^{\nu}}\mathbb{C}=\mathcal{O}_{\mathbb{P}^{1}}(\nu)$,
(\ref{eq: higher direct image}) follows.

Now taking the Euler characteristic of both sides of (\ref{eq: higher direct image}),
using $\chi(\mathbb{P}^{1},\mathcal{O}_{\mathbb{P}^{1}}(\nu))=\nu+1$,
we have
\begin{eqnarray*}
\chi\left(\mathcal{X}_{\beta},k\mathcal{L}_{\beta}\right) & = & \sum_{q}(-1)^{q}\sum_{\nu}(\nu+1)\dim\textrm{H}^{q}(\mathcal{X},k\mathcal{L})_{\nu}\\
 & = & \chi_{1}^{\beta}(\mathcal{X},k\mathcal{L})+\chi(\mathcal{X},k\mathcal{L}).
\end{eqnarray*}
The last statement follows by $\chi(\mathcal{X},k\mathcal{L})$ is
a polynomial of $k$ with degree $\leq n+1$.
\end{proof}

\subsection{Coincidence of two definitions for inner products }

\begin{theorem} \label{thm: coincide of inner product}Let $(\mathcal{X},\mathcal{L})$
be a $\mathbb{C}^{*}$-equivariant test-configuration for $(M,L)$.
If $\mathcal{L}$ is relatively ample, then $\left\langle \alpha,\beta\right\rangle $
defined by (\ref{eq:def inner product}) is same to $\left\langle \alpha,\beta\right\rangle _{0}$
defined by (\ref{eq:old def inner product}). \end{theorem}
\begin{proof}
We only need to show $c_{0}$ in expansion (\ref{eq: def of c_0})
is equal to $\mathcal{L}_{\beta}^{n+2}$. By Lemma \ref{lem: intersec num vs equi Euler num},
this is same to show $c_{0}$ is the leading coefficient of $\chi_{1}^{\beta}(\mathcal{X},k\mathcal{L})$.

Consider the Leray's spectral sequence associated to $\pi:\mathcal{X}\rightarrow\mathbb{P}^{1}$,
\[
E_{2}^{p,q}=\textrm{H}^{p}(\mathbb{P}^{1},\textrm{R}^{q}\pi_{*}(k\mathcal{L}))\Rightarrow E^{p+q}=\textrm{H}^{p+q}(\mathcal{X},k\mathcal{L}).
\]
Without loss of the information of actions, the spectral sequence
is taking in the abelian category of $\mathbb{C}^{*}\times\mathbb{C}^{*}$-modules,
namely each term $E_{r}^{p,q}$ is equipped with a $\mathbb{C}^{*}\times\mathbb{C}^{*}$-action
induced by $\alpha\times\beta$ and the differentials $d_{r}^{p,q}$
are equivariant w.r.t. these actions.

Since the differentials on page-$2$ are $d_{2}^{p,q}:E_{2}^{p,q}\rightarrow E_{2}^{p+2,q-1}$
and $E_{2}^{p,q}=0$ when $p\notin\{0,1\}$ (since $\dim\mathbb{P}^{1}=1$),
above spectral sequence actually degenerates at page-$2$. Combining
this with the fact that abelian category of $\mathbb{C}^{*}\times\mathbb{C}^{*}$-modules
is semi-simple (i.e. every short exact sequence splits), we obtain
\begin{eqnarray*}
\textrm{H}^{q}(\mathcal{X},k\mathcal{L}) & \cong & \textrm{H}^{0}(\mathbb{P}^{1},\textrm{R}^{q}\pi_{*}(k\mathcal{L}))\bigoplus\textrm{H}^{1}(\mathbb{P}^{1},\textrm{R}^{q-1}\pi_{*}(k\mathcal{L})),\ \textrm{for}\ q\geq1;\\
\textrm{H}^{0}(\mathcal{X},k\mathcal{L}) & \cong & \textrm{H}^{0}(\mathbb{P}^{1},\pi_{*}(k\mathcal{L})),
\end{eqnarray*}
all these isomorphisms are between $\mathbb{C}^{*}\times\mathbb{C}^{*}$-modules.

Now we use the assumption of ampleness. When $k$ is sufficiently
large and divisible, we have $\textrm{R}^{q}\pi_{*}(k\mathcal{L})=0$
for any $q>0$. Then above isomorphisms imply
\begin{equation}
\textrm{H}^{q}(\mathcal{X},k\mathcal{L})\cong\textrm{H}^{q}(\mathbb{P}^{1},\pi_{*}(k\mathcal{L})),\ \textrm{for}\ q=0,1;\ \textrm{H}^{q}(\mathcal{X},k\mathcal{L})=0,\ \textrm{for}\ q>1\label{eq:Leray spec isomorp}
\end{equation}
as $\mathbb{C}^{*}\times\mathbb{C}^{*}$-modules when $k\gg1$.

Let $\mathcal{E}\coloneqq\pi_{*}(k\mathcal{L})$ be the direct image,
since $\pi$ is flat, by \cite{Hartshorne} Corollary 12.9, when $k\gg1$,
$\mathcal{E}$ is a vector bundle over $\mathbb{P}^{1}$ with rank
$N_{k}$. Moreover, $\mathcal{E}$ is equipped with two commutative
$\mathbb{C}^{*}$-actions, still denoted by $\alpha$ and $\beta$.
By action $\beta$, $\mathcal{E}$ is decomposed to subbundles,
\[
\mathcal{E}=\bigoplus_{\nu\in\mathbb{Z}}\mathcal{E}_{\nu},
\]
where for each point $\tau\in\mathbb{P}^{1}$, $\mathcal{E}_{\nu}\vert_{\tau}$
is the weight-$\nu$ subspace of $\textrm{H}^{0}(\mathcal{X}_{\tau},k\mathcal{L}_{\tau})$.
On the other hand, since $\alpha$ commutes with $\beta$, thus $\alpha$
preserves each bundle $\mathcal{E}_{\nu}$.

Direct sum $\textrm{H}^{p}(\mathbb{P}^{1},\mathcal{E})=\bigoplus_{\nu}\textrm{H}^{p}(\mathbb{P}^{1},\mathcal{E}_{\nu})$
can be seen as the weight decomposition w.r.t. $\beta$, then isomorphism
(\ref{eq:Leray spec isomorp}) tells us $\textrm{H}^{p}(\mathcal{X},k\mathcal{L})_{\nu}=\textrm{H}^{p}(\mathbb{P}^{1},\mathcal{E}_{\nu})$.
This follows that
\begin{equation}
\chi_{1}^{\beta}(\mathcal{X},k\mathcal{L})=\sum_{p}(-1)^{p}\sum_{\nu}\nu\cdot\dim\textrm{H}^{p}(\mathbb{P}^{1},\mathcal{E}_{\nu})=\sum_{\nu}\nu\cdot\chi(\mathbb{P}^{1},\mathcal{E}_{\nu}).\label{eq:Euler number identity}
\end{equation}
Next we consider each subbundle $\mathcal{E}_{\nu}\rightarrow\mathbb{P}^{1}$
equipped with a $\mathbb{C}^{*}$-action $\alpha$. Note the action
on $\mathcal{E}_{\nu}\vert_{\infty}$ is trivial since the compactification
we made. Using the following lemma, take the vector bundle $E$ therein
to be $\mathcal{E}_{\nu}$. Since $\mathcal{E}_{\nu}\vert_{0}=\textrm{H}^{0}(\mathcal{X}_{0},k\mathcal{L}_{0})_{\nu}$,
it yields
\[
\chi(\mathbb{P}^{1},\mathcal{E}_{\nu})=\sum_{\mu}\mu\cdot\dim\textrm{H}^{0}(\mathcal{X}_{0},k\mathcal{L}_{0})_{\mu,\nu}+\textrm{rank}\mathcal{E}_{\nu}.
\]
Put this into (\ref{eq:Euler number identity}), we obtain
\[
\chi_{1}^{\beta}(\mathcal{X},k\mathcal{L})=\sum_{\mu,\nu}\mu\nu\cdot\dim\textrm{H}^{0}(\mathcal{X}_{0},k\mathcal{L}_{0})_{\mu,\nu}+\sum_{\nu}\nu\cdot\dim\textrm{H}^{0}(M,kL)_{\nu}.
\]
Since the second term on the RHS has lower degree than others, the
leading coefficient of $\chi_{1}^{\beta}(\mathcal{X},k\mathcal{L})$
is $c_{0}$, defined by (\ref{eq: def of c_0}).
\end{proof}
\begin{lemma} Let $E$ be a rank $r$ holomorphic vector bundle over
$\mathbb{P}^{1}=\mathbb{C}\cup\{\infty\}$. There is a $\mathbb{C}^{*}$-action
on $E$ covers the multiplication action on $\mathbb{P}^{1}$. Then
\[
\chi(\mathbb{P}^{1},E)=r+w_{0}-w_{\infty},
\]
where $w_{0}$ and $w_{\infty}$ is the total weight of action on
$E_{0}$ and $E_{\infty}$ respectively. \end{lemma}
\begin{proof}
The Hirzebruch-Riemann-Roch theorem says
\[
\chi(\mathbb{P}^{1},E)=\int_{\mathbb{P}^{1}}\textrm{ch}(E)\cdot\textrm{td}(T_{\mathbb{P}^{1}}).
\]
Since the Chern character $\textrm{ch}(E)=r+c_{1}(E)$, and the Todd
class $\textrm{td}(T_{\mathbb{P}^{1}})=1+\frac{1}{2}c_{1}(T_{\mathbb{P}^{1}})=1+c_{1}(\mathcal{O}(1))$.
Thus
\[
\chi(\mathbb{P}^{1},E)=r+\int_{\mathbb{P}^{1}}c_{1}(E)=r+\textrm{deg}(\wedge^{r}E).
\]
Denote the line bundle $\wedge^{r}E$ by $K$. Take any nonzero $s_{1}\in K_{1}$,
define a meromorphic section of $K$ by $s(\tau)=\tau.s_{1}$, for
$\tau\in\mathbb{C}^{*}$. It is easy to see that the divisor defined
by $s$ is $w_{0}\cdot0-w_{\infty}\cdot\infty$, where $w_{0}$ and
$w_{\infty}$ are weights of $\mathbb{C}^{*}$-action on $K_{0}$
and $K_{\infty}$ respectively, which are also the total weights of
$\mathbb{C}^{*}$-action on $E_{0}$ and $E_{\infty}$. The degree
of $K$ is $w_{0}-w_{\infty}$.
\end{proof}

\section{From Monge-Amp\`{e}re measures to Duistermaat-Heckman measures \label{sec:From-Monge-Ampere-measures}}

Let $(\mathcal{X},\mathcal{L})$ be a test-configuration for $(M,L)$
with a metric $\Phi$ on $\mathcal{L}$, and $\{\phi_{t}\}$ be the
induced ray of metrics on $L$. Suppose the curvature forms $\omega_{u_{t}}=i\partial\bar{\partial}\phi_{t}$
are nonnegative. Then the pushforwards of $\frac{1}{V}\omega_{u_{t}}^{n}$
by function $\dot{u}_{t}$ are probability measures on $\mathbb{R}$.
When $\{\phi_{t}\}$ is Phong-Sturm's geodesic ray, Hisamoto \cite{Hisa measure}
showed them weakly converge to $\textrm{DH}(\mathcal{X},\mathcal{L})$.
In this section, we extend this result to some general rays, such
as the ray induced by a smooth metric $\Phi$ with no curvature restriction.
The proof runs along a similar route as Theorem \ref{thm: limit slope},
hence is different from \cite{Hisa measure}.

\subsection{From the limit of integrals to the equivariant Euler characteristic}

In this subsection, $(\mathcal{X},\mathcal{L})$ is not necessarily
ample. The following proposition convert the limit into an integral
over the total space.

\begin{proposition} \label{prop: DH limit to integral}Let $(\mathcal{X},\mathcal{L})$
be a test-configuration for $(M,L)$ with smooth $\mathcal{X}$. Let
$\Phi$ be a $\alpha(\mathbb{S}^{1})$-invariant metric on $\mathcal{L}$
satisfying condition (A) or (B) in Definition \ref{def: type of metric}.
It induces a ray of metrics $\{\phi_{t}=\psi_{0}+u_{t}\}_{t\geq0}$
on $L$. Let $\Omega=i\partial\bar{\partial}\Phi$ be the curvature
current, and $\Theta_{W}$ be the Hamiltonian function of vector field
$W$ induced by the lifted action $\alpha$. Let $f$ be any $C^{1}$
function, then we have
\begin{equation}
\lim_{t\rightarrow\infty}\frac{1}{V}\int_{M}f(\dot{u}_{t})\omega_{u_{t}}^{n}=\frac{1}{L^{n}/n!}\int_{\mathcal{X}}f'(\Theta_{W})\frac{(\Omega/2\pi)^{n+1}}{(n+1)!}+f(0).\label{eq:limit to integral general f}
\end{equation}
\end{proposition}
\begin{proof}
The proof is similar with Proposition \ref{Prop: limit to global integral}.
We define $F(\tau)\coloneqq\int_{\mathcal{X}_{\tau}}f(\Theta_{W})\Omega^{n}$
for $\tau\in\mathbb{P}^{1}\backslash\{0\}$, it only depends on $\left\vert \tau\right\vert $.
We know $F(e^{-t/2})=\int_{M}f(\dot{u}_{t})\omega_{u_{t}}^{n}$, $F(\infty)=f(0)V$.
A direct computation shows that
\[
i\bar{\partial}\left(F(\tau)\frac{d\tau}{\tau}\right)=\frac{1}{n+1}\pi_{*}\left(f'(\Theta_{W})\Omega^{n+1}\right),\ \textrm{on}\ \mathbb{C}^{*}.
\]
Let $A_{r,R}\coloneqq\{\tau\mid r\leq\left\vert \tau\right\vert \leq R\}$,
by Stokes' formula, we have
\[
\int_{A_{r,R}}i\bar{\partial}\left(F(\tau)\frac{d\tau}{\tau}\right)=i\int_{A_{r,R}}d\left(F(\tau)\frac{d\tau}{\tau}\right)=2\pi\left(F(r)-F(R)\right),
\]
then let $r\rightarrow0$ and $R\rightarrow+\infty$, (\ref{eq:limit to integral general f})
follows.
\end{proof}
We will take $f=x^{d+1}$ and apply equivariant HRR formula to relate
the above integral over $\mathcal{X}$ to the equivariant Euler characteristic.

\begin{definition}For a test-configuration $(\mathcal{X},\mathcal{L})$,
the degree-$d$ equivariant Euler characteristic is defined by
\[
\chi_{d}^{\alpha}(\mathcal{X},k\mathcal{L})\coloneqq\sum_{q=0}^{n+1}(-1)^{q}\sum_{\mu\in\mathbb{Z}}\frac{\mu^{d}}{d!}\cdot\dim\textrm{H}^{q}(\mathcal{X},k\mathcal{L})_{\mu},
\]
where $\textrm{H}^{q}(\mathcal{X},k\mathcal{L})_{\mu}$ is the weight-$\mu$
subspace w.r.t. structure action $\alpha$. Note that we use different
action comparing with $\chi_{1}^{\beta}$. \end{definition}

Since we only apply HRR formula on smooth $\mathcal{X}$, same with
the proof of Theorem \ref{thm: limit slope}, we pass to the resolution
and then come back. For this, we need the following invariance of
the leading term of $\chi_{d}^{\alpha}(\mathcal{X},k\mathcal{L})$.

\begin{proposition} \label{prop: DH invariance of euler number}When
$k\gg1$, $\chi_{d}^{\alpha}(\mathcal{X},k\mathcal{L})$ is a polynomial
of $k$ with degree at most $n+1+d$. If $(\mathcal{X}',\mathcal{L}')$
is a pullback of $(\mathcal{X},\mathcal{L})$, the coefficient of
term $k^{n+1+d}$ of $\chi_{d}^{\alpha}(\mathcal{X},k\mathcal{L})$
is same with that of $\chi_{d}^{\alpha}(\mathcal{X}',k\mathcal{L}')$.
\end{proposition}
\begin{proof}
The proof is similar with Lemma \ref{lem: intersec num vs equi Euler num}.
We express the leading coefficient by intersection number.

Similar to the bundle construction (\ref{eq: projection F}), we take
the tautological $\mathbb{C}^{*}$-principal bundle $\mathbb{C}^{d+1}\backslash\{0\}\rightarrow\mathbb{P}^{d}$,
define fiber bundle
\[
F:\mathcal{X}_{\alpha,d}\coloneqq\left(\mathbb{C}^{d+1}\backslash\{0\}\right)\times_{\alpha}\mathcal{X}\rightarrow\mathbb{P}^{d}
\]
and line bundle $\mathcal{L}_{\alpha,d}$ on $\mathcal{X}_{\alpha,d}$
in the similar way.

Then consider the Leray spectral sequence associated to $F$, the
invariance of Euler characteristic for each page gives us
\[
\chi(\mathcal{X}_{\alpha,d},k\mathcal{L}_{\alpha,d})=\sum_{q}(-1)^{q}\chi\left(\mathbb{P}^{d},\textrm{R}^{q}F_{*}(k\mathcal{L}_{\alpha,d})\right).
\]
Similar with (\ref{eq: higher direct image}), we can show
\[
\textrm{R}^{q}F_{*}(k\mathcal{L}_{\alpha,d})=\bigoplus_{\mu\in\mathbb{Z}}\textrm{H}^{q}(\mathcal{X},k\mathcal{L})_{\mu}\otimes\mathcal{O}_{\mathbb{P}^{d}}(\mu).
\]
It follows that
\[
\chi(\mathcal{X}_{\alpha,d},k\mathcal{L}_{\alpha,d})=\sum_{\mu}\sum_{q}(-1)^{q}\chi(\mathbb{P}^{d},\mathcal{O}_{\mathbb{P}^{d}}(\mu))\cdot\dim\textrm{H}^{q}(\mathcal{X},k\mathcal{L})_{\mu}.
\]
Since $\chi(\mathbb{P}^{d},\mathcal{O}_{\mathbb{P}^{d}}(\mu))=\frac{\mu^{d}}{d!}+O(\mu^{d-1})$
is a polynomial of $\mu$, and $\chi(\mathcal{X}_{\alpha,d},k\mathcal{L}_{\alpha,d})$
is a polynomial of $k$ with degree $\leq n+1+d$, by induction on
$d$, we can show that $\chi_{d}^{\alpha}(\mathcal{X},k\mathcal{L})$
is a polynomial of $k$ with degree $\leq n+1+d$ when $k\gg1$. Moreover,
by above identity, we see the coefficient of term $k^{n+1+d}$ of
$\chi_{d}^{\alpha}(\mathcal{X},k\mathcal{L})$ is $\mathcal{L}_{\alpha,d}^{n+1+d}/(n+1+d)!$,
then the invariance of leading coefficients under pullback follows.
\end{proof}
Now we connect the limit with the equivariant Euler characteristic
via equivariant HRR formula.

\begin{proposition} \label{prop: DH limit to euler}With same assumptions
as Proposition \ref{prop: DH limit to integral}, except that $\mathcal{X}$
might be singular. For any integer $d\geq0$, we have
\begin{equation}
\lim_{t\rightarrow\infty}\frac{1}{V}\int_{M}\frac{(\dot{u}_{t})^{d+1}}{(d+1)!}\omega_{u_{t}}^{n}=\lim_{k\rightarrow\infty}\frac{\chi_{d}^{\alpha}(\mathcal{X},k\mathcal{L})}{k^{d+1}N_{k}}.\label{eq: DH limit to euler num}
\end{equation}
\end{proposition}
\begin{proof}
If $\mathcal{X}$ is singular, then take an equivariant resolution
$p:\mathcal{X}'\rightarrow\mathcal{X}$. Let $\mathcal{L}'=p^{*}\mathcal{L}$,
then $(\mathcal{X}',\mathcal{L}')$ is a pullback of $(\mathcal{X},\mathcal{L})$.
Endowing $\mathcal{L}'$ with pullback metric $p^{*}\Phi$, then the
LHS limits of (\ref{eq: DH limit to euler num}) are same for $\mathcal{X}'$
and $\mathcal{X}$. On the other hand, by Proposition \ref{prop: DH invariance of euler number},
the RHS leading coefficients are also same for $\mathcal{X}'$ and
$\mathcal{X}$. Hence it is sufficient to consider the case of smooth
$\mathcal{X}$.

Take $f(x)=x^{d+1}/(d+1)!$ in Proposition \ref{prop: DH limit to integral},
it yields
\[
\lim_{t\rightarrow\infty}\frac{1}{V}\int_{M}\frac{(\dot{u}_{t})^{d+1}}{(d+1)!}\omega_{u_{t}}^{n}=\frac{1}{L^{n}/n!}\int_{\mathcal{X}}\frac{\Theta_{W}^{d}}{d!}\frac{(\Omega/2\pi)^{n+1}}{(n+1)!}.
\]
Then apply equivariant HRR formula (\ref{eq: HRR formula} ) to datum
$(\mathcal{X},\mathcal{L},\Phi,\alpha)$, it yields
\[
\chi_{d}^{\alpha}(\mathcal{X},k\mathcal{L})=\int_{\mathcal{X}}\frac{\Theta_{W}^{d}}{d!}\frac{(\Omega/2\pi)^{n+1}}{(n+1)!}\cdot k^{n+1+d}+O(k^{n+d}),
\]
thus (\ref{eq: DH limit to euler num}) follows.
\end{proof}
The next step is to show the leading coefficient of $\chi_{d}^{\alpha}(\mathcal{X},k\mathcal{L})$
equals to $(d+1)$-moment of DH measure $\textrm{DH}(\mathcal{X},\mathcal{L})$.
For this, we need the associated filtration of section ring introduced
by Witt Nystr\"{o}m \cite{Nystrom TC-body}.

\subsection{Filtrations of the section ring}

The references of this subsection are \cite{Chen huayi,BHJ1,Nystrom TC-body}.

\begin{definition}[Filtrations]Let $(M,L)$ be a polarized manifold.
We call $\mathcal{\mathcal{F}}=\{F^{t}\textrm{H}^{0}(kL)\}_{t\in\mathbb{R},k\geq0}$
is a \textit{filtration} of the section ring $R(M,L)=\bigoplus_{k\geq0}\textrm{H}^{0}(M,kL)$,
if for each $k$, $\{F^{t}\textrm{H}^{0}(kL)\}_{t\in\mathbb{R}}$
is a family of subspaces of $\textrm{H}^{0}(M,kL)$ and nonincreasing
and left-continuous in $t$. A filtration $\mathcal{\mathcal{F}}$
is called \textit{admissible} if it is

(1) multiplicative: $F^{t}\textrm{H}^{0}(kL)\cdot F^{s}\textrm{H}^{0}(lL)\subset F^{t+s}\textrm{H}^{0}((k+l)L)$;

(2) pointwise left-bounded: for each $k$, $F^{t}\textrm{H}^{0}(kL)=\textrm{H}^{0}(kL)$
for $t$ small enough;

(3) linear right-bounded: $\exists C>0$ such that $F^{t}\textrm{H}^{0}(kL)=\{0\}$
for $t>Ck$ and $k\geq0$. \end{definition}

\begin{remark} Since $R(M,L)$ is finitely generated, by \cite{Chen huayi}
Lemma 1.5, any multiplicative and pointwise left-bounded filtration
is automatically linear left-bounded, i.e. there exists $C>0$ such
that $F^{t}\textrm{H}^{0}(kL)=\textrm{H}^{0}(kL)$ for $t<-Ck$ and
$k\geq0$. \end{remark}

Let $\mathcal{\mathcal{F}}$ be an admissible filtration of section
ring. For each $k$, to study the jumping numbers of filtration, we
consider probability measure
\[
\mathfrak{m}_{k}\coloneqq-\frac{d}{d\lambda}\left(\frac{1}{N_{k}}\dim F^{k\lambda}\textrm{H}^{0}(kL)\right),
\]
called the weight distribution. The linear boundness of $\mathcal{\mathcal{F}}$
implies that $\{\mathfrak{m}_{k}\}$ have uniformly bounded support.
Its weak limit can be obtained in the following way. For $\lambda\in\mathbb{R}$,
we define a graded subalgebra $R^{(\lambda)}$ of $R(M,L)$,
\[
R^{(\lambda)}\coloneqq\bigoplus_{k\geq0}F^{k\lambda}\textrm{H}^{0}(kL).
\]
Its volume is defined by
\begin{equation}
\textrm{vol}(R^{(\lambda)})\coloneqq\limsup_{k\rightarrow\infty}\frac{1}{N_{k}}\dim F^{k\lambda}\textrm{H}^{0}(kL).\label{eq: DH distribut function}
\end{equation}
By the linear boundness, $\textrm{vol}(R^{(\lambda)})=0$ when $\lambda\gg0$
and $\textrm{vol}(R^{(\lambda)})=1$ when $\lambda\ll0$.

\begin{theorem}[see Theorem 5.3 \cite{BHJ1}] \label{thm: LM is deriv of Vol}For
an admissible filtration $\mathcal{\mathcal{F}}$ of the section ring
$R(M,L)$, the sequence of weight distributions $\{\mathfrak{m}_{k}\}$
weakly converge to $-\frac{d}{d\lambda}\textrm{vol}(R^{(\lambda)})$,
where the derivative is taken in the sense of distributions. \end{theorem}

We call this weak limit of $\mathfrak{m}_{k}$ the \textit{limit measure}
of admissible filtration $\mathcal{\mathcal{F}}$, and denote it by
$\textrm{LM}(\mathcal{\mathcal{F}})$.

\subsection{Filtrations associated to test-configurations \label{subsec: Filtra of a TC}}

In \cite{Nystrom TC-body}, Witt Nystr\"{o}m define an admissible filtration
of $R(M,L)$ for a test-configuration for $(M,L)$.

Let $(\mathcal{X},\mathcal{L})$ be an ample test-configuration for
$(M,L)$. We use $(\mathcal{X},\mathcal{L})\vert_{\mathbb{C}}$ to
denote the restricted family over $\mathbb{C}$. For each $\mu\in\mathbb{Z}$
and $k\geq0$, let $F^{\mu}\textrm{H}^{0}(kL)$ be the image of following
restriction map,
\[
\textrm{H}^{0}(\mathcal{X}\vert_{\mathbb{C}},k\mathcal{L})_{\mu}\hookrightarrow\textrm{H}^{0}(\mathcal{X}_{1},k\mathcal{L})=\textrm{H}^{0}(M,kL).
\]
For $s\in\textrm{H}^{0}(M,kL)$, let $\bar{s}$ be the equivariant
(w.r.t. action $\alpha$) extension on $\mathcal{X}\vert_{\mathbb{C}^{*}}$.
Then we have an equivalent description:
\begin{equation}
F^{\mu}\textrm{H}^{0}(kL)=\{s\in\textrm{H}^{0}(M,kL)\mid\tau^{-\mu}\cdot\bar{s}\in\textrm{H}^{0}(\mathcal{X}\vert_{\mathbb{C}},k\mathcal{L})\}.\label{eq: filtration ind by t.c.}
\end{equation}
Namely, $F^{\mu}\textrm{H}^{0}(kL)$ is constituted of $s$ such that
$\tau^{-\mu}\cdot\bar{s}$ can be extended onto $\mathcal{X}\vert_{\mathbb{C}}$.
We can consider the restriction map to $\mathcal{X}_{0}$,
\[
F^{\mu}\textrm{H}^{0}(kL)\longrightarrow\textrm{H}^{0}(\mathcal{X}_{0},k\mathcal{L}_{0})_{\mu},\ s\mapsto\left(\tau^{-\mu}\bar{s}\right)\vert_{\mathcal{X}_{0}}.
\]
Since its kernel is $F^{\mu+1}\textrm{H}^{0}(kL)$, we have
\begin{equation}
F^{\mu}\textrm{H}^{0}(kL)/F^{\mu+1}\textrm{H}^{0}(kL)\cong\textrm{H}^{0}(\mathcal{X}_{0},k\mathcal{L}_{0})_{\mu}.\label{eq: quot of filtrat to X0}
\end{equation}

\begin{definition}[associated filtrations]For an ample test-configuration
$(\mathcal{X},\mathcal{L})$ of $(M,L)$, define the \textit{associated
filtration} is
\[
\mathcal{\mathcal{F}}(\mathcal{X},\mathcal{L})=\{F^{\left\lceil t\right\rceil }\textrm{H}^{0}(kL)\}_{t\in\mathbb{R},k\geq0},
\]
where $F^{\left\lceil t\right\rceil }\textrm{H}^{0}(kL)$ is given
by (\ref{eq: filtration ind by t.c.}). By \cite{Nystrom TC-body},
it is an admissible filtration of $R(M,L)$. \end{definition}

By (\ref{eq: quot of filtrat to X0}) and (\ref{eq: weight measure}),
we see the limit measure of $\mathcal{\mathcal{F}}(\mathcal{X},\mathcal{L})$
is DH measure $\textrm{DH}(\mathcal{X},\mathcal{L})$, and Theorem
\ref{thm: LM is deriv of Vol} implies
\begin{equation}
\textrm{DH}(\mathcal{X},\mathcal{L})=-\frac{d}{d\lambda}\textrm{vol}(R^{(\lambda)}).\label{eq: volume and DH}
\end{equation}

\begin{remark}\label{rem: Lost of nega part}Consider the relation
between $\mathcal{\mathcal{F}}(\mathcal{X},\mathcal{L})$ and sections
$\textrm{H}^{0}(\mathcal{X},k\mathcal{L})$ over the compactified
total space. For $s\in\textrm{H}^{0}(kL)$, since action $\alpha$
on $(\mathcal{X},\mathcal{L})_{\infty}$ is trivial, it is easy to
see $\tau^{-\mu}\cdot\bar{s}$ can be extended onto $\mathcal{X}_{\infty}$
only if $\mu\geq0$. Hence the relation is
\begin{equation}
\textrm{H}^{0}(\mathcal{X},k\mathcal{L})_{\mu}=F^{\mu}\textrm{H}^{0}(kL),\ \textrm{if}\ \mu\geq0\ \textrm{and}\ \textrm{H}^{0}(\mathcal{X},k\mathcal{L})_{\mu}=0,\ \textrm{if}\ \mu<0.\label{eq: H0 of total space to filtration}
\end{equation}
Thus if only take the sections over the compactified space, we can
not obtain the negative part of $\mathcal{\mathcal{F}}(\mathcal{X},\mathcal{L})$.
\end{remark}

\subsection{Proof of Theorem \ref{thm: lim to DH them}}

\begin{proof} It is sufficient to show that for any integer $d\geq0$,
we have
\begin{equation}
\lim_{t\rightarrow\infty}\frac{1}{V}\int_{M}\frac{(\dot{u}_{t})^{d+1}}{(d+1)!}\omega_{u_{t}}^{n}=\int_{\mathbb{R}}\frac{\lambda^{d+1}}{(d+1)!}d\textrm{DH}(\mathcal{X},\mathcal{L}).\label{eq: DH limit norm to d-moment}
\end{equation}
Using Proposition \ref{prop: DH limit to euler} as a bridge, we only
need to connect the leading coefficient of $\chi_{d}^{\alpha}(\mathcal{X},k\mathcal{L})$
to $\textrm{DH}(\mathcal{X},\mathcal{L})$. The associated filtration
$\mathcal{\mathcal{F}}(\mathcal{X},\mathcal{L})$ can connect them.
But there are two issues, one is that $\mathcal{\mathcal{F}}(\mathcal{X},\mathcal{L})$
does not involve $\textrm{H}^{q}(\mathcal{X},k\mathcal{L})_{\mu}$
($q>0$) but it appears in $\chi_{d}^{\alpha}$; another issue has
been mentioned in Remark \ref{rem: Lost of nega part}, the negative
part of $\mathcal{\mathcal{F}}(\mathcal{X},\mathcal{L})$ is lost
if we only take sections over $\mathcal{X}$. Fortunately, by the
relative ampleness of $\mathcal{L}$, both issues can be avoided.

Take a large integer $c$, let $\mathcal{L}_{c}\coloneqq\mathcal{L}+\pi^{*}\mathcal{O}_{\mathbb{P}^{1}}(c)$
equipped the product $\mathbb{C}^{*}$-action (action on $\mathcal{O}_{\mathbb{P}^{1}}(c)$
is defined in Remark \ref{rem: action on O(-1)}). Actually, $(\mathcal{X},\mathcal{L}_{c})$
is the compactification of $(\mathcal{X},\mathcal{L})\vert_{\mathbb{C}}$
with modified structure action by multiplying character $\tau^{c}$.

We endow $\mathcal{L}_{c}$ with metric $\Phi_{c}\coloneqq\Phi+c\cdot\pi^{*}\phi_{FS}$,
here $\phi_{FS}$ is the Fubini-Study metric on $\mathcal{O}_{\mathbb{P}^{1}}(-1)$.
The curvature current is $\Omega_{c}\coloneqq\Omega+c\pi^{*}\omega_{FS}$,
$\omega_{FS}=i\partial\bar{\partial}\log(1+\left\vert \tau\right\vert ^{2})$.
We can check the Hamiltonian function associated to datum $(\mathcal{X},\mathcal{L}_{c},\Phi_{c})$
is
\[
\Theta_{W,c}\coloneqq\Theta_{W}+\frac{c}{1+\left\vert \tau\right\vert ^{2}}.
\]
It vanishes on $\mathcal{X}_{\infty}$ as we expected.

Now applying Proposition \ref{prop: DH limit to euler} to datum $(\mathcal{X},\mathcal{L}_{c},\Phi_{c})$,
it yields
\[
\lim_{t\rightarrow\infty}\frac{1}{V}\int_{M}\frac{1}{(d+1)!}(\dot{u}_{t}+\frac{c}{1+e^{-t}})^{d+1}\omega_{u_{t}}^{n}=\lim_{k\rightarrow\infty}\frac{\chi_{d}^{\alpha}(\mathcal{X},k\mathcal{L}_{c})}{k^{d+1}N_{k}}.
\]
Since $\mathcal{L}$ is relatively ample, we can choose $c$ is sufficiently
large such that $\mathcal{L}_{c}$ is ample over $\mathcal{X}$. Then
for $k\gg1$, we have
\[
\chi_{d}^{\alpha}(\mathcal{X},k\mathcal{L})=\sum_{\mu\in\mathbb{Z}}\frac{\mu^{d}}{d!}\cdot\dim\textrm{H}^{0}(\mathcal{X},k\mathcal{L}_{c})_{\mu}.
\]
Let $F^{\mu}\textrm{H}^{0}(kL)$ be the filtration associated to $(\mathcal{X},\mathcal{L})$,
then it is easy to see the filtration associated to $(\mathcal{X},\mathcal{L}_{c})$
is
\[
\{F^{\mu-kc}\textrm{H}^{0}(kL)\}_{\mu\in\mathbb{Z}}.
\]
By the relation (\ref{eq: H0 of total space to filtration}) for $(\mathcal{X},\mathcal{L}_{c})$,
we have
\begin{eqnarray*}
\frac{1}{k^{d+1}N_{k}}\chi_{d}^{\alpha}(\mathcal{X},k\mathcal{L}_{c}) & = & \frac{1}{k^{d+1}N_{k}}\sum_{\mu\geq0}\frac{\mu^{d}}{d!}\cdot\dim F^{\mu-kc}\textrm{H}^{0}(kL)\\
(\textrm{let}\ \mu=k\lambda) & = & \sum_{0\leq\lambda\in\frac{1}{k}\mathbb{Z}}\frac{\lambda^{d}}{d!}\cdot\frac{1}{N_{k}}\dim F^{k(\lambda-c)}\textrm{H}^{0}(kL)\cdot\frac{1}{k}.
\end{eqnarray*}
Let $k\rightarrow\infty$, by definition (\ref{eq: DH distribut function}),
above sum converges to
\[
\int_{0}^{\infty}\frac{\lambda^{d}}{d!}\textrm{vol}(R^{(\lambda-c)})d\lambda=\int_{-c}^{\infty}\frac{(\lambda+c)^{d}}{d!}\textrm{vol}(R^{(\lambda)})d\lambda.
\]
We take $c$ to be large enough such that $\textrm{supp}\textrm{DH}(\mathcal{X},\mathcal{L})\subset[-c+1,\infty)$.
Integrating by parts and use (\ref{eq: volume and DH}), it further
equals to
\[
\int_{-c}^{\infty}\frac{(\lambda+c)^{d+1}}{(d+1)!}d\textrm{DH}(\mathcal{X},\mathcal{L})=\int_{\mathbb{R}}\frac{(\lambda+c)^{d+1}}{(d+1)!}d\textrm{DH}(\mathcal{X},\mathcal{L}).
\]
In a summary, we have showed that
\[
\lim_{t\rightarrow\infty}\frac{1}{V}\int_{M}\frac{1}{(d+1)!}(\dot{u}_{t}+\frac{c}{1+e^{-t}})^{d+1}\omega_{u_{t}}^{n}=\int_{\mathbb{R}}\frac{(\lambda+c)^{d+1}}{(d+1)!}d\textrm{DH}(\mathcal{X},\mathcal{L})
\]
for all integer $c\gg1$. Expand the LHS, we see it is a polynomial
in $c$, so it is for the RHS. Hence the constant terms of two sides
are same, that is (\ref{eq: DH limit norm to d-moment}). \end{proof}

\section{Mabuchi solitons and relative D-semistability \label{sec: Mabuchi-solitons-and Stability}}

\subsection{Extremal vector fields }

Let $M$ be a Fano manifold, set $L=-K_{M}$. There is a canonical
lifting of action $\textrm{Aut}(M)\hookrightarrow\textrm{Aut}(M,L)$.
By Kodaira's imbedding, the induced homomorphism $\textrm{Aut}(M,L)\rightarrow\textrm{GL}\left(\textrm{H}^{0}(M,kL)\right)$
will be injective when $k\gg1,$ thus $\textrm{Aut}(M)$ is a linear
algebraic group.

Let $K\subset\textrm{Aut}^{0}(M)$ be a maximal compact subgroup with
Lie algebra $\mathfrak{k}$. By Cartan-Iwasawa-Malcev theorem, $K$
is connected and unique up to a conjugation. Let $G\subset\textrm{Aut}^{0}(M)$
be its complexification with Lie algebra $\mathfrak{k}\otimes\mathbb{C}$.
In this paper, elements in Lie algebra are $(1,0)$-type holomorphic
vector fields, and the exponential map is
\[
\mathfrak{k}\rightarrow K,\ X\mapsto\exp(\textrm{Im}X),
\]
the latter is the time-$1$ map generated by $\textrm{Im}X$.

Then we take a $K$-invariant reference metric $\omega\in2\pi c_{1}(M)$,
let
\[
\mathcal{H}_{\omega}^{K}=\{u\in C^{\infty}(M)\mid\omega_{u}>0,\ \textrm{Im}X.u=0,\ \forall X\in\mathfrak{k}\}
\]
be the space of $K$-invariant K\"{a}hler potentials. Since $K$ is maximal,
for each $u\in\mathcal{H}_{\omega}^{K}$, the identity component of
the isometry group of $\omega_{u}$ is exactly $K$.

For any $X\in\mathfrak{k}$, let $\theta_{X}(\omega)$ be the normalized
Hamiltonian function of $X$ w.r.t. $\omega$, which satisfies $\iota_{X}\omega=i\bar{\partial}\theta_{X}(\omega)$
and $\int\theta_{X}(\omega)\omega^{n}=0$. For any $u\in\mathcal{H}_{\omega}^{K}$,
let $\theta_{X}(u)\coloneqq\theta_{X}(\omega)+X(u)$, we have $\iota_{X}\omega_{u}=i\bar{\partial}\theta_{X}(u)$.
By the below proposition, $\theta_{X}(u)$ is the normalized Hamiltonian
function of $X$ w.r.t. $\omega_{u}$.

\begin{proposition} (1) For each $X\in\mathfrak{k}$, the pushforward
of measure $\omega_{u}^{n}$ under the map $\theta_{X}(u):M\rightarrow\mathbb{R}$
is independent of $u\in\mathcal{H}_{\omega}^{K}$. In particular,
its barycenter is $\int\theta_{X}(u)\omega_{u}^{n}=\int\theta_{X}(\omega)\omega^{n}=0$.

(2) For each $X,\ Y\in\mathfrak{k}$, the pushforward of measure $\omega_{u}^{n}$
under the map $\left(\theta_{X}(u),\theta_{Y}(u)\right):M\rightarrow\mathbb{R}^{2}$
is also independent of $u\in\mathcal{H}_{\omega}^{K}$. \end{proposition}
\begin{proof}
For (1), we only need to show that for any $f\in C^{1}(\mathbb{R})$
and $u\in\mathcal{H}_{\omega}^{K}$, we have
\[
\int_{M}f\left(\theta_{X}(u)\right)\omega_{u}^{n}=\int_{M}f\left(\theta_{X}(\omega)\right)\omega^{n}.
\]
For $t\in[0,1]$, $tu\in\mathcal{H}_{\omega}^{K}$, we compute
\begin{eqnarray*}
\frac{d}{dt}\int_{M}f\left(\theta_{X}(tu)\right)\omega_{tu}^{n} & = & \int f'X(u)\omega_{tu}^{n}+\int f\left(\theta_{X}(tu)\right)i\partial\bar{\partial}u\wedge n\omega_{tu}^{n-1}\\
 & = & \int f'X(u)\omega_{tu}^{n}+\int f'i\bar{\partial}\theta_{X}(tu)\wedge\partial u\wedge n\omega_{tu}^{n-1}\\
 & = & \int f'X(u)\omega_{tu}^{n}+\int f'\iota_{X}\omega_{tu}\wedge\partial u\wedge n\omega_{tu}^{n-1}=0.
\end{eqnarray*}
The proof of (2) is similar with (1).
\end{proof}
\begin{definition}[DH measure]Let $K\subset\textrm{Aut}^{0}(M)$
be a maximal compact subgroup with Lie algebra $\mathfrak{k}$. Since
the part (1) of above proposition, we denote the pushforward measure
by
\[
\textrm{DH}_{K}(X)\coloneqq\theta_{X}(\omega)_{\#}\frac{\omega^{n}}{V},
\]
which is independent of the choice of $K$-invariant metric $\omega\in2\pi c_{1}(M)$.
We call it the \textit{DH measure} of vector field $X\in\mathfrak{k}$.

For $X,\ Y\in\mathfrak{k}$ and any $u\in\mathcal{H}_{\omega}^{K}$,
we define
\[
B_{K}(X,Y)\coloneqq\frac{1}{V}\int_{M}\theta_{X}(u)\theta_{Y}(u)\omega_{u}^{n}.
\]
By part (2) of the above proposition, it is independent of the choice
of $u$ and an inner product on $\mathfrak{k}$. This is the restriction
of the bilinear form introduced by Futaki-Mabuchi \cite{Futaki-Mabuchi}.
\end{definition}

\begin{remark} \label{rem: change K}Let $K'=gKg^{-1}$, $g\in\textrm{Aut}^{0}(M)$
be another choice of maximal compact subgroups. Then $\mathfrak{k}'=g_{*}\mathfrak{k}$
and $\omega'=(g^{-1})^{*}\omega$ is $K'$-invariant. For any $X\in\mathfrak{k}$,
let $X'=g_{*}X\in\mathfrak{k}'$, it is easy to see that
\[
\theta_{X'}(\omega')=\theta_{X}(\omega)\circ g^{-1}.
\]
This follows that $\textrm{DH}_{K'}(X')=\textrm{DH}_{K}(X)$. For
the inner product, let $X,\ Y\in\mathfrak{k}$ and $X'=g_{*}X,\ Y'=g_{*}Y\in\mathfrak{k}'$,
we have $B_{K'}(X',Y')=B_{K}(X,Y)$. \end{remark}

For a holomorphic vector field $X$, let $\omega\in2\pi c_{1}(M)$
with Ricci potential $h_{\omega}$, the Futaki invariant is defined
by
\[
F(X)\coloneqq\frac{1}{V}\int X(h_{\omega})\omega^{n}.
\]
It is independent of the choice of $\omega$. For a maximal compact
subgroup $K\subset\textrm{Aut}^{0}(M)$, by taking $\omega$ to be
$K$-invariant, we see the restriction of $F$ on $\mathfrak{k}$
is real-valued.

\begin{definition}[extremal vector fields and invariant $\vartheta(M)$]\label{def: VarTheta M}

Let $K\subset\textrm{Aut}^{0}(M)$ be a maximal compact subgroup,
there is a unique element $Z_{K}\in\mathfrak{k}$ such that
\[
F(X)=-B_{K}(X,Z_{K}),\ \textrm{for}\ \forall X\in\mathfrak{k}.
\]
We call $Z_{K}$ the \textit{extremal vector field} associated to
$K$. If $K'=gKg^{-1}$, $g\in\textrm{Aut}^{0}(M)$ is another maximal
compact subgroup, since $F(g_{*}X)=F(X)$, we have $Z_{K'}=g_{*}Z_{K}$.
Then by Remark \ref{rem: change K}, $\textrm{DH}_{K}(Z_{K})$ is
independent of the choice of $K$. Thus we omit ``$K$'' and denote
it by
\[
\textrm{DH}(M)\coloneqq\theta_{Z_{K}}(\omega)_{\#}\frac{\omega^{n}}{V},
\]
where $\omega$ is any $K$-invariant metric $\omega\in2\pi c_{1}(M)$.
$\textrm{DH}(M)$ is an invariant for Fano manifolds. In this paper,
we mostly concern the invariant
\[
\vartheta(M)\coloneqq\max_{M}\theta_{Z_{K}}(\omega)=\sup\textrm{supp}\textrm{DH}(M)-\int_{\mathbb{R}}\lambda d\textrm{DH}(M),
\]
where the barycenter of $\textrm{DH}(M)$ is the origin. It had been
introduced by Mabuchi \cite{Mabuchi} with notation $\alpha_{M}$.

\end{definition}

The extremal vector field has many important properties. By \cite{Futaki-Mabuchi},
$Z_{K}$ belongs to the center of $\mathfrak{k}$ and its imaginary
part generates a circle action, which we call the \textit{extremal
action}.

There is another equivalent definition of $Z_{K}$ in \cite{Futaki-Mabuchi}
in terms of scalar curvatures. Fix $K\subset\textrm{Aut}^{0}(M)$,
and take a $K$-invariant metric $\omega$. We endow the space $C^{\infty}(M,\mathbb{R})$
with $L^{2}$-inner product $\frac{1}{V}\int f_{1}f_{2}\omega^{n}$,
and then consider the orthogonal projection to the subspace
\[
P_{\omega}:C^{\infty}(M,\mathbb{R})\rightarrow\{\theta_{X}(\omega)\mid X\in\mathfrak{k}\}.
\]
Since for any $X\in\mathfrak{k}$, we have
\[
F(X)=-\int h_{\omega}\mathcal{L}_{X}\frac{\omega^{n}}{V}=-\int h\triangle\theta_{X}(\omega)\frac{\omega^{n}}{V}=-\int\theta_{X}(\omega)\left(S(\omega)-n\right)\frac{\omega^{n}}{V}.
\]
It follows that $\theta_{Z_{K}}(\omega)=P_{\omega}\left(S(\omega)-n\right)$,
thus
\[
Z_{K}=\textrm{grad}_{\omega}P_{\omega}\left(S(\omega)-n\right)
\]
Moreover, Theorem 2.1 in \cite{Mabuchi} says $P_{\omega}\left(S(\omega)-n\right)=P_{\omega}\left(1-e^{h_{\omega}}\right)$,
hence $Z_{K}$ also can be defined by
\begin{equation}
Z_{K}=\textrm{grad}_{\omega}P_{\omega}\left(1-e^{h_{\omega}}\right).\label{eq: def of Z by proj}
\end{equation}

\subsection{Mabuchi solitons }

There is an analogue of Calabi energy. Let $\omega\in2\pi c_{1}(M)$
be a reference metric and $\mathcal{H}_{\omega}=\{u\in C^{\infty}(M)\mid\omega_{u}>0\}$.
The Ding energy is defined by
\[
Din(u)=\frac{1}{V}\int_{M}\left(e^{h(u)}-1\right)^{2}\omega_{u}^{n},\ u\in\mathcal{H}_{\omega},
\]
 where $h(u)$ is the Ricci potential of $\omega_{u}$.

\begin{definition}If $u\in\mathcal{H}_{\omega}$ is a critical point
of Ding energy, then we say the associated metric $\omega_{u}$ is
a \textit{Mabuchi soliton.} \end{definition}

In this sense, Mabuchi solitons are analogues of extremal metrics
which are critical points of Calabi energy. This definition is equivalent
to that one given by Mabuchi \cite{Mabuchi}.

\begin{theorem} \cite{Yao} $\omega_{u}$ is a Mabuchi soliton if
and only if the vector field $\textrm{grad}_{\omega_{u}}(1-e^{h(u)})$
is holomorphic. \end{theorem}

Moreover, Mabuchi soliton is unique modulo the action of $\textrm{Aut}^{0}(M)$,
see Theorem C in \cite{Mabuchi-multiplier} and Corollary B in \cite{Mabuchi-CM}.
As with extremal metrics, the intersection of the isometry group of
a Mabuchi soliton with $\textrm{Aut}^{0}(M)$ is maximally compact
in $\textrm{Aut}^{0}(M)$, see Theorem D in \cite{Mabuchi-CM}. Clearly,
if $M$ admits a Mabuchi soliton, then for any maximal compact subgroup
$K$, there exists a $K$-invariant Mabuchi soliton, since all maximal
compact subgroups are conjugated with each other. Hence in the search
for a Mabuchi soliton, we firstly fix a maximal compact subgroup $K\subset\textrm{Aut}^{0}(M)$,
and take the reference metric $\omega\in2\pi c_{1}(M)$ is $K$-invariant.
Let $\mathcal{H}_{\omega}^{K}$ be the space of $K$-invariant K\"{a}hler
potentials.

Suppose $u\in\mathcal{H}_{\omega}^{K}$ such that $X=\textrm{grad}_{u}(1-e^{h(u)})$
is holomorphic. Since $\theta_{X}(u)=1-e^{h(u)}$ is real-valued,
$\textrm{Im}X$ preserves $\omega_{u}$. Since the isometry group
of $\omega_{u}$ is $K$, so $X\in\mathfrak{k}$. Then by (\ref{eq: def of Z by proj})
we conclude $X=Z_{K}$. After adding $u$ by a constant, we see $u$
satisfies a Monge-Amp\`{e}re type equation
\[
(1-\theta_{Z}(u))\omega_{u}^{n}=e^{h_{\omega}-u}\omega^{n}.
\]
A necessary condition for it admits a smooth solution is $\max_{M}\theta_{Z}(u)<1$,
i.e. $\vartheta(M)<1$.

\subsection{Modified Ding functionals \label{subsec:Modi Ding func}}

Fix $K\subset\textrm{Aut}^{0}(M)$ and a $K$-invariant reference
metric $\omega\in2\pi c_{1}(M)$. Let $Z=Z_{K}$ be the associated
extremal vector field. We define
\[
\mathcal{H}_{\omega}^{Z}=\{u\in C^{\infty}(M)\mid\omega_{u}>0,\ \textrm{Im}Z.u=0\}\supset\mathcal{H}_{\omega}^{K}.
\]
Then $\theta_{Z}(u)$ is real-valued for each $u\in\mathcal{H}_{\omega}^{Z}$.

\begin{definition} \label{Def of modi Ding}(1) The \textit{modified
Monge-Amp\`{e}re Energy} $E_{Z}:\mathcal{H}_{\omega}^{Z}\rightarrow\mathbb{R}$
is defined such that $E_{Z}(0)=0$ and
\[
\delta E_{Z}\vert_{u}(\delta u)=\frac{1}{V}\int_{M}\delta u\left(1-\theta_{Z}(u)\right)\omega_{u}^{n}.
\]

(2) The \textit{modified Ding functional} $D_{Z}:\mathcal{H}_{\omega}^{Z}\rightarrow\mathbb{R}$
is defined by
\begin{equation}
D_{Z}(u)\coloneqq-E_{Z}(u)+L(u),\ L(u)\coloneqq-\log\left(\frac{1}{V}\int_{M}e^{h_{\omega}-u}\omega^{n}\right).\label{eq: modif Ding}
\end{equation}
When $Z=0$, $D_{Z}$ is the ordinary Ding functional. \end{definition}

For $u\in\mathcal{H}_{\omega}^{Z}$, integrating along the path $u_{t}=tu$
for $t\in[0,1]$, we have
\[
E_{Z}(u)=\int_{0}^{1}dt\int_{M}u\left(1-\theta_{Z}(\omega)-t\cdot Zu\right)\left(t\omega_{u}+(1-t)\omega\right)^{n}.
\]
Expanding the volume form and integrating $t$, we can obtain an explicit
formula for $E_{Z}$.

By definition, the variation of $D_{Z}$ is
\[
\delta D_{Z}\vert_{u}(\delta u)=-\frac{1}{V}\int_{M}\delta u\left(1-e^{h(u)}-\theta_{Z}(u)\right)\omega_{u}^{n},
\]
thus its critical points are Mabuchi solitons.

Moreover, $D_{Z}$ can be extended to the space of $\textrm{Im}Z$-invariant
bounded psh metrics on $K_{M}^{-1}$, which is denoted by
\[
\mathcal{H}_{b}^{Z}=\{\phi\mid\textrm{Im}Z\textrm{-invariant, bounded psh metric on }K_{M}^{-1}\}.
\]
By \cite{Berman-Nystrom} Proposition 2.17, $E_{Z}$ is convex along
bounded sub-geodesics in $\mathcal{H}_{b}^{Z}$ and affine along geodesics.
Since $L(u)$ is also convex along bounded geodesics due to Berndtsson,
thus $D_{Z}$ is convex along bounded geodesics.

\subsection{Modified K-energies}

We only consider extremal metrics in $c_{1}(M)$ for a Fano manifold
$M$. Let $K\subset\textrm{Aut}^{0}(M)$ be a maximal compact subgroup
with associated extremal vector field $Z$. \textit{The modified K-energy}
$M_{Z}:\mathcal{H}_{\omega}^{Z}\rightarrow\mathbb{R}$ is defined
such that $M_{Z}(0)=0$ and
\[
\delta M_{Z}\vert_{u}(\delta u)=-\frac{1}{V}\int_{M}\delta u\left(S(\omega_{u})-n-\theta_{Z}(u)\right)\omega_{u}^{n}.
\]
Integrating along a smooth path $\{u_{t}\}$ from $0$ to $u$ (e.g.
$u_{t}=tu$), we have
\[
M_{Z}(u)=-\frac{1}{V}\int_{0}^{1}dt\int_{M}\dot{u_{t}}\left(S(\omega_{u_{t}})-n-1\right)\omega_{u_{t}}^{n}-E_{Z}(u)
\]
We derive the Chen-Tian formula of $M_{Z}$. Let $\Omega_{u}=e^{h_{\omega}-u}\omega^{n}$,
then $Ric(\Omega_{u})=\omega_{u}$. Note that $-\int_{M}\dot{u_{t}}\left(S(\omega_{u_{t}})-n\right)\omega_{u_{t}}^{n}$
is equal to
\begin{eqnarray*}
 &  & -n\int\dot{u_{t}}\left(Ric(\omega_{u_{t}})-Ric(\Omega_{u_{t}})\right)\wedge\omega_{u_{t}}^{n-1}=n\int\dot{u_{t}}\ i\partial\bar{\partial}\log\frac{\omega_{u_{t}}^{n}}{\Omega_{u_{t}}}\wedge\omega_{u_{t}}^{n-1}\\
 & = & \int\log\frac{\omega_{u_{t}}^{n}}{\Omega_{u_{t}}}\frac{d}{dt}\omega_{u_{t}}^{n}=\frac{d}{dt}\left[\int\log\left(\frac{\omega_{u_{t}}^{n}}{\Omega_{u_{t}}}\right)\omega_{u_{t}}^{n}\right]-\int\dot{u_{t}}\omega_{u_{t}}^{n}.
\end{eqnarray*}
Putting this back and integrate over $t$, we obtain
\begin{eqnarray}
M_{Z}(u) & = & \frac{1}{V}\int\log\left(\frac{\omega_{u}^{n}}{\Omega_{u}}\right)\omega_{u}^{n}+\frac{1}{V}\int h_{\omega}\omega^{n}-E_{Z}(u)\label{eq: modif K and Ding}\\
 & = & D_{Z}(u)+\int\log\left(\frac{\omega_{u}^{n}/V}{\Omega_{u}/\int\Omega_{u}}\right)\frac{\omega_{u}^{n}}{V}+\frac{1}{V}\int h_{\omega}\omega^{n}\nonumber \\
 & \geq & D_{Z}(u)+\frac{1}{V}\int h_{\omega}\omega^{n},\nonumber
\end{eqnarray}
where the entropy term is nonnegative by Jensen's inequality.

\begin{remark}\label{rem extK and MS} The above inequality for $M_{Z}$
and $D_{Z}$ implies an implication relation between extremal metrics
and Mabuchi solitons. If Fano manifold $M$ admits a Mabuchi soliton,
by \cite{LiY ZhouB} or \cite{Li-HanJY} $D_{Z}$ is proper in the
sense therein, the above inequality implies $M_{Z}$ is also proper.
Then by the recent breakthrough by Chen-Cheng \cite{Chen-Cheng} (also
see \cite{HeWY}), $M$ admits extremal metrics. Conversely, toric
3-fold $\mathbb{P}\left(\mathcal{O}_{\mathbb{P}^{2}}\oplus\mathcal{O}_{\mathbb{P}^{2}}(2)\right)$
does not admit Mabuchi solitons (it is relative D-unstable, see \cite{Yotsutani}),
but as a projective bundle it admits extremal metrics by the construction
in \cite{Apostolov}. \end{remark}

\subsection{Limit slopes and relative D-semistability \label{subsec: rel BD inv}}

By \cite{Futaki-Mabuchi} Theorem F, $\textrm{Im}Z$ generates a circle
action, then there exists a minimum $\tau=\tau_{Z}>0$ such that $\exp4\pi\tau\cdot\textrm{Im}Z=Id$.
Thus $\tau Z$ generates a $\mathbb{C}^{*}$-action $\beta_{\tau Z}$
satisfying
\begin{equation}
\beta_{\tau Z}(e^{-\frac{1}{2}(t+is)})=\exp\left(t\cdot\tau\textrm{Re}Z-s\cdot\tau\textrm{Im}Z\right).\label{eq: def extremal action}
\end{equation}
We call it the extremal action associated to maximal compact subgroup
$K$. We use ``$\beta_{Z}$'' to denote the \textit{formal} $\mathbb{C}^{*}$-action
generated by $Z$.

Assume $(\mathcal{X},\mathcal{L})$ is a $\mathbb{C}^{*}$-equivariant
test-configuration with fiberwise $\mathbb{C}^{*}$-action $\beta_{\tau Z}$,
we define the inner product of the structure $\mathbb{C}^{*}$-action
$\alpha$ and the formal $\mathbb{C}^{*}$-action $\beta_{Z}$ is
\[
\left\langle \alpha,\beta_{Z}\right\rangle \coloneqq\frac{1}{\tau}\left\langle \alpha,\beta_{\tau Z}\right\rangle .
\]
For simplicity, we will not explicitly mention the period $\tau$
and just take $\beta_{Z}$ to be a real $\mathbb{C}^{*}$-action.

\begin{proposition} Let $(\mathcal{X},\mathcal{L})$ be a $\mathbb{C}^{*}$-equivariant
ample test-configuration for $(M,-K_{M})$ with structure action $\alpha$
and fiberwise extremal action $\beta_{Z}$. Let $\Phi$ be a $\alpha(\mathbb{S}^{1})\times\beta(\mathbb{S}^{1})$-invariant
metric on $\mathcal{L}$ satisfying condition (B) in Definition \ref{def: type of metric}.
Pulling back by $\alpha$, $\Phi$ induces a ray $\{u_{t}\}_{t\geq0}\subset\mathcal{H}_{b}^{Z}$.
Then we have
\[
\lim_{t\rightarrow\infty}\frac{d}{dt}E_{Z}(u_{t})=E^{NA}(\mathcal{X},\mathcal{L})-\left\langle \alpha,\beta_{Z}\right\rangle
\]
and
\[
\lim_{t\rightarrow\infty}\frac{d}{dt}D_{Z}(u_{t})=D^{NA}(\mathcal{X},\mathcal{L})+\left\langle \alpha,\beta_{Z}\right\rangle .
\]
\end{proposition}
\begin{proof}
The first identity follows from Theorem \ref{thm: limit slope} and
Theorem \ref{thm: lim to DH them}. Then by the formula of limit slope
of $L(u_{t})$ due to Berman \cite{Berman}, the second identity follows.
\end{proof}
\begin{definition}[relative D-semistability] \label{def: def of D-semistable}
Let $M$ be a Fano manifold. Let $K\subset\textrm{Aut}^{0}(M)$ be
a compact subgroup with complexification $G=K^{\mathbb{C}}\subset\textrm{Aut}^{0}(M)$.
Assume $K^{m}\supset K$ is a maximal compact subgroup such that $Z\coloneqq Z_{K^{m}}$
belongs to $\textrm{Lie}(K)$. In this paper, we take $G$ to be a
torus.

(1) Let $(\mathcal{X},\mathcal{L})$ be a $G$-equivariant ample test-configuration
for $(M,-K_{M})$, the \textit{relative Berman-Ding invariant} is
defined by
\[
D_{Z}^{NA}(\mathcal{X},\mathcal{L})\coloneqq D^{NA}(\mathcal{X},\mathcal{L})+\left\langle \alpha,\beta_{Z}\right\rangle ,
\]
where $\left\langle \alpha,\beta_{Z}\right\rangle $ is defined by
taking $(\mathcal{X},\mathcal{L})$ as a $\mathbb{C}^{*}$-equivariant
test-configuration via $\beta_{Z}:\mathbb{C}^{*}\rightarrow G$.

(2) We say $M$ is \textit{D-semistable relative to} $G$, if for
any $G$-equivariant ample test-configuration $(\mathcal{X},\mathcal{L})$,
we always have $D_{Z}^{NA}(\mathcal{X},\mathcal{L})\geq0$. \end{definition}

The existence of Mabuchi solitons can imply that $M$ is D-semistable
relative to $G$. The argument is similar to the KE situation \cite{Berman}.
If $M$ admits a $K$-invariant Mabuchi solitons, we consider the
Phong-Sturm's geodesic ray $\{u_{t}\}$ associated to a $G$-equivariant
test-configurations and starting from that Mabuchi solitons. Since
$\frac{d}{dt}D_{Z}(u_{t})$ is nondecreasing, then we have
\[
0\leq\frac{d}{dt}D_{Z}(u_{t})\vert_{t=0}\leq\lim_{t\rightarrow\infty}\frac{d}{dt}D_{Z}(u_{t})=D_{Z}^{NA}(\mathcal{X},\mathcal{L}).
\]

\section{Relative D-semistability implies $\vartheta(M)\protect\leq1$ \label{sec: theta<=00003D1}}

In this section, we show that $M$ is relatively D-semistable can
imply $\vartheta(M)=\sup_{M}\theta_{Z}(\omega)\leq1$.

Let $T\subset\textrm{Aut}^{0}(M)$ be a $m$-dimensional complex torus
which complexifies compact torus $S$. Suppose $K$ is a maximal compact
subgroup containing $S$ and such that the associated extremal vector
field $Z$ belongs to $\textrm{Lie}(S)$. Take a $K$-invariant K\"{a}hler
metric $\omega\in2\pi c_{1}(M)$. Let $\theta$ the normalized Hamiltonian
function of $Z$ satisfying $\iota_{Z}\omega=i\bar{\partial}\theta$
and $\int\theta\omega^{n}=0$.

\begin{theorem} \label{thm: semistable to theta 1}Let $M$ be a
Fano manifold, $T\subset\textrm{Aut}^{0}(M)$ be a complex torus.
If $M$ is D-semistable relative to $T$, then $\vartheta(M)\leq1$.
\end{theorem}

The proof is by constructing a specific $T$-equivariant test-configuration:
the deformation to normal cone. This kind of test-configurations have
been systematically studied in \cite{Ross-Thomas}. First we need
to choose the center of blowing-up.

\begin{proposition} \label{prop: exist fixed pt}There exists a $T$-fixed
point $z_{*}\in M$ attains the maximum value of $\theta$ i.e. $\vartheta(M)$.
\end{proposition}
\begin{proof}
Since $T$-action can be canonically lifted to $-K_{M}$ and $\omega$
is $S$-invariant, it induces a moment map $\mu:M\rightarrow\triangle$
for $S$-action, where $\triangle\subset\mathbb{R}^{m}$ is a convex
set. Since $Z\in\textrm{Lie}(S)$, $\theta$ can be represented as
$l\circ\mu$ by an affine function $l$ on $\mathbb{R}^{m}$. By \cite{Atiyah},
$\triangle$ is the convex hull of image of $T$-fixed points under
$\mu$, thus the maximum of $l\vert_{\triangle}$ can be attained
at the image of some $T$-fixed point, say $z_{*}$, it has desired
properties.
\end{proof}
\begin{proof}[Proof of Theorem \ref{thm: semistable to theta 1}]Consider
the deformation to normal cone of $z_{*}$ (provided by above lemma).
First, take product $(M\times\mathbb{P}^{1},p_{1}^{*}K_{M}^{-1})$,
the structure action $\alpha$ is the multiplication on $\mathbb{P}^{1}$,
the fiberwise $T$-action $\beta$ is the given torus $T\subset\textrm{Aut}^{0}(M)$.
Let
\[
\Pi:\mathcal{X}\rightarrow M\times\mathbb{P}^{1}
\]
be the blowup at point $(z_{*},0)$ with exceptional divisor $P\cong\mathbb{P}^{n}$.
Let
\[
\pi=p_{2}\circ\Pi:\mathcal{X}\rightarrow\mathbb{P}^{1},\ p=p_{1}\circ\Pi:\mathcal{X}\rightarrow M.
\]
For $c\in\mathbb{Q}_{+}$, define a $\mathbb{Q}$-line bundle $\mathcal{L}_{c}=p^{*}K_{M}^{-1}-cP$
on $\mathcal{X}$. When $c\ll1$, $\mathcal{L}_{c}$ is relatively
ample. Since action $\alpha$ and $\beta$ fixes $(z_{*},0)\in M\times\mathbb{P}^{1}$,
they can be lifted to $\mathcal{X}$. Since $P$ is invariant, they
can be further lifted to $\mathcal{L}_{c}$ such that $\alpha$ trivially
acts on the fiber over $\infty$ and $\beta$ is the given torus $T\subset\textrm{Aut}^{0}(M)$
on fibers $(\mathcal{X},\mathcal{L}_{c})\vert_{\tau\neq0}\cong(M,K_{M}^{-1})$.

In a summary, $(\mathcal{X},\mathcal{L}_{c})$ is a family of $T$-equivariant
test-configurations for $(M,-K_{M})$ with parameter $c\ll1$. It
can be seen as a perturbation of the trivial product test-configuration.

We need to compute $D_{Z}^{NA}(\mathcal{X},\mathcal{L}_{c})$. The
$D^{NA}$-part is easy. Since $K_{\mathcal{X}/\mathbb{P}^{1}}=p^{*}K_{M}+nP$,
let $B=(n-c)P$, we have $\mathcal{L}_{c}+K_{\mathcal{X}/\mathbb{P}^{1}}=\mathcal{O}(B)$.
Note $\mathcal{X}_{0}=\hat{M}+P$, and $\hat{M}$ is the blowup of
$M$ at $z_{*}$. The lct-part of $D^{NA}$ is
\[
\textrm{lct}(\mathcal{X},-B;\mathcal{X}_{0})=\sup\{t\in\mathbb{R}\mid(\mathcal{X},-B+t\mathcal{X}_{0})\ \textrm{is\ log\ canonical}\}.
\]
Since $\mathcal{X}$ is smooth and $-B+t\mathcal{X}_{0}=t\hat{M}+(c-n+t)P$
is simple normal crossing, thus we have
\[
\textrm{lct}(\mathcal{X},-B;\mathcal{X}_{0})=\min\{1,n-c+1\}=1,\ \textrm{when}\ c\ll1.
\]
On the other hand, $\mathcal{L}_{c}^{n+1}=\left(p^{*}K_{M}^{-1}-cP\right)^{n+1}=-c^{n+1}$.
Put these into the definition of $D^{NA}$ (\ref{eq: def of Ding inv}),
we have
\begin{equation}
D^{NA}(\mathcal{X},\mathcal{L}_{c})=\frac{c^{n+1}}{(n+1)c_{1}(M)^{n}}.\label{eq:ding inv c}
\end{equation}
Next we analyze the inner product part of $D_{Z}^{NA}(\mathcal{X},\mathcal{L}_{c})$
by the integral formula (\ref{eq: integral formula}) for inner products.

If we have a form $\Omega\in2\pi c_{1}(\mathcal{L}_{c})$ and function
$\Theta$ on $\mathcal{X}$ satisfying $\iota_{Z}\Omega=i\bar{\partial}\Theta$
and $\int_{\mathcal{X}_{1}}\Theta\Omega^{n}=0$, the integral formula
tells us
\[
\left\langle \alpha,\beta_{Z}\right\rangle =\frac{1}{(n+1)2\pi V}\int_{\mathcal{X}}\Theta\Omega^{n+1}.
\]
To construct such $\Omega$ and $\Theta$, we start with $\omega,\ \theta$.
First pulling back by $p$, we have $\iota_{Z}p^{*}\omega=i\bar{\partial}p^{*}\theta$
and $\int_{\mathcal{X}_{1}}p^{*}\theta\cdot\left(p^{*}\omega\right)^{n}=0$.

Let $s$ be a global section of $\mathcal{O}_{\mathcal{X}}(P)$ such
that $(s)=P$. Take a $\beta_{Z}(\mathbb{S}^{1})$-invariant metric
$h$ on $\mathcal{O}(P)$ such that $h(s)\equiv1$ outside of a neighborhood
of $P$ such as $U_{P}\coloneqq\pi^{-1}\{\left\vert \tau\right\vert <\frac{1}{2}\}$.
Denote the curvature form by
\[
-\eta\coloneqq-i\partial\bar{\partial}\log h\in2\pi c_{1}(\mathcal{O}(P)),\ \textrm{supp}\eta\subset U_{P}.
\]
The linearization of $\beta_{Z}$ induces a Hamiltonian function $f$
satisfying $\iota_{Z}\eta=i\bar{\partial}f$. We can write down $f$
on $\mathcal{X}\backslash P$. Since $\eta=i\partial\bar{\partial}\log h(s)$
on $\mathcal{X}\backslash P$, taking $\iota_{Z}$, we have $i\bar{\partial}f=i\bar{\partial}\left(Z\left(\log h(s)\right)\right)$.
It is easy to see $\log h(s)$ is $\beta_{Z}(\mathbb{S}^{1})$-invariant,
thus $Z\left(\log h(s)\right)$ is real. It follows that
\[
f=Z\left(\log h(s)\right)+\textrm{const},\ \textrm{on}\ \mathcal{X}\backslash P.
\]
Then we modify $f$ by a constant such that $f=Z\left(\log h(s)\right)$
on $\mathcal{X}\backslash P$. Since our choice of $h$, $f\equiv0$
outside $U_{P}$.

Now we let $\Omega\coloneqq p^{*}\omega+c\eta\in2\pi c_{1}(\mathcal{L}_{c})$
and $\Theta\coloneqq p^{*}\theta+cf$, they satisfy
\[
\iota_{Z}\Omega=i\bar{\partial}\Theta,\ \int_{\mathcal{X}_{1}}\Theta\Omega^{n}=0.
\]
The second identity since on $\mathcal{X}_{1}$ we have $f\equiv0$,
$\eta\equiv0$ and $\int_{M}\theta\omega^{n}=0$.

Putting these into integral formula (\ref{eq: integral formula}),
and rewrite the integrand in terms of equivariant curvatures, e.g.
$\Theta+\Omega$, $f+\eta$. Then expand it in $c$, we obtain

\begin{eqnarray*}
\left\langle \alpha,\beta_{Z}\right\rangle  & = & \frac{1}{(n+1)2\pi V}\int_{\mathcal{X}}\Theta\Omega^{n+1}=\frac{1}{(n+1)(n+2)2\pi V}\int_{\mathcal{X}}\left(\Theta+\Omega\right)^{n+2}\\
 & = & \frac{1}{(n+1)(n+2)2\pi V}\sum_{i=0}^{n+2}\dbinom{n+2}{i}\int_{\mathcal{X}}(p^{*}\theta+p^{*}\omega)^{n+2-i}\wedge(f+\eta)^{i}\cdot c^{i}.
\end{eqnarray*}
Denote $I_{i}\coloneqq\int_{\mathcal{X}}(p^{*}\theta+p^{*}\omega)^{n+2-i}\wedge(f+\eta)^{i}$.
Obviously, the coefficients $I_{i}$ (of $c^{i}$) are independent
of the choices of $\omega$ and $\eta$. We need to compute them.
Note that $p^{*}\omega\vert_{P}\equiv0$ and $p^{*}\theta\vert_{P}\equiv\theta(z_{*})=\vartheta(M)$,
if we can localize these integrals $I_{i}$ onto $P$, the computation
would be much simpler.

\begin{claim}\label{mclaim} For $0\leq i\leq n$, we have $I_{i}=0$,
and $I_{n+1}=-(2\pi)^{n+1}\vartheta(M).$ \end{claim}

We verify this claim by localization techniques. Firstly, we see
\[
I_{0}=(n+2)\int_{\mathcal{X}}p^{*}\theta\left(p^{*}\omega\right)^{n+1}=0
\]
 since $\omega^{n+1}=0$. So in the following, we assume $1\leq i\leq n+1$.

Let $U_{\delta}=\{x\in\mathcal{X}\mid h\left(s(x)\right)>\delta^{2}\}$
for small $\delta>0$, which is a tubular neighborhood of $P$. Since
\[
f+\eta=Z\left(\log h(s)\right)+i\partial\bar{\partial}\log h(s)=-(i\bar{\partial}-\iota_{Z})\partial\log h(s),\ \textrm{on}\ \mathcal{X}\backslash P,
\]
so in the following second row, we can replace one factor ``$f+\eta$''
by $-(i\bar{\partial}-\iota_{Z})\partial\log h(s)$.
\begin{eqnarray*}
I_{i} & = & \lim_{\delta\rightarrow0}\int_{U_{\delta}}(p^{*}\theta+p^{*}\omega)^{n+2-i}\wedge(f+\eta)^{i}\\
 & = & -\lim_{\delta\rightarrow0}\int_{U_{\delta}}(p^{*}\theta+p^{*}\omega)^{n+2-i}\wedge(f+\eta)^{i-1}\wedge(i\bar{\partial}-\iota_{Z})\partial\log h(s)\\
 & = & -\lim_{\delta\rightarrow0}\int_{\partial U_{\delta}}(p^{*}\theta+p^{*}\omega)^{n+2-i}\wedge(f+\eta)^{i-1}\wedge i\partial\log h(s).
\end{eqnarray*}
The third row is by integrating by parts (take $i\bar{\partial}-\iota_{Z}$
as a whole) and use
\[
(i\bar{\partial}-\iota_{Z})(p^{*}\theta+p^{*}\omega)=0,\ (i\bar{\partial}-\iota_{Z})(f+\eta)=0.
\]
We continue to expand it, when $1<i\leq n+1$, we have
\begin{eqnarray}
I_{i} & = & -\lim_{\delta\rightarrow0}(i-1)\int_{\partial U_{\delta}}f\cdot p^{*}\omega^{n+2-i}\wedge\eta^{i-2}\wedge i\partial\log h(s)\nonumber \\
 &  & -\lim_{\delta\rightarrow0}(n+2-i)\int_{\partial U_{\delta}}p^{*}\theta\cdot p^{*}\omega^{n+1-i}\wedge\eta^{i-1}\wedge i\partial\log h(s);\label{eq: I_i formula}
\end{eqnarray}
when $i=1$, we have
\[
I_{1}=-\lim_{\delta\rightarrow0}(n+1)\int_{\partial U_{\delta}}p^{*}\theta\cdot p^{*}\omega^{n}\wedge i\partial\log h(s).
\]
These limits will be integrals on $P$. The below lemma can be found
in standard textbooks in complex geometry, e.g. \cite{Huybrechts}
Page 203. We repeat the proof for completeness.

\begin{lemma}[Localization] With the above setting and notations,
let $\alpha$ be a smooth $2n$-form, we have
\begin{equation}
\lim_{\delta\rightarrow0}\int_{\partial U_{\delta}}\alpha\wedge i\partial\log h(s)=2\pi\int_{P}\alpha.\label{eq: localization}
\end{equation}
\end{lemma}
\begin{proof}
Essentially, this is a local result. Take a local coordinates $\{z_{0},\cdots,z_{n}\}$
and a local frame $e$ of $\mathcal{O}(P)$ around a point on $P$
such that $s=z_{0}\cdot e$. Hence $h(s)=\left\vert z_{0}\right\vert ^{2}u$,
$u\coloneqq h(e)$ is a positive smooth function. The boundary is
\[
\partial U_{\delta}=\{(z_{0},z')\mid\left\vert z_{0}\right\vert \cdot u(z_{0},z')^{1/2}=\delta\},\ z'=(z_{1},\cdots,z_{n}).
\]
We compute
\[
\partial\log h(s)=\frac{dz_{0}}{z_{0}}+\partial\log u.
\]
Note that $\lim_{\delta\rightarrow0}\int_{\partial U_{\delta}}\alpha\wedge i\partial\log u=0$
since $\log u$ is smooth. Thus we only consider $\int_{\partial U_{\delta}}\alpha\wedge\frac{dz_{0}}{z_{0}}$.
Denote $dz'=dz_{1}\wedge\cdots\wedge dz_{n}$, then $\alpha$ can
be decomposed as
\[
\alpha=g(z)dz'\wedge d\bar{z}'+dz_{0}\wedge\beta+d\bar{z}_{0}\wedge\gamma,
\]
where $g$ is a function, $\beta$ and $\gamma$ are forms. It can
be seen that only the first part has contribution to integral. Hence
we have
\begin{eqnarray*}
\lim_{\delta\rightarrow0}\int_{\partial U_{\delta}}\alpha\wedge i\partial\log h(s) & = & \int_{P}\left(\lim_{\delta\rightarrow0}\int g(z_{0},z')\frac{idz_{0}}{z_{0}}\right)dz'\wedge d\bar{z}'\\
 & = & 2\pi\int_{P}g(0,z')dz'\wedge d\bar{z}'=2\pi\int_{P}\alpha.
\end{eqnarray*}
\end{proof}
Now we apply (\ref{eq: localization}) to (\ref{eq: I_i formula}),
note $p^{*}\omega\vert_{P}\equiv0$ and $p^{*}\theta\vert_{P}\equiv\theta(z_{*})=\vartheta(M)$,
\[
\int_{P}p^{*}\theta\cdot\eta^{n}=\vartheta(M)\cdot\int_{P}\eta^{n}=-\vartheta(M)\cdot(2\pi)^{n}(-P)^{.n+1}=(2\pi)^{n}\vartheta(M),
\]
then Claim \ref{mclaim} follows.

We obtain an expansion
\[
\left\langle \alpha,\beta_{Z}\right\rangle =-\frac{\vartheta(M)}{(n+1)c_{1}(M)^{n}}c^{n+1}+A\cdot c^{n+2},
\]
where $A$ is a constant. Combining this with (\ref{eq:ding inv c}),
we obtain
\begin{equation}
D_{Z}^{NA}(\mathcal{X},\mathcal{L}_{c})=\frac{1-\vartheta(M)}{(n+1)c_{1}(M)^{n}}c^{n+1}+A\cdot c^{n+2}.\label{eq: expansion of Ding invari}
\end{equation}
If $\vartheta(M)>1$, $D_{Z}^{NA}(\mathcal{X},\mathcal{L}_{c})$ will
be negative when $c\ll1$, contradicting with relative D-semistability.
Proof of Theorem \ref{thm: semistable to theta 1} is completed. \end{proof}

\section{Uniformly relative D-stability and reduced NA J-functionals \label{sec: Unif stability}}

To define the uniform stability, we need a ``norm'' to measure how
far an equivariant test-configuration from a product. That is the
reduced non-Archimedean J-functional introduced by Hisamoto in \cite{Hisa toric}.
We restrict to the case of torus action.

\subsection{Twisting a $\mathbb{C}^{*}$-equivariant test-configuration}

Let $(M,L)$ be a polarized manifold. Let $T\subset\textrm{Aut}(M,L)$
be a $m$-dimensional torus whose elements are denoted by $\sigma=(\sigma_{1},\cdots,\sigma_{m})$,
$\sigma_{i}\in\mathbb{C}^{*}$. All the 1-parameter subgroups of $T$
constitute a lattice $\mathbb{Z}^{m}$. Specifically, each $\rho=(\rho_{1},\cdots,\rho_{m})\in\mathbb{Z}^{m}$
gives a subgroup
\[
\rho(\tau)\coloneqq(\tau^{\rho_{1}},\cdots,\tau^{\rho_{m}}):\mathbb{C}^{*}\rightarrow T.
\]
All the characters of $T$ constitute the dual lattice, each $\nu=(\nu_{1},\cdots,\nu_{m})\in\mathbb{Z}^{m}$
gives a character $\nu(\sigma)\coloneqq\sigma^{\nu}=\sigma_{1}^{\nu_{1}}\cdots\sigma_{m}^{\nu_{m}}$.
For any $\rho,\ \nu\in\mathbb{Z}^{m}$, the weight of composition
$\mathbb{C}^{*}\overset{\rho}{\rightarrow}T\overset{\nu}{\rightarrow}\mathbb{C}^{*}$
is $\left\langle \rho,\nu\right\rangle =\sum_{i=1}^{m}\rho_{i}\nu_{i}$.

\begin{definition}[twist a $T$-equivariant test-configuration] \label{def: twist test-config}Let
$(\mathcal{X},\mathcal{L})$ be a $T$-equivariant test-configuration
for $(M,L)$ with structure action $\alpha$ and fiberwise action
$\beta$. Given a 1-parameter subgroup $\rho:\mathbb{C}^{*}\rightarrow T$,
we modify the $\mathbb{C}^{*}$-action $\alpha$ on $(\mathcal{X},\mathcal{L})\vert_{\mathbb{C}}$
by
\[
\alpha_{\rho}(\tau)\coloneqq\alpha(\tau)\circ\beta(\rho(\tau)),\ \tau\in\mathbb{C}^{*}
\]
and keep the other data unchanged. We denote this new test-configuration
by $(\mathcal{X},\mathcal{L})\vert_{\mathbb{C}}^{\rho}$ (as a family
over $\mathbb{C}$), and call it the \textit{$\rho$-twisting} of
$(\mathcal{X},\mathcal{L})$.\end{definition}

\begin{remark}\label{Rm: prevent use inter formu}This definition
is same to \cite{Hisa toric}, merely using a different formulation.
Note the total space of compactification of $(\mathcal{X},\mathcal{L})\vert_{\mathbb{C}}^{\rho}$
will be different from $\mathcal{X}$, since action $\alpha$ on the
$\infty$-fiber of compactification space is always trivial. Thus
if $\mathcal{X}$ dominates the trivial product $M\times\mathbb{P}^{1}$
($\alpha$ only acts on $\mathbb{P}^{1}$), then its twisting no longer
dominates trivial product. This makes it difficult to apply the intersection
number formula (\ref{eq: J intersection}) of $J^{NA}$ for twisted
test-configurations. \end{remark}

Next we consider the effect of twisting on the associated filtration
$\mathcal{F}(\mathcal{X},\mathcal{L})$. By its definition (\ref{eq: filtration ind by t.c.}),
we see $F^{\mu}\textrm{H}^{0}(kL)$ is preserved by $T$-action on
$R(M,L)$. We call $\mathcal{F}(\mathcal{X},\mathcal{L})$ is a $T$-invariant
filtration. Consider the weight decomposition
\[
F^{\mu}\textrm{H}^{0}(kL)=\bigoplus_{\nu\in\mathbb{Z}^{m}}F^{\mu}\textrm{H}^{0}(kL)_{\nu},
\]
where $\nu$ runs over all characters of $T$. Let $\textrm{H}^{0}(kL)_{\nu}$
be the weight-$\nu$ subspace of $\textrm{H}^{0}(kL)$, then
\begin{equation}
F^{\mu}\textrm{H}^{0}(kL)_{\nu}\coloneqq\{s\in\textrm{H}^{0}(kL)_{\nu}\mid\tau^{-\mu}\cdot\bar{s}\in\textrm{H}^{0}(\mathcal{X}\vert_{\mathbb{C}},k\mathcal{L})\}.\label{eq: induce filtra beta}
\end{equation}
We will denote a $T$-invariant filtration by $\{F^{\mu}\textrm{H}^{0}(kL)_{\nu}\}$.
Note that for fixed $k,\ \nu$, it constitutes a filtration of $\textrm{H}^{0}(kL)_{\nu}$.

\begin{proposition}Let $(\mathcal{X},\mathcal{L})$ be a $T$-equivariant
ample test-configuration with associated filtration $\mathcal{F}(\mathcal{X},\mathcal{L})=\{F^{\mu}\textrm{H}^{0}(kL)_{\nu}\}$.
For 1-parameter subgroup $\rho\in\mathbb{Z}^{m}$, denote by $\mathcal{F}(\mathcal{X},\mathcal{L})^{\rho}=\{F_{\rho}^{\mu}\textrm{H}^{0}(kL)_{\nu}\}$
the filtration associated to twisting test-configuration $(\mathcal{X},\mathcal{L})\vert_{\mathbb{C}}^{\rho}$.
Then we have

\[
F_{\rho}^{\mu}\textrm{H}^{0}(kL)_{\nu}=F^{\mu-\left\langle \rho,\nu\right\rangle }\textrm{H}^{0}(kL)_{\nu}.
\]
\end{proposition}
\begin{proof}
Let $s\in\textrm{H}^{0}(kL)_{\nu}$, it satisfies $\sigma.s(\sigma^{-1}.z)=\sigma^{\nu}s(z)$
for $\forall\sigma\in T$, $z\in M$. Its equivariant extension by
the new action $\alpha_{\rho}$ is
\begin{eqnarray*}
\bar{s}^{\rho}(x)\coloneqq\alpha_{\rho}(\tau)s\left(\alpha_{\rho}(\tau^{-1})x\right) & = & \alpha(\tau)\beta(\rho(\tau))s\left(\beta(\rho(\tau^{-1}))\alpha(\tau^{-1})x\right)\\
 & = & \tau^{\left\langle \rho,\nu\right\rangle }\cdot\alpha(\tau)s\left(\alpha(\tau^{-1})x\right)=\tau^{\left\langle \rho,\nu\right\rangle }\cdot\bar{s}(x),
\end{eqnarray*}
where $x\in\mathcal{X}_{\tau}$, $\alpha_{\rho}(\tau^{-1})x\in\mathcal{X}_{1}\cong M$,
the restriction of $\beta$ on $(\mathcal{X}_{1},\mathcal{L}_{1})$
is the given $T$-action. Then the formula follows by (\ref{eq: induce filtra beta}).
\end{proof}
By taking base change, twisting operation can be extended to all rational
$\rho$. For irrational $\rho\in\mathbb{R}^{m}$, we can define $\rho$-twisting
of the associated filtration by the formula in above proposition.

\begin{definition}[twist the associated filtration]Let $(\mathcal{X},\mathcal{L})$
be a $T$-equivariant ample test-configuration. Its associated $T$-invariant
filtration is $\mathcal{F}(\mathcal{X},\mathcal{L})=\{F^{t}\textrm{H}^{0}(kL)_{\nu}\}$.
For any $\rho\in\mathbb{R}^{m}$, we define a new $T$-invariant filtration
(called the \textit{$\rho$-twisting} of $\mathcal{F}(\mathcal{X},\mathcal{L})$)
denoted by $\mathcal{F}(\mathcal{X},\mathcal{L})^{\rho}=\{F_{\rho}^{t}\textrm{H}^{0}(kL)_{\nu}\}$,
\begin{equation}
F_{\rho}^{t}\textrm{H}^{0}(kL)_{\nu}\coloneqq F^{t-\left\langle \rho,\nu\right\rangle }\textrm{H}^{0}(kL)_{\nu},\ F_{\rho}^{t}\textrm{H}^{0}(kL)\coloneqq\bigoplus_{\nu\in\mathbb{Z}^{m}}F_{\rho}^{t}\textrm{H}^{0}(kL)_{\nu}.\label{eq: twist filtra formula}
\end{equation}
Recall that $F^{t}\textrm{H}^{0}(kL)=F^{\left\lceil t\right\rceil }\textrm{H}^{0}(kL)$
for all $t\in\mathbb{R}$. \end{definition}

\subsection{Uniformly relative D-stability}

For a probability measure $\mathfrak{m}$ on $\mathbb{R}$ with bounded
support, we set
\[
j(\mathfrak{m})\coloneqq\sup\textrm{supp}\mathfrak{m}-\int_{\mathbb{R}}\lambda d\mathfrak{m}\geq0.
\]

\begin{definition}[reduced NA J-functionals]\label{def of J^NA_T}

(1) Let $\mathcal{F}$ be an admissible filtration of the section
ring $R(M,L)$. We define
\begin{equation}
J^{NA}(\mathcal{F})\coloneqq j\left(\textrm{LM}(\mathcal{F})\right),\label{eq: JNA on filtration}
\end{equation}
where $\textrm{LM}(\mathcal{F})$ is the limit measure of $\mathcal{F}$,
see Theorem \ref{thm: LM is deriv of Vol}. In particular, take $\mathcal{F}$
to be the associated filtration $\mathcal{F}(\mathcal{X},\mathcal{L})$
to an ample test-configuration $(\mathcal{X},\mathcal{L})$ for $(M,L)$,
it gives $J^{NA}(\mathcal{X},\mathcal{L})$ (\ref{eq: NA J}).

(2) Let $(\mathcal{X},\mathcal{L})$ be a $T$-equivariant ample test-configuration
for $(M,L)$. The \textit{reduced} (means modulo $T$-action) NA J-functional
is defined by
\[
J_{T}^{NA}(\mathcal{X},\mathcal{L})\coloneqq\inf_{\rho\in\mathbb{R}^{m}}J^{NA}\left(\mathcal{F}(\mathcal{X},\mathcal{L})^{\rho}\right),
\]
where $\mathcal{F}(\mathcal{X},\mathcal{L})^{\rho}$ is $\rho$-twisting
of $\mathcal{F}(\mathcal{X},\mathcal{L})$ (\ref{eq: twist filtra formula}).
\end{definition}

By \cite{BHJ1} Corollary B, $J^{NA}(\mathcal{X},\mathcal{L})$ is
zero if and only if $(\mathcal{X},\mathcal{L})$ is trivial product
$M\times\mathbb{P}^{1}$ (i.e. action $\alpha$ is trivial). Roughly
speaking, $J_{T}^{NA}$ measures how far is $(\mathcal{X},\mathcal{L})$
from product test-configurations (may have a nontrivial action $\alpha$
via $T$-action).

\begin{definition} \label{def of Uni D-stable}With the same assumptions
in Definition \ref{def: def of D-semistable}, here we take $G$ is
torus $T$. $M$ is said to be \textit{uniformly D-stable relative
to} $T$ if there exists a $\delta>0$ such that
\begin{equation}
D_{Z}^{NA}(\mathcal{X},\mathcal{L})\geq\delta\cdot J_{T}^{NA}(\mathcal{X},\mathcal{L})\label{eq: def of unif Dstable}
\end{equation}
for all $T$-equivariant ample test-configuration for $(M,-K_{M})$.
\end{definition}

We want to show uniformly relative D-stability implies $\vartheta(M)<1$.
Continue the proof for semistability implies $\vartheta(M)\leq1$,
next we need to expand $J_{T}^{NA}(\mathcal{X},\mathcal{L}_{c})$
in parameter $c$. As noted in Remark \ref{Rm: prevent use inter formu},
the intersection number formula of $J^{NA}$ is difficult to use when
we make twisting. Hence in the sequel, we give a convex-geometry description
for $J_{T}^{NA}$ via Okounkov bodies. Finally, we will see $J_{T}^{NA}(\mathcal{X},\mathcal{L}_{c})$
is attained by the trivial twisting (i.e. $\rho=0$) when $c\ll1$,
then we are done.

\subsection{Infinitesimal Okounkov bodies}

In \cite{Okounkov}, Okounkov associated a convex body $\triangle\subset\mathbb{R}^{n}$
to a polarized algebraic variety $(X,L)$ (only need $L$ is big),
its Euclidean volume gives the volume of $L$. Boucksom-Chen \cite{Chen huayi}
showed an admissible filtration $\mathcal{\mathcal{F}}$ of $R(X,L)$
gives rise to a concave function $G[\mathcal{F}]$ on $\triangle$,
called the concave transform of $\mathcal{\mathcal{F}}$, and the
pushforward of Lebesgue measure by $G[\mathcal{F}]$ gives the limit
measure $\textrm{LM}(\mathcal{\mathcal{F}})$. In particular, take
$\mathcal{\mathcal{F}}=\mathcal{\mathcal{F}}(\mathcal{X},\mathcal{L})$
for a test-configuration $(\mathcal{X},\mathcal{L})$, the associated
concave transform $G[\mathcal{X},\mathcal{L}]$ encodes many information
of $(\mathcal{X},\mathcal{L})$, such as $\textrm{DH}(\mathcal{X},\mathcal{L})$,
$J^{NA}(\mathcal{X},\mathcal{L})$.

In Okounkov's construction, we need to choose a flag of subvarieties
in $M$, the resulting $\triangle$ depends heavily on the choices.
In our setting, we need a flag adapting to the $T$-action. For this,
a local version called \textit{infinitesimal Okounkov body} introduced
in \cite{L-Mustata} is more suitable for our setting. Its construction
is as following.

Let $(M,L)$ be a polarized manifold with a lifted $m$-dimensional
torus action $\beta:T\rightarrow\textrm{Aut}(M,L)$. By \cite{Atiyah},
there exists a $T$-fixed point $z_{*}\in M$. Let $V_{\bullet}$
be a full flag of subspaces in the tangent space $T_{z_{*}}M$ (called
an \textit{infinitesimal flag} at $z_{*}$). Let
\[
q:\hat{M}\rightarrow M
\]
be the blowup at $z_{*}$. Then $V_{\bullet}$ induces a flag of subvarieties
$Y_{\bullet}$ in $\hat{M}$ by taking projectivization of each subspace
in $V_{\bullet}$ (subvarieties in $Y_{\bullet}$ are contained in
the exceptional divisor). Then the infinitesimal Okounkov body is
the Okounkov body associated to datum $(\hat{M},q^{*}L,Y_{\bullet})$.
More specifically, note $q^{*}L$ is big, and the pullback
\[
q^{*}:\textrm{H}^{0}(M,kL)\rightarrow\textrm{H}^{0}(\hat{M},kq^{*}L)
\]
is an isomorphism ($q_{*}\mathcal{O}_{\hat{M}}=\mathcal{O}_{M}$).
For $s\in\textrm{H}^{0}(M,kL)$, the vanishing orders of $q^{*}s$
along $Y_{\bullet}$ give us a rank-$n$ valuation $\mathcal{V}(s)$,
then $\triangle$ will be defined by these values.

Next we go through the construction in details. First we choose an
infinitesimal flag $V_{\bullet}$ adapted to $T$-action.

\textbf{Step 1}: take a $\beta$-invariant local frame $e$ around
$z_{*}$.

Suppose $\beta$ acts on $L_{z_{*}}$ via character $\nu_{0}^{*}\in\mathbb{Z}^{m}$.
Since $L$ is ample, there exists $\kappa_{0}\in\mathbb{N}$ and $s\in\textrm{H}^{0}(M,\kappa_{0}L)$
such that $s(z_{*})\neq0$. Decompose $s=\sum_{\nu}s_{\nu}$ w.r.t.
action $\beta$, where $s_{\nu}\in\textrm{H}^{0}(M,\kappa_{0}L)_{\nu}$
satisfying
\[
\sigma.s_{\nu}(\sigma^{-1}.z)=\sigma^{\nu}\cdot s_{\nu}(z),\ \textrm{for}\ \sigma\in T,\ z\in M.
\]
Take $z=z_{*}$ we see $s_{\nu}(z_{*})\neq0$ only if $\nu=\nu_{0}^{*}$.
Now we set
\[
e\coloneqq s_{\nu_{0}^{*}},
\]
which is a $\beta$-invariant local frame of $L^{\kappa_{0}}$ on
the $\beta$-invariant open set $\{s_{\nu_{0}^{*}}\neq0\}\ni z_{*}$.

\textbf{Step 2}: take a local coordinates $(z_{i})_{i=1}^{n}$ centered
at $z_{*}$ (i.e. $z_{i}(z_{*})=0$) such that action $\beta$ is
diagonal under $(z_{i})$ up to higher order terms.

which diagonalize s action $\beta$ at first order.

Since $\beta$ induces a $T$-action on $T_{z_{*}}^{(1,0)}M$, let
$\nu_{1}^{*},\cdots,\nu_{n}^{*}\in\mathbb{Z}^{m}$ be the associated
characters. By a linear transformation, we can find a coordinates
$(z_{i})$ centered at $z_{*}$ such that
\begin{equation}
z_{i}\left(\beta(\sigma).z\right)=\nu_{i}^{*}(\sigma)\cdot z_{i}+\textrm{higher\ order terms},\ \textrm{for}\ \forall\sigma\in T,\ 1\leq i\leq n.\label{eq: lineariza of beta action}
\end{equation}
Where ``higher order terms'' means monomials of $(z_{i})$ with
order $\geq2$.

\textbf{Step 3}: choose an infinitesimal flag $V_{\bullet}$.

We describe the blowup $q:\hat{M}\rightarrow M$ explicitly by coordinates.
Let $E=\mathbb{P}(T_{z_{*}}M)$ be the exceptional divisor. Near $E$,
using coordinate $(z_{i})$, $\hat{M}$ is given by
\[
\left\{ \left((z_{i}),[w_{1},\cdots,w_{n}]\right)\in M\times\mathbb{P}^{n-1}\mid z_{i}w_{j}=z_{j}w_{i},\ \forall i,j\right\} .
\]
On the open set $\{w_{1}\neq0\}$, we define $u_{1}=z_{1}$, $u_{i}=w_{i}/w_{1}$
for $2\leq i\leq n$, then $(u_{i})_{i=1}^{n}$ constitutes a local
coordinates on $\hat{M}$ around $E$. Under this coordinates, $q$
is given by
\begin{equation}
q:(u_{1},\cdots,u_{n})\longmapsto(z_{i})=(u_{1},u_{1}u_{2},\cdots,u_{1}u_{n}).\label{eq: blowup}
\end{equation}
Now we choose a full flag of subspaces in $T_{z_{*}}M$,
\begin{equation}
V_{\bullet}:\ T_{z_{*}}M\supset\textrm{span}\{\partial_{1},\partial_{3},\cdots,\partial_{n}\}\supset\textrm{span}\{\partial_{1},\partial_{4},\cdots,\partial_{n}\}\supset\cdots\supset\textrm{span}\{\partial_{1}\},\label{eq: infi flag}
\end{equation}
where $\partial_{i}\coloneqq\frac{\partial}{\partial z_{i}}$. Taking
projectivization, it induces a flag of subvarieties in $\hat{M}$,
\begin{eqnarray}
Y_{\bullet}: &  & Y_{1}=E=\overline{\{u_{1}=0\}}\supset Y_{2}=\overline{\{u_{1}=u_{2}=0\}}\supset\cdots\label{eq: flag}\\
 &  & \cdots\supset Y_{n}=\{u_{1}=\cdots=u_{n}=0\}.\nonumber
\end{eqnarray}

\textbf{Step 4}: define valuation $\mathcal{V}$.

Take a nonzero $s\in\textrm{H}^{0}(M,kL)$, then $s^{\kappa_{0}}\in\textrm{H}^{0}(M,\kappa_{0}kL)$
(see Step 1 for $\kappa_{0}$), express $s^{\kappa_{0}}=f(z)\cdot e^{k}$
in a neighborhood of $z_{*}$. We expand $f$ as power series,
\begin{equation}
f(z)=\sum_{a\geq0}c_{a}\cdot z^{a},\ a=(a_{1},\cdots,a_{n}),\ z^{a}=z_{1}^{a_{1}}\cdots z_{n}^{a_{n}}.\label{eq: expand of f}
\end{equation}
Pulling back by $q$, we have $q^{*}s^{\kappa_{0}}=(f\circ q)q^{*}e^{k}$.
By (\ref{eq: blowup}), we have
\[
(f\circ q)(u)=\sum_{a\geq0}c_{a}\cdot u_{1}^{\sum a_{i}}u_{2}^{a_{2}}\cdots u_{n}^{a_{n}}.
\]
The vanishing orders of $f\circ q$ along the subvarieties in flag
(\ref{eq: flag}) define a rank-$n$ valuation
\begin{equation}
\mathcal{V}(s)\coloneqq\frac{1}{\kappa_{0}}\min\{(\sum_{i=1}^{n}a_{i},a_{2}\cdots,a_{n})\mid c_{a}\neq0\}\in\frac{1}{\kappa_{0}}\mathbb{N}^{n},\ \mathcal{V}(0)\coloneqq\infty,\label{eq: define  valuation}
\end{equation}
where the minimizer is taken w.r.t. the lexicographic order. It satisfies
\[
\mathcal{V}(s\otimes t)=\mathcal{V}(s)+\mathcal{V}(t),\ \textrm{for}\ s\in\textrm{H}^{0}(kL),\ t\in\textrm{H}^{0}(lL);
\]
\[
\mathcal{V}(s+t)\geq\min\{\mathcal{V}(s),\mathcal{V}(t)\},\ \textrm{for}\ s,t\in\textrm{H}^{0}(kL),
\]
where $\geq$ is the lexicographic order.

\begin{remark}Note the first component of $\mathcal{V}(s)$ is the
multiplicity $\textrm{ord}_{z_{*}}(s)=\textrm{ord}_{E}(q^{*}s)$.
In the definition of $\mathcal{V}(s)$, the order we used takes the
multiplicity $\sum_{i=1}^{n}a_{i}$ as priority, contrasting with
the lexicographic order on $(a_{1},\cdots,a_{n})$. We choose this
order because it is compatible with expansion of power series, this
is required for Proposition \ref{pro: weight and valu} in the sequel.
This is also the reason why we chose infinitesimal flag $V_{\bullet}$
(\ref{eq: infi flag}). \end{remark}

\textbf{Step 5}: define the (infinitesimal) Okounkov body $\triangle$.

Collecting all the values of $\mathcal{V}(s)$, we obtain a semigroup
\[
\Gamma(L)\coloneqq\{(\mathcal{V}(s),k)\mid0\neq s\in\textrm{H}^{0}(M,kL),\ k\in\mathbb{N}\}\subset\frac{1}{\kappa_{0}}\mathbb{N}^{n+1}.
\]
Let $\Sigma(L)\subset\mathbb{R}^{n+1}$ be the convex cone spanned
by $\Gamma(L)$.

\begin{definition}[infinitesimal Okounkov body]\label{Def: infi Okounkov body}
The infinitesimal Okounkov body associated to the infinitesimal flag
$V_{\bullet}$ (\ref{eq: infi flag}) at the $T$-fixed point $z_{*}$
is defined by
\[
\triangle(L)\coloneqq\{x\mid(x,1)\in\Sigma(L)\}\subset\mathbb{R}^{n}.
\]
It is the ordinary Okounkov body of $q^{*}L$ associated to the flag
of subvarieties (\ref{eq: flag}) in $\hat{M}$. \end{definition}

\begin{remark} $\triangle(L)$ only depends on the infinitesimal
flag (\ref{eq: infi flag}), but this flag depends on the action $\beta$
around $z_{*}$. Hence $\triangle(L)$ depends on action $\beta$.
\end{remark}

On the other hand, $\triangle(L)$ also can be defined as closure
\[
\triangle(L)=\overline{\bigcup_{k\geq1}\triangle_{k}(L)},\ \triangle_{k}(L)\coloneqq\{\mathcal{V}(s)/k\mid0\neq s\in\textrm{H}^{0}(M,kL)\}.
\]
Okounkov body $\triangle(L)$ encodes the volume of $L$. By \cite{L-Mustata}
Theorem A, we have
\[
\left\vert \triangle(L)\right\vert =\frac{(q^{*}L)^{n}}{n!}=\frac{L^{n}}{n!},
\]
where $\left\vert \cdot\right\vert $ is the standard Lebesgue measure.

In the following, for an eigen-section $s\in\textrm{H}^{0}(M,kL)_{\nu}$,
we see its weight $\nu\in\mathbb{Z}^{m}$ can be linearly determined
by its valuation $\mathcal{V}(s)$.

\begin{proposition}[weight v.s. valuation]\label{pro: weight and valu}For
any weight $\nu\in\mathbb{Z}^{m}$ and nonzero $s\in\textrm{H}^{0}(M,kL)_{\nu}$.
Suppose $\mathcal{V}(s)/k=(x_{1},\cdots,x_{n})$, then we have
\begin{equation}
\frac{\mathbf{\nu}}{k}=-x_{1}\nu_{1}^{*}+\sum_{i=2}^{n}x_{i}(\nu_{1}^{*}-\nu_{i}^{*})+\frac{\nu_{0}^{*}}{\kappa_{0}}\eqqcolon H_{\beta}(x).\label{eq: weight and valu}
\end{equation}
See above Step 1-2 for constants $\kappa_{0}$, $\nu_{0}^{*},\cdots,\nu_{n}^{*}$,
which depend on action $\beta$ around $z_{*}$. In particular, the
weight $\nu$ of $s$ is an affine function $H_{\beta}$ of the valuation
$\mathcal{V}(s)$. \end{proposition}
\begin{proof}
By the local frame $e$ and local coordinates $(z_{i})$ near $z_{*}$
in Step 1-2, assume $s^{\kappa_{0}}=f(z)\cdot e^{k}$ and expand $f=\sum_{a\geq0}c_{a}\cdot z^{a}$
near $z_{*}$. Since $s$ and $e$ satisfy equivariant condition,
i.e.
\[
\sigma.s(\sigma^{-1}.z)=\sigma^{\nu}s(z),\ \sigma.e(\sigma^{-1}.z)=\sigma^{\nu_{0}^{*}}e(z),\ \textrm{for}\ \forall\sigma\in T.
\]
It implies
\begin{equation}
f(\sigma.z)=\sigma^{k\nu_{0}^{*}-\kappa_{0}\nu}f(z),\ \textrm{for}\ \forall\sigma\in T.\label{eq: inv property of f}
\end{equation}
Let $z^{a}=z_{1}^{a_{1}}\cdots z_{n}^{a_{n}}$ be the minimum monomial
in the expansion of $f$ w.r.t. the lexicographic order of $(\sum_{i=1}^{n}a_{i},a_{2}\cdots,a_{n})$.
By the definition (\ref{eq: define  valuation}) of $\mathcal{V}(s)$,
we have
\[
(x_{1},\cdots,x_{n})=\frac{1}{\kappa_{0}k}(\sum_{i=1}^{n}a_{i},a_{2}\cdots,a_{n}).
\]
Put the expansion (\ref{eq: lineariza of beta action}) of action
$\beta$ into (\ref{eq: inv property of f}), we obtain
\[
c_{a}\cdot\sigma^{a\cdot\nu^{*}}z^{a}+\cdots=c_{a}\cdot\sigma^{k\nu_{0}^{*}-\kappa_{0}\nu}z^{a}+\cdots,\ \forall\sigma\in T,
\]
where $a\cdot\nu^{*}\coloneqq\sum_{i=1}^{n}a_{i}\nu_{i}^{*}$. It
follows that $\sum_{i=1}^{n}a_{i}\nu_{i}^{*}=k\nu_{0}^{*}-\kappa_{0}\nu$.
Then replace $(a_{i})$ by $(x_{i})$, we obtain (\ref{eq: weight and valu}).
\end{proof}

\subsection{Concave transforms of filtrations}

It is showed in \cite{Chen huayi}, an admissible filtration of section
ring gives rise to a concave function on the Okounkov body. In our
definition, $\triangle(L)$ is the Okounkov body of $q^{*}L$. Since
\[
q^{*}:R(M,L)\rightarrow R(\hat{M},q^{*}L)
\]
is an isomorphism, an admissible filtration of $R(M,L)$ induces another
one for $R(\hat{M},q^{*}L)$, thus by \cite{Chen huayi} it gives
rise to a concave function on $\triangle(L)$. Next we review the
construction in \cite{Chen huayi}.

Let $\mathcal{F}=\{F^{t}\textrm{H}^{0}(kL)\}$ be an admissible filtration
of $R(M,L)$. Let
\[
\triangle_{k}^{t}(L,\mathcal{F})\coloneqq\{\mathcal{V}(s)/k\mid0\neq s\in F^{t}\textrm{H}^{0}(kL)\}\subset\triangle_{k}(L),
\]
then $\{\triangle_{k}^{t}(L,\mathcal{F})\}_{t}$ constitutes a filtration
of subsets of $\triangle_{k}(L)$, which is decreasing in $t$ and
becoming empty when $t>Ck$ for a constant $C>0$. We also have (see
\cite{L-Mustata} Lemma 1.4)
\[
\left\vert \triangle_{k}^{t}(L,\mathcal{F})\right\vert =\dim F^{t}\textrm{H}^{0}(kL).
\]
For each $k\geq1$, we define a function $G_{k}$ on $\triangle_{k}(L)$
(called $k$-th approximation function),
\[
G_{k}(x)\coloneqq\sup\{t/k\mid x\in\triangle_{k}^{t}(L,\mathcal{F})\},\ \textrm{for}\ x\in\triangle_{k}(L).
\]
With these $G_{k}$, we define a function on $\bigcup_{k\geq1}\triangle_{k}(L)$,
\begin{equation}
G[\mathcal{F}](x)\coloneqq\sup\{G_{k}(x)\mid k\ \textrm{s.t.}\ x\in\triangle_{k}(L)\},\ \textrm{for}\ x\in\bigcup_{k\geq1}\triangle_{k}(L).\label{eq: concave transform}
\end{equation}

\begin{theorem}[concave transform \cite{Chen huayi}]\label{thm: push-out measure okounkov}Let
$\mathcal{F}$ be an admissible filtration of the section ring $R(M,L)$.
The function $G[\mathcal{F}]$ defined above can be extended to a
bounded concave function on the interior of $\triangle(L)$ (called
the \textit{concave transform} of $\mathcal{F}$). Moreover, the pushforward
\[
G[\mathcal{F}]_{\#}\left(\frac{dx}{\left\vert \triangle(L)\right\vert }\right)
\]
of the normalized Lebesgue measure is the limit measure $\textrm{LM}(\mathcal{F})$.
\end{theorem}

By this theorem, we can use $G[\mathcal{F}]$ to express $J^{NA}$,
\begin{equation}
J^{NA}(\mathcal{F})=j\left(\textrm{LM}(\mathcal{F})\right)=\sup_{\triangle(L)}G[\mathcal{F}]-\int_{\triangle(L)}G[\mathcal{F}]\frac{dx}{\left\vert \triangle(L)\right\vert }.\label{eq: J^NA Okounkov}
\end{equation}

\subsection{\label{subsec: convex-geo descrip}A convex-geometry description
for $J_{T}^{NA}$}

Let $(\mathcal{X},\mathcal{L})$ be a $T$-equivariant ample test-configuration
for $(M,L)$. In order to express $J_{T}^{NA}(\mathcal{X},\mathcal{L})$
via (\ref{eq: J^NA Okounkov}), we consider the concave transform
of twisted filtration $\mathcal{F}(\mathcal{X},\mathcal{L})^{\rho}$.

For each $\nu\in\mathbb{Z}^{m}$, let
\[
\triangle_{k}(L)_{\nu}\coloneqq\{\mathcal{V}(s)/k\mid0\neq s\in\textrm{H}^{0}(M,kL)_{\nu}\}.
\]
Proposition \ref{pro: weight and valu} ensures $\triangle_{k}(L)_{\nu}$
are contained in different hyperplanes for different $\nu$, hence
we have disjoint union
\[
\triangle_{k}(L)=\bigsqcup_{\nu\in\mathbb{Z}^{m}}\triangle_{k}(L)_{\nu}.
\]
Let $\mathcal{F}(\mathcal{X},\mathcal{L})=\{F^{t}\textrm{H}^{0}(kL)_{\nu}\}$
be the associated $T$-invariant filtration to $(\mathcal{X},\mathcal{L})$.
Similarly, for each $\nu\in\mathbb{Z}^{m}$, let
\[
\triangle_{k}^{t}(L,\mathcal{F})_{\nu}\coloneqq\{\mathcal{V}(s)/k\mid0\neq s\in F^{t}\textrm{H}^{0}(kL)_{\nu}\}\subseteq\triangle_{k}(L)_{\nu}.
\]
Then $\{\triangle_{k}^{t}(L,\mathcal{F})_{\nu}\}_{t}$ constitutes
a filtration of subsets of $\triangle_{k}(L)_{\nu}$ and we have disjoint
union
\[
\triangle_{k}^{t}(L,\mathcal{F})=\bigsqcup_{\nu\in\mathbb{Z}^{m}}\triangle_{k}^{t}(L,\mathcal{F})_{\nu}.
\]
By the definition of $k$-th approximation function $G_{k}$, for
$x\in\triangle_{k}(L)_{\nu}$, we have
\[
G_{k}(x)=\sup\{t/k\mid x\in\triangle_{k}^{t}(L,\mathcal{F})_{\nu}\}.
\]

\begin{proof}[Proof of Theorem \ref{Intro convex-geometry of JNA}]Denote
by $G_{k}^{\rho}$ be the $k$-th approximation function associated
to twisted filtration $\mathcal{F}^{\rho}\coloneqq\mathcal{F}(\mathcal{X},\mathcal{L})^{\rho}$,
i.e. for $x\in\triangle_{k}(L)_{\nu}$, we define
\[
G_{k}^{\rho}(x)=\sup\{t/k\mid x\in\triangle_{k}^{t}(L,\mathcal{F}^{\rho})_{\nu}\}.
\]
By the definition of $\rho$-twisting (\ref{eq: twist filtra formula}),
we have
\[
G_{k}^{\rho}(x)=\sup\{s/k\mid x\in\triangle_{k}^{s}(L,\mathcal{F})_{\nu}\}+\frac{\left\langle \rho,\nu\right\rangle }{k}=G_{k}(x)+\left\langle \rho,\frac{\nu}{k}\right\rangle .
\]
We can express $\left\langle \rho,\frac{\nu}{k}\right\rangle $ by
$x$. In fact, since $x\in\triangle_{k}(L)_{\nu}$, suppose $x=\mathcal{V}(s)/k$
for a nonzero $s\in\textrm{H}^{0}(M,kL)_{\nu}$. By the relation (\ref{eq: weight and valu}),
we have
\begin{equation}
\left\langle \rho,\frac{\nu}{k}\right\rangle =\left\langle \rho,H_{\beta}(x)\right\rangle =-\left\langle \rho,\nu_{1}^{*}\right\rangle x_{1}+\sum_{i=2}^{n}\left\langle \rho,\nu_{1}^{*}-\nu_{i}^{*}\right\rangle x_{i}+\frac{1}{\kappa_{0}}\left\langle \rho,\nu_{0}^{*}\right\rangle .\label{eq: twist plus affine func}
\end{equation}
The RHS is an affine function of $x$ with coefficients independent
of $k$. By the definition of concave transforms, we obtain $G[\mathcal{F}^{\rho}]=G+\left\langle \rho,H_{\beta}\right\rangle $.
Put this into (\ref{eq: J^NA Okounkov}), we obtain (\ref{eq: Intro JNA formula}).
The part (2) of this theorem is obvious. \end{proof}

\begin{remark} By (\ref{eq: Intro JNA formula}), we see $J^{NA}\left(\mathcal{F}(\mathcal{X},\mathcal{L})^{\rho}\right)$
is proper as a function of $\rho\in\mathbb{R}^{m}$, thus the infimum
$J_{T}^{NA}$ can be attained by some $\rho_{0}\in\mathbb{R}^{m}$,
but $\rho_{0}$ may be irrational. \end{remark}

\section{Uniformly relative D-stability implies $\vartheta(M)<1$\label{sec: theta<1}}

Now we can complete the proof of our main theorem. The concave transform
associated to the deformation to normal cone had been considered in
\cite{Nystrom TC-body}.

\begin{theorem} \label{thm: unif stable theta 1} Let $M$ be a Fano
manifold, $L=-K_{M}$. Let $T\subset\textrm{Aut}^{0}(M)$ be a torus.
If $M$ is uniformly D-stable relative to $T$ in the sense of Definition
\ref{def of Uni D-stable}, then $\vartheta(M)\leq1-\delta$. \end{theorem}
\begin{proof}
We take the canonical lifting of $T$-action $\beta:T\rightarrow\textrm{Aut}(M,-K_{M})$.
Let $z_{*}$ be the $T$-fixed point provided by Proposition \ref{prop: exist fixed pt},
and $\triangle(L)$ is the infinitesimal Okounkov body associated
to $z_{*}$ and infinitesimal flag $V_{\bullet}$ (\ref{eq: infi flag}).
Let $\{(\mathcal{X},\mathcal{L}_{c})\mid\mathbb{Q}_{+}\ni c\ll1\}$
be the family of $T$-equivariant test-configurations constructed
in the proof of Theorem \ref{thm: semistable to theta 1}.

Let $\mathcal{F}_{c}=\{F_{c}^{t}\textrm{H}^{0}(kL)\}$ be the filtration
associated to $(\mathcal{X},\mathcal{L}_{c})$ (where subscript ``$c$''
does not mean twisting!). The Lemma 5.17 in \cite{BHJ1} tells us
\[
F_{c}^{t}\textrm{H}^{0}(kL)=\{s\in\textrm{H}^{0}(kL)\mid\textrm{ord}_{z_{*}}(s)\geq t+ck\},\ \textrm{when}\ t\leq0;
\]
and $F_{c}^{t}\textrm{H}^{0}(kL)=0$ when $t>0$. Let $q:\hat{M}\rightarrow M$
be the blowup at $z_{*}$ with exceptional divisor $E$, we have
\[
F_{c}^{t}\textrm{H}^{0}(kL)\cong\textrm{H}^{0}(\hat{M},kq^{*}L-\left\lceil t+ck\right\rceil E),\ \textrm{when}\ -ck\leq t\leq0;
\]
$F_{c}^{t}\textrm{H}^{0}(kL)\cong\textrm{H}^{0}(\hat{M},kq^{*}L)$
when $t<-ck$; and $F_{c}^{t}\textrm{H}^{0}(kL)=0$ when $t>0$.

From this we obtain the DH measure (see Proposition 8.5 \cite{BHJ1}),
\begin{equation}
\textrm{DH}(\mathcal{X},\mathcal{L}_{c})=\frac{n}{L^{n}}(\lambda+c)^{n-1}\boldsymbol{1}_{[-c,0]}d\lambda+(1-\frac{c^{n}}{L^{n}})\delta_{0},\label{eq: DH of Lc}
\end{equation}
where $\boldsymbol{1}_{[-c,0]}$ is characteristic function, $d\lambda$
is Lebesgue measure and $\delta_{0}$ is the Dirac measure at $0$.
It follows that
\[
J^{NA}(\mathcal{X},\mathcal{L}_{c})=j\left(\textrm{DH}(\mathcal{X},\mathcal{L}_{c})\right)=\frac{c^{n+1}}{(n+1)L^{n}}.
\]

Next we consider the concave transform of $\mathcal{F}_{c}$. Since
the first component of $\mathcal{V}(s)$ is $\textrm{ord}_{z_{*}}(s)=\textrm{ord}_{E}(q^{*}s)$,
thus when $t\leq0$ we have
\[
\triangle_{k}^{t}(L,\mathcal{F}_{c})=\{\mathcal{V}(s)/k\mid0\neq s\in F_{c}^{t}\textrm{H}^{0}(kL)\}=\triangle_{k}(L)\cap\{x_{1}\geq\frac{t}{k}+c\};
\]
and $\triangle_{k,t}(L,\mathcal{F}_{c})=\emptyset$ when $t>0$. It
follows that for $x\in\triangle_{k}(L)$,
\[
G_{k}(x)=\sup\{t/k\mid x\in\triangle_{k}^{t}(L,\mathcal{F}_{c})\}=\min\{x_{1}-c,0\}.
\]
Hence the concave transform is
\[
G[\mathcal{F}_{c}]=\min\{x_{1}-c,0\}.
\]

Next we consider $J_{T}^{NA}(\mathcal{X},\mathcal{L}_{c})$, denote
by $\mathcal{F}_{c}^{\rho}$ the $\rho$-twisting of $\mathcal{F}_{c}$.
We claim the infimum
\[
J_{T}^{NA}(\mathcal{X},\mathcal{L}_{c})\coloneqq\inf_{\rho\in\mathbb{R}^{m}}J^{NA}(\mathcal{F}_{c}^{\rho})
\]
is attained at $\rho=0$ when $c\ll1$. First we note the following
facts:

(1) $\triangle(L)$ is a convex body with nonempty interior, and by
the definition of $\mathcal{V}$ (\ref{eq: define  valuation}) it
is contained in region $\{x\in\mathbb{R}_{\geq}^{n}\mid x_{1}\geq x_{2}+\cdots+x_{n}\}.$

(2) $\inf_{\triangle}x_{1}=0$. If not, $\triangle\subset\{x_{1}\geq c\}$
for some $1\gg c\in\mathbb{Q}_{+}$, thus $G[\mathcal{F}_{c}]\equiv0$
on $\triangle$. However, since $\textrm{DH}(\mathcal{X},\mathcal{L}_{c})$
is the pushforward of Lebesgue measure by $G[\mathcal{F}_{c}]$, this
contradicts with (\ref{eq: DH of Lc}). Then combine $\inf_{\triangle}x_{1}=0$
with (1), we know $0\in\triangle(L)$.

With these facts, the claim follows from the convex-geometry description
of $J_{T}^{NA}$, see Theorem \ref{Intro convex-geometry of JNA}
(2) and refer to Figure \ref{fig: NA J}.

Hence we have
\[
J_{T}^{NA}(\mathcal{X},\mathcal{L}_{c})=J^{NA}(\mathcal{X},\mathcal{L}_{c})=\frac{c^{n+1}}{(n+1)L^{n}},\ \textrm{when}\ c\ll1.
\]
Combining this with the expansion of Berman-Ding invariant (\ref{eq: expansion of Ding invari}),
thus the uniform stability condition (\ref{eq: def of unif Dstable})
implies $\vartheta(M)\leq1-\delta$.
\end{proof}

\textbf{Acknowledgments:} The author would like to thank Feng Wang, Jia-xiang Wang and Naoto Yotsutani for helpful discussions, and specially thank Mingchen Xia for valuable revision suggestions. The author is supported by Grants: National Natural Science Foundation of China (NSFC) (No. 751203123) and Fundamental Research Funds for the Central Universities (No. 531118010149).


\begin{thebibliography}{42}
\bibitem{Atiyah}Atiyah, M.: Convexity and commuting Hamiltonians.
Bull. Lond. Math. Soc. \textbf{14}(1), 1-15 (1982)

\bibitem{Apostolov}Apostolov, V., Calderbank, D., Gauduchon, P.,
T\textrm{\o}nnesen-Friedman, C.: Hamiltonian 2-forms in K\"{a}hler
geometry, III Extremal metrics and stability. Invent. Math. \textbf{173}(3),
547-601 (2008)

\bibitem{Berman}Berman, R.: K-polystability of Q-Fano varieties admitting
K\"{a}hler-Einstein metrics. Invent. Math. \textbf{203}(3), 973-1025 (2016)

\bibitem{BBJ}Berman, R., Boucksom, S., Jonsson, M.: A variational
approach to the Yau-Tian-Donaldson conjecture. To appear in J. Amer.
Math. Soc. arXiv:1509.04561v3

\bibitem{Berman-Nystrom}Berman, R., Witt Nystr\"{o}m, D.: Complex optimal
transport and the pluripotential theory of K\"{a}hler-Ricci solitons.
arXiv: 1401.8264

\bibitem{Chen huayi}Boucksom, S., Chen, H.: Okounkov bodies of filtered
linear series. Compos. Math. \textbf{147}(4), 1205-1229 (2011)

\bibitem{BHJ1}Boucksom, S., Hisamoto, T., Jonsson. M.: Uniform K-stability,
Duistermaat-Heckman measures and singularities of pairs. Ann. Inst.
Fourier (Grenoble). \textbf{67}(2), 743-841 (2017)

\bibitem{BHJ2}Boucksom, S., Hisamoto, T., Jonsson. M.: Uniform K-stability
and asymptotics of energy functionals in K\"{a}hler geometry. J. Eur.
Math. Soc. \textbf{21}(9), 2905-2944 (2019)

\bibitem{Chen-Cheng}Chen, X., Cheng, J.: On the constant scalar curvature
K\"{a}hler metrics, general automorphism group. arXiv:1801.05907

\bibitem{Chu C11}Chu, J., Tosatti, V., Weinkove, B.: C1,1 regularity
for degenerate complex Monge-Amp\`{e}re equations and geodesic rays. Comm.
Partial Differential Equations. \textbf{43}(2), 292-312 (2018)

\bibitem{inv MA}Collins, T., Hisamoto, T., Takahashi, R.: The inverse
Monge-Amp\`{e}re flow and applications to K\"{a}hler-Einstein metrics. To
appear in J. Differential Geom. arXiv:1712.01685v2

\bibitem{Dervan relative}Dervan, R.: Relative K-stability for K\"{a}hler
manifolds. Math. Ann. \textbf{372}(3-4), 859-889 (2018)

\bibitem{Dona lower}Donaldson, S.: Lower bounds on the Calabi functional.
J. Differential Geom. \textbf{70}(3), 453-472 (2005)

\bibitem{Dona Ding}Donaldson, S.: The Ding functional, Berndtsson
convexity and moment maps. Geometry, Analysis and Probability (In
Honor of Jean-Michel Bismut). Birkh\"{a}user, Cham (2017)

\bibitem{Futaki-Mabuchi}Futaki, A., Mabuchi, T.: Bilinear forms and
extremal K\"{a}hler vector fields associated with K\"{a}hler classes. Math.
Ann. \textbf{301}(1), 199-210 (1995)

\bibitem{Li-HanJY}Han, J., Li, C.: On the Yau-Tian-Donaldson conjecture
for generalized K\"{a}hler-Ricci soliton equations. arXiv:2006.00903v3

\bibitem{Harder}Harder, G.: Lectures on algebraic geometry I. sheaves,
cohomology of sheaves, and applications to Riemann surfaces. Springer,
Wiesbaden (2011)

\bibitem{Hartshorne}Hartshorne, R.: Algebraic geometry. Springer,
New York (1977)

\bibitem{HeWY}He, W.: On Calabi's extremal metric and properness.
Trans. Am. Math. Soc. \textbf{372}(8), 5595-5619 (2019)

\bibitem{Hisa measure}Hisamoto, T.: On the limit of spectral measures
associated to a test configuration of a polarized K\"{a}hler manifold.
J. Reine Angew. Math. \textbf{713}, 129-148 (2016)

\bibitem{Hisa ortho proj}Hisamoto, T.: Orthogonal projection of a
test configuration to vector fields. arXiv: 1610.07158v3

\bibitem{Hisa toric}Hisamoto, T.: Stability and coercivity for toric
polarizations. arXiv: 1610.07998v3

\bibitem{Hisa MS}Hisamoto, T.: Mabuchi's soliton metric and relative
D-stability. arXiv: 1905.05948v2

\bibitem{Huybrechts}Huybrechts, D.: Complex geometry: an introduction.
Springer, Berlin (2005)

\bibitem{Lazarsfeld}Lazarsfeld, R.: Positivity in algebraic geometry
I: Classical setting: line bundles and linear series. Springer, Berlin
(2004)

\bibitem{L-Mustata}Lazarsfeld, R., Musta\c{t}\u{a},
M.: Convex bodies associated to linear series. Ann. Sci. \'{E}c. Norm.
Sup\'{e}r. \textbf{42}(5), 783-835 (2009)

\bibitem{Chi Li G-uniform}Li, C.: G-uniform stability and K\"{a}hler-Einstein
metrics on Fano varieties. To appear in Invent. Math. https://doi.org/10.1007/s00222-021-01075-9
(2021)

\bibitem{Chi Li cscK}Li, C.: Geodesic rays and stability in the cscK
problem. arXiv: 2001.01366

\bibitem{LiY ZhouB}Li, Y., Zhou, B.: Mabuchi metrics and properness
of the modified Ding functional. Pacific J. Math. \textbf{302}(2),
659-692 (2019)

\bibitem{Mabuchi}Mabuchi, T.: K\"{a}hler-Einstein metrics for manifolds
with nonvanishing Futaki character. Tohoku Math. J. \textbf{53}(2),
171-182 (2001)

\bibitem{Mabuchi-CM}Mabuchi, T.: A theorem of Calabi--Matsushima\textquoteright s
type. Osaka J. Math. \textbf{39}, 49-57 (2002)

\bibitem{Mabuchi-multiplier}Mabuchi, T.: Multiplier Hermitian structures
on K\"{a}hler manifolds. Nagoya Math. J. \textbf{170}, 73-115 (2003)

\bibitem{equi R-R}Meinrenken, E.: On Riemann-Roch Formulas for Multiplicities.
J. Amer. Math. Soc. \textbf{9}(2), 373-389 (1996)

\bibitem{Yotsutani}Nitta,Y., Saito S., Yotsutani, N.: Relative Ding
stability of toric Fano manifolds in low dimensions. arXiv: 1712.01131v3

\bibitem{Okounkov}Okounkov, A.: Brunn-Minkowski inequality for multiplicities.
Invent. Math. \textbf{125}(3), 405-411 (1996)

\bibitem{Phong-Sturm}Phong, D., Sturm, J.: Test configurations for
K-stability and geodesic rays. J. Symplectic Geom. \textbf{5}(2),
221-247 (2007)

\bibitem{Ross-Thomas}Ross, J., Thomas, R.: A study of the Hilbert-Mumford
criterion for the stability of projective varieties. J. Algebraic
Geom. \textbf{16}(2), 201-255 (2007)

\bibitem{Dyrefelt}Sj\"{o}str\"{o}m Dyrefelt, Z.: K-semistability of cscK
manifolds with transcendental cohomology class. J. Geom. Anal. \textbf{28}(4),
2927-2960 (2018)

\bibitem{Gabor thesis}Sz\'{e}kelyhidi, G.: Extremal metrics and K-stability.
Bull. Lond. Math. Soc. \textbf{39}(1), 76-84 (2007)

\bibitem{Wang Zhu}Wang, X., Zhu, X.: K\"{a}hler-Ricci solitons on toric
manifolds with positive first Chern class, Adv. Math. \textbf{188}(1),
87-103 (2004)

\bibitem{Nystrom TC-body}Witt Nystr\"{o}m, D.: Test configurations and
Okounkov bodies. Compos. Math. \textbf{148}(6), 1736-1756 (2012)

\bibitem{Yao}Yao, Y.: Mabuchi solitons and relative Ding stability
of toric Fano varieties. To appear in Int. Math. Res. Not. IMRN. https://doi.org/10.1093/imrn/rnab226
(2021)
\end{thebibliography}
\end{document}